\documentclass{amsart}
\usepackage{amssymb}
\usepackage{amscd}
\usepackage{verbatim}
\usepackage{epsfig}
\input xy
\xyoption{all}
\begin{document}
\newcommand\Mand{\ \text{and}\ }
\newcommand\Mwith{\ \text{with}\ }
\newcommand\Mst{\ \text{such that}\ }
\newcommand\Mor{\ \text{or}\ }
\newcommand\Mif{\ \text{if}\ }
\newcommand\Miff{\ \text{iff}\ }
\newcommand\nin{\notin}
\newcommand\identity{\operatorname{id}}
\newcommand\Id{\operatorname{Id}}
\newcommand\Real{\mathbb{R}}
\newcommand\pos{\Real^+}
\newcommand\Rnp{\Real\setminus\{0\}}
\newcommand\nzero{\setminus\{0\}}
\newcommand\Cx{\mathbb{C}}
\newcommand\Cxp{\Cx^+}
\newcommand\Cxm{\Cx^-}
\newcommand\Nat{\mathbb{N}}
\newcommand\halfNat{{\frac{1}{2}}\mathbb{N}}
\newcommand\intgr{\mathbb{Z}}
\newcommand\im{\operatorname{Im}}
\newcommand\re{\operatorname{Re}}
\newcommand\sign{\operatorname{sign}}
\newcommand\codim{\operatorname{codim}}
\newcommand\End{\operatorname{End}}
\newcommand\Ker{\operatorname{Ker}}
\newcommand\Hom{\operatorname{Hom}}
\newcommand\ideal{{\mathcal I}}
\newcommand\Span{\operatorname{span}}
\newcommand\image{\operatorname{image}}
\newcommand\Range{\operatorname{Ran}}
\newcommand\Graph{\operatorname{graph}}
\newcommand\slim{\operatornamewithlimits{s-lim}}
\newcommand\Rn{\Real^n}
\newcommand\Rm{\Real^m}
\newcommand\RN{\Real^N}
\newcommand\RtN{\Real^{2N}}
\newcommand\RM{\Real^M}
\newcommand\sphere{\mathbb{S}}
\newcommand\Sn{\sphere^{n-1}}
\newcommand\Sm{\sphere^{m-1}}
\newcommand\Snp{\sphere^n_+}
\newcommand\Smp{\sphere^m_+}
\newcommand\SN{\sphere^{N-1}}
\newcommand\SNp{\sphere^N_+}
\newcommand\circlep{\sphere^1_+}
\newcommand\Phom{P_{h}}
\newcommand\Shom{S_{h}}
\newcommand\distance{\operatorname{dist}}
\newcommand\cl{\operatorname{cl}}
\newcommand\interior{\operatorname{int}}
\newcommand\Fa{\operatorname{Fa}}
\newcommand\ff{\operatorname{ff}}
\newcommand\mf{\operatorname{mf}}
\newcommand\cf{\operatorname{cf}}
\newcommand\scf{\operatorname{sf}}
\newcommand\lf{\operatorname{lf}}
\newcommand\rf{\operatorname{rf}}
\newcommand\indfam{{\mathcal K}}
\newcommand\calA{{\mathcal A}}
\newcommand\calB{{\mathcal B}}
\newcommand\calR{{\mathcal R}}
\newcommand\calO{{\mathcal O}}
\newcommand\calX{{\mathcal X}}
\newcommand\calF{{\mathcal F}}
\newcommand\calG{{\mathcal G}}
\newcommand\calT{{\mathcal T}}
\newcommand\calC{{\mathcal C}}
\newcommand\calCt{{\tilde \mathcal C}}
\newcommand\calCL{{\mathcal C}_{\text L}}
\newcommand\calCR{{\mathcal C}_{\text R}}
\newcommand\Cinf{{\mathcal C}^{\infty}}
\newcommand\dist{{\mathcal C}^{-\infty}}
\newcommand\dCinf{\dot\Cinf}
\newcommand\ddist{\dot\dist}
\newcommand\Cj{{\mathcal C}^j}
\newcommand\Linf{L^{\infty}}
\newcommand\phg{{\text{phg}}}
\newcommand\bcon{{\mathcal A}}
\newcommand\bconc{{\mathcal A}_{\text{phg}}}
\newcommand\Sch{{\mathcal S}}
\newcommand\temp{\Sch^{\prime}}
\newcommand\Diff{\operatorname{Diff}}
\newcommand\Diffb{\operatorname{Diff}_{\text{b}}}
\newcommand\Diffc{\operatorname{Diff}_{\text{c}}}
\newcommand\Diffsc{\operatorname{Diff}_{\text{sc}}}
\newcommand\DiffI{\operatorname{Diff}_{\text{I}}}
\newcommand\DiffIq{\operatorname{Diff}_{\text{I},q}}
\newcommand\sing{\text{sing}}
\newcommand\supp{\operatorname{supp}}
\newcommand\ssupp{\operatorname{sing\ supp}}
\newcommand\csupp{\operatorname{cone\ supp}}
\newcommand\esupp{\operatorname{ess\ supp}}
\newcommand\Fr{{\mathcal F}}
\newcommand\Frinv{\Fr^{-1}}
\newcommand\bop{{\mathcal B}}
\newcommand\spec{\operatorname{spec}}
\newcommand\pspec{\spec_{pp}}
\newcommand\cspec{\spec_{c}}
\newcommand\FIO{{\mathcal I}}
\newcommand\SP{\operatorname{RC}}
\newcommand\RC{\operatorname{RC}}
\newcommand\Symc{S_c}
\newcommand\Symca{S_c^{\alpha}}
\newcommand\Symczero{S_c^{0,...,0}}
\newcommand\sci{{}^{\text{sc}}}
\newcommand\sct{\sci T^*}
\newcommand\scdt{\sci \dot T^*}
\newcommand\dS{\dot S^*}
\newcommand\dT{\dot T^*}
\newcommand\dSreg{\dot\Sigma_{\text reg}}
\newcommand\scct{\sci\bar{T}^*}
\newcommand\Csc{C_{\text{sc}}}
\newcommand\SNpscd{(\SNp)^2_{\text{sc}}}
\newcommand\scdiag{\Delta_{\text{sc}}}
\newcommand\projscl{\pi^L_{\text{sc}}}
\newcommand\projscr{\pi^R_{\text{sc}}}
\newcommand\scHL{\sci H^{2,0}_{|\zeta|^2-\lambda^2}}
\newcommand\scHrg{\sci H^{2,0}_{\sqrt{g}}}
\newcommand\Hsc{H_{\text{sc}}}
\newcommand\WF{\operatorname{WF}}
\newcommand\WFp{\operatorname{WF^{\prime}}}
\newcommand\WFsc{\operatorname{WF}_{\text{sc}}}
\newcommand\WFscp{\operatorname{WF_{sc}^{\prime}}}
\newcommand\WFC{\operatorname{WF}_C}
\newcommand\WFCi{\operatorname{WF}_{C_i}}
\newcommand\elliptic{\operatorname{ell}}
\newcommand\Psop{\operatorname{\Psi}}
\newcommand\Psiscrs{\operatorname{\Psi_{sc}^{-2,\infty}}}
\newcommand\Psiscr{\operatorname{\Psi_{sc}^{-2,0}}}
\newcommand\Psiscrm{\operatorname{\Psi_{sc}^{0,2}}}
\newcommand\PsiscHam{\operatorname{\Psi_{sc}^{2,0}}}
\newcommand\Psisci{\operatorname{\Psi_{sc}^{*,*}}}
\newcommand\Psiscid{\operatorname{\Psi_{sc}^{0,0}}}
\newcommand\Psiscis{\operatorname{\Psi_{sc}^{0,\infty}}}
\newcommand\Psiscsi{\operatorname{\Psi_{sc}^{-\infty,0}}}
\newcommand\Psiscs{\operatorname{\Psi_{sc}^{-\infty,\infty}}}
\newcommand\Psiscalg{\operatorname{\Psi_{sc}^{\infty,-\infty}}}
\newcommand\nullHam{{\mathcal N}}
\newcommand\charD{\Sigma_{\Delta-\lambda^2}}
\newcommand\charLap{\Sigma_{\Delta-\lambda}}
\newcommand\Snl{\Sn_{\lambda}}
\newcommand\SNl{\SN_{\lambda}}
\newcommand\gammat{\tilde\gamma}
\newcommand\gammasc{\gamma}
\newcommand\Tau{\mathcal{T}}
\newcommand\taut{\tilde\tau}
\newcommand\taub{\bar\tau}
\newcommand\Nout{N^+_{\lambda}}
\newcommand\Nin{N^-_{\lambda}}
\newcommand\Nio{N^{\pm}_{\lambda}}
\newcommand\El{E_{\lambda}}
\newcommand\Elt{\tilde E_{\lambda}}
\newcommand\Eil{E^i_{\lambda}}
\newcommand\Ejl{E^j_{\lambda}}
\newcommand\Eajl{E^{\alpha_j}_{\lambda}}
\newcommand\Eilt{\tilde E^i_{\lambda}}
\newcommand\Np{N^+}
\newcommand\Nm{N^-}
\newcommand\Npm{N^{\pm}}
\newcommand\Fin{F^-(\lambda)}
\newcommand\Fini{F^-_i(\lambda)}
\newcommand\Fout{F^+(\lambda)}
\newcommand\Fouti{F^+_i(\lambda)}
\newcommand\Foutj{F^+_j(\lambda)}
\newcommand\Rout{R^+_{\lambda}}
\newcommand\Routl{R^+_{\lambda^2}}
\newcommand\Routsgnl{R^{\sign\lambda}_{\lambda^2}}
\newcommand\Rin{R^-_{\lambda}}
\newcommand\Rinl{R^-_{\lambda^2}}
\newcommand\Rinsgnl{R^{-\sign\lambda}_{\lambda^2}}
\newcommand\Rio{R^{\pm}_{\lambda}}
\newcommand\Riol{R^{\pm}_{\lambda^2}}
\newcommand\Roi{R^{\mp}_{\lambda}}
\newcommand\Roil{R^{\mp}_{\lambda^2}}
\newcommand\Riob{R^{\pm}}
\newcommand\Roib{R^{\mp}}
\newcommand\Tio{T^{\pm}}
\newcommand\Tiob{T^{\pm}_{\ff}}
\newcommand\Toi{T^{\mp}}
\newcommand\Toib{T^{\mp}_{\ff}}
\newcommand\TIiob{T_I^{\pm}}
\newcommand\Rinb{R^-}
\newcommand\Rinbsgnl{R^{-\sign\lambda}}
\newcommand\Tin{T^-}
\newcommand\Tinb{T^-_{\ff}}
\newcommand\TIinb{T^-_I}
\newcommand\Routb{R^+}
\newcommand\Routbsgnl{R^{\sign\lambda}}
\newcommand\Tout{T^+}
\newcommand\Toutb{T^+_{\ff}}
\newcommand\TIoutb{T^+_I}
\newcommand\Rlkf{(|\xib|^2-(\lambda-i0)^2)^{-1}}
\newcommand\Rlk{\rho_0(\lambda)}
\newcommand\Rmlk{\rho_0(-\lambda)}
\newcommand\Rpmlk{\rho_0(\pm\lambda)}
\newcommand\Rlka{\rho_1(\lambda)}
\newcommand\Rlkb{\rho_2(\lambda)}
\newcommand\Rilk{\rho_i(\lambda)}
\newcommand\reduced{\natural}
\newcommand\Rlf{R_0(\lambda)}
\newcommand\Rla{R_1(\lambda)}
\newcommand\Rlb{R_2(\lambda)}
\newcommand\Ril{R_i(\lambda)}
\newcommand\Rlj{R_j(\lambda)}
\newcommand\Rlft{R_0(\lambda)}
\newcommand\Rflambda{R_0^{\reduced}(\sigma)}
\newcommand\RV{R^{\reduced}_V}
\newcommand\Rfsigma{R_0^{\reduced}(\sigma)}
\newcommand\Rfsigmah{R_0^{\reduced}(\sigma^{1/2})}
\newcommand\Rfzero{R_0^{\reduced}(0)}
\newcommand\RlV{R^{\reduced}_V(\sigma)}
\newcommand\RlVi{R^{\reduced}_{V_i}(\sigma)}
\newcommand\RlVt{R_V(\lambda)}
\newcommand\RlVtL{{R}_V^L(\lambda)}
\newcommand\RlVtR{{R}_V^R(\lambda)}
\newcommand\RlVit{{R}_{V_i}(\lambda)}
\newcommand\RlVta{{R}_V^{(1)}(\lambda)}
\newcommand\RlVtk{{R}_V^{(k)}(\lambda)}
\newcommand\RlVatV{{R}_{V_{\alpha}}(\lambda)V_{\alpha}}
\newcommand\RlVatVa{{R}_{V_{\alpha_1}}(\lambda)V_{\alpha_1}}
\newcommand\RlVatVb{{R}_{V_{\alpha_2}}(\lambda)V_{\alpha_2}}
\newcommand\RlVatVk{{R}_{V_{\alpha_k}}(\lambda)V_{\alpha_k}}
\newcommand\RlVatVkk{{R}_{V_{\alpha_{k+1}}}(\lambda)V_{\alpha_{k+1}}}
\newcommand\RlVaptV{{R}_{V_{\alpha'}}(\lambda)V_{\alpha'}}
\newcommand\RlVapptV{{R}_{V_{\alpha''}}(\lambda)V_{\alpha''}}
\newcommand\RlVajtV{{R}_{V_{\alpha_j}}(\lambda)V_{\alpha_j}}
\newcommand\RlVaktV{{R}_{V_{\alpha_k}}(\lambda)V_{\alpha_k}}
\newcommand\RlVakktV{{R}_{V_{\alpha_{k+1}}}(\lambda)V_{\alpha_{k+1}}}
\newcommand\Tl{T(\lambda)}
\newcommand\Tlt{\tilde\Tl}
\newcommand\Tltp{\tilde T'(\lambda)}
\newcommand\Tltpp{\tilde T''(\lambda)}
\newcommand\Tli{T_i(\lambda)}
\newcommand\Tlit{\tilde\Tli}
\newcommand\Tlip{T_i'(\lambda)}
\newcommand\Tlipp{T_i''(\lambda)}
\newcommand\Tlj{T_j(\lambda)}
\newcommand\Tla{T_{\alpha}(\lambda)}
\newcommand\Tlaa{T_{\alpha_1}(\lambda)}
\newcommand\Tlab{T_{\alpha_2}(\lambda)}
\newcommand\Tlak{T_{\alpha_k}(\lambda)}
\newcommand\Tlakt{\tilde\Tlak}
\newcommand\Tlaj{T_{\alpha_j}(\lambda)}
\newcommand\Tlajj{T_{\alpha_{j+1}}(\lambda)}
\newcommand\Tlajp{T_{\alpha_j}'(\lambda)}
\newcommand\Tlajpt{\tilde\Tlajp}
\newcommand\Tlajt{\tilde\Tlaj}
\newcommand\Tlakk{T_{\alpha_{k+1}}(\lambda)}
\newcommand\Tlakkp{T_{\alpha_{k+1}}'(\lambda)}
\newcommand\Tlap{T_{\alpha'}(\lambda)}
\newcommand\Tlapt{\tilde\Tlap}
\newcommand\Tlapp{T_{\alpha''}(\lambda)}
\newcommand\Tkl{T^{(k)}(\lambda)}
\newcommand\Tcl{T^{\flat}(\lambda)}
\newcommand\Fl{F(\lambda)}
\newcommand\BlVt{\tilde B_V(\lambda)}
\newcommand\KBlVt{K_{\BlVt}}
\newcommand\BlVaat{B_{V_{\alpha_1}}(\lambda)}
\newcommand\BV{B_V}
\newcommand\Bone{B_1}
\newcommand\Btwo{B_2}
\newcommand\Bthree{B_3}
\newcommand\Banyj{B_j}
\newcommand\PlV{P_V(\lambda)}
\newcommand\PlVc{P_V^{\flat}(\lambda)}
\newcommand\Pl{P_0(\lambda)}
\newcommand\SVl{S_V(\lambda)}
\newcommand\Sjr{S_j^{\reduced}}
\newcommand\Rkp{{\mathcal R}^k_+}
\newcommand\Rkm{{\mathcal R}^k_-}
\newcommand\Rkpm{{\mathcal R}^k_{\pm}}
\newcommand\Phys{{\mathcal P}}
\newcommand\Pc{\overline{\mathcal P}}
\newcommand\pip{\pi^{\perp}}
\newcommand\pipa{\pi_\partial}
\newcommand\gammapa{\gamma_\partial}
\newcommand\pipah{\hat\pi_\partial}
\newcommand\pit{\tilde\pi}
\newcommand\xit{\tilde\xi}
\newcommand\zetat{\tilde\zeta}
\newcommand\etat{\tilde\eta}
\newcommand\sigmat{\tilde\sigma}
\newcommand\sigmahat{\hat\sigma}
\newcommand\thetat{\tilde\theta}
\newcommand\psit{\tilde\psi}
\newcommand\phit{\tilde\phi}
\newcommand\chit{\tilde\chi}
\newcommand\rhot{\tilde\rho}
\newcommand\xib{\bar\xi}
\newcommand\zetab{\bar\zeta}
\newcommand\thetab{\bar\theta}
\newcommand\etab{\bar\eta}
\newcommand\iotal{\iota_{\lambda}}
\newcommand\rhoat{\rhot_{\alpha_1}}
\newcommand\Lambdat{\tilde\Lambda}
\newcommand\poles{\Lambda'}
\newcommand\rpoles{\Lambda_p}
\newcommand\thresholds{\Lambda}
\newcommand\Vt{\tilde V}
\newcommand\It{\tilde I}
\newcommand\half{{\frac{1}{2}}}
\newcommand\sigmah{\sigma^{1/2}}
\newcommand\bX{\partial X}
\newcommand\Deltabt{\tilde\Delta_0}
\newcommand\strip{\Omega_T}
\newcommand\Kf{K^{\flat}}
\newcommand\Gs{G^{\sharp}}
\newcommand\Gt{\tilde G}
\newcommand\Osc{\sci\Omega}
\newcommand\OSc{{}^\Scl\Omega}
\newcommand\Osch{\sci\Omega^{\half}}
\newcommand\Oscmh{\sci\Omega^{-\half}}
\newcommand\Isc{I_{sc}}
\newcommand\Qzl{Q^0_{-\lambda}}
\newcommand\Lie{{\mathcal L}}
\newcommand\bl{{\text b}}
\newcommand\scl{{\text{sc}}}
\newcommand\sccl{{\text{scc}}}
\newcommand\Scl{{\text{Sc}}}
\newcommand\ScLl{{\text{Sc,L}}}
\newcommand\ScRl{{\text{Sc,R}}}
\newcommand\Sccl{{\text{Scc}}}
\newcommand\sus{{\text{sus}}}
\newcommand\XXb{X^2_\bl}
\newcommand\XXbt{\Xt^2_\bl}
\newcommand\XXsc{X^2_\scl}
\newcommand\XXsct{\Xt^2_\scl}
\newcommand\XXSc{X^2_\Scl}
\newcommand\XXSct{\Xt^2_\Scl}
\newcommand\XXScL{X^2_\ScLl}
\newcommand\XXScR{X^2_\ScRl}
\newcommand\MMsc{M^2_\scl}
\newcommand\Deltab{\Delta_\bl}
\newcommand\Deltasc{\Delta_\scl}
\newcommand\DeltaSc{\Delta_\Scl}
\newcommand\DeltaScL{\Delta_\ScLl}
\newcommand\DeltaScR{\Delta_\ScRl}
\newcommand\prs{\sigma}
\newcommand\Nsc{N_\scl}
\newcommand\Nscp{N_{\scl,p}}
\newcommand\Nff{N_{\ff}}
\newcommand\Nffz{N_{\ff,0}}
\newcommand\Nffzp{N_{\ff,0,p}}
\newcommand\Nffl{N_{\ff,l}}
\newcommand\Nffml{N_{\ff,-l}}
\newcommand\Nmf{N_{\mf}}
\newcommand\Nmfz{N_{\mf,0}}
\newcommand\Nmfl{N_{\mf,l}}
\newcommand\Nmfml{N_{\mf,-l}}
\newcommand\ffb{\operatorname{bf}}
\newcommand\Ffb{\operatorname{bf'}}
\newcommand\ffsc{\operatorname{sf}}
\newcommand\ffSc{\operatorname{sf_C}}
\newcommand\Ffsc{\operatorname{sf'}}
\newcommand\rff{\rho_{\ff}}
\newcommand\rmf{\rho_{\mf}}
\newcommand\rffb{\rho_{\ffb}}
\newcommand\rffsc{\rho_{\ffsc}}
\newcommand\rFfsc{\rho_{\Ffsc}}
\newcommand\rffSc{\rho_{\ffSc}}
\newcommand\rinf{\rho_{\infty}}
\newcommand\CL{C_L}
\newcommand\CR{C_R}
\newcommand\betab{\beta_\bl}
\newcommand\betasc{\beta_\scl}
\newcommand\betaSc{\beta_\Scl}
\newcommand\BetaSc{\bar\beta_\Scl}
\newcommand\betaScL{\beta_\ScLl}
\newcommand\betaScR{\beta_\ScRl}
\newcommand\ScT{{}^\Scl T^*}
\newcommand\SccT{{}^\Scl \bar T^*}
\newcommand\ScS{{}^\Scl S^*}
\newcommand\Tb{{}^\bl T}
\newcommand\Tsc{{}^\scl T}
\newcommand\TSc{{}^\Scl T}
\newcommand\CSc{C_\Scl}
\newcommand\Lambdasc{{}^\scl\Lambda}
\newcommand\XXXb{X^3_\bl}
\newcommand\XXXsc{X^3_\scl}
\newcommand\XXXSc{X^3_\Scl}
\newcommand\XXXScO{X^3_{\Scl,O}}
\newcommand\XXXScF{X^3_{\Scl,F}}
\newcommand\XXXScS{X^3_{\Scl,S}}
\newcommand\XXXScC{X^3_{\Scl,C}}
\newcommand\KDsc{\operatorname{KD^{\half}_\scl}}
\newcommand\KDSc{\operatorname{KD^{\half}_\Scl}}
\newcommand\KDScEF{\operatorname{KD^{E,F}_\Scl}}
\newcommand\Oh{\operatorname{\Omega^{\half}}}
\newcommand\WFSc{\WF_\Scl}
\newcommand\WFtSc{\WF_{\text 3sc}}
\newcommand\WFScmf{\WF_{\Scl,\mf}}
\newcommand\WFScff{\WF_{\Scl,\ff}}
\newcommand\WFScs{\WF_{\Scl,\prs}}
\newcommand\WFScp{\WF'_\Scl}
\newcommand\WFScmfp{\WF'_{\Scl,\mf}}
\newcommand\WFScffp{\WF'_{\Scl,\ff}}
\newcommand\WFScsp{\WF'_{\Scl,\prs}}
\newcommand\Diffscc{\Diff_\sccl}
\newcommand\DiffSc{\Diff_\Scl}
\newcommand\DiffScc{\Diff_\Sccl}
\newcommand\DiffscI{\Diff_{\scl,\text{I}}}
\newcommand\VscI{\Vf_{\scl,\text{I}}}
\newcommand\DiffsV{\operatorname{Diff}_{\sus(V)}}
\newcommand\DiffsVsc{\operatorname{Diff}_{\sus(V),\scl}}
\newcommand\DiffsVCsc{\operatorname{Diff}_{\sus(V)-C,\scl}}   
\newcommand\Psisc{\Psop_\scl}
\newcommand\Psiscc{\Psop_\sccl}
\newcommand\PsiSc{\Psop_\Scl}
\newcommand\PsiScph{\Psop_{\Scl,\phi}}
\newcommand\PsiScra{\Psop_{\Scl,\rho^\sharp_a}}
\newcommand\PsiScc{\Psop_\Sccl}
\newcommand\PsiSccml{\Psop^{m,l}_\Sccl}
\newcommand\PsiScxx{\Psop^{*,*}_\Scl}
\newcommand\PsiScml{\Psop^{m,l}_\Scl}
\newcommand\PsiScmz{\Psop^{m,0}_\Scl}
\newcommand\PsiScmmz{\Psop^{-m,0}_\Scl}
\newcommand\PsiSckz{\Psop^{k,0}_\Scl}
\newcommand\PsiScmmml{\Psop^{-m,-l}_\Scl}
\newcommand\Psiscmkk{\Psop^{-k,k}_\scl}
\newcommand\Psiscmmmkk{\Psop^{-m-k,k}_\scl}
\newcommand\Psiscmoo{\Psop^{-1,1}_\scl}
\newcommand\Psiscmz{\Psop^{m,0}_\scl}
\newcommand\Psiscmmz{\Psop^{-m,0}_\scl}
\newcommand\PsiSckmkl{\Psop^{km,kl}_\Scl}
\newcommand\PsiScmplp{\Psop^{m',l'}_\Scl}
\newcommand\PsiScmmpllp{\Psop^{m+m',l+l'}_\Scl}
\newcommand\Psiscml{\Psop^{m,l}_\scl}
\newcommand\PsiScid{\Psop^{0,0}_\Scl}
\newcommand\PsiSczo{\Psop^{0,1}_\Scl}
\newcommand\PsiScmii{\Psop^{-\infty,\infty}_\Scl}
\newcommand\PsiScmiz{\Psop^{-\infty,0}_\Scl}
\newcommand\PsiScmoo{\Psop^{-1,1}_\Scl}
\newcommand\PsisCid{\Psop^{0,0}_{\scl-C}}
\newcommand\PsisC{\Psop_{\scl-C}}
\newcommand\Psiinf{\Psop_{\infty}}
\newcommand\Psiinfid{\Psop_{\infty}^0}
\newcommand\PsiFinf{\Psop_{\infty-\Fr}}
\newcommand\PsisVscml{\Psop^{m,l}_{\sus(V),\scl}}
\newcommand\PsisVsc{\Psop_{\sus(V),\scl}}
\newcommand\PsisVpsc{\Psop_{\sus(V_p),\scl}}
\newcommand\PsisVCSc{\Psop_{\sus(V)-C,\scl}}
\newcommand\SFinf{S_{\infty-\Fr}}
\newcommand\YsVC{Y^2_{\sus(V)-C,\scl}}
\newcommand\ffYsc{\ffsc_{\sus(V)}}
\newcommand\SXC{S(X;C)}
\newcommand\Ios{I_{\text{os}}}
\newcommand\pbL{\pi^2_{\bl,{\text L}}}
\newcommand\pbR{\pi^2_{\bl,{\text R}}}
\newcommand\pscL{\pi^2_{\scl,{\text L}}}
\newcommand\pscR{\pi^2_{\scl,{\text R}}}
\newcommand\PbO{\pi^3_{\bl,{\text O}}}
\newcommand\PscO{\pi^3_{\scl,{\text O}}}
\newcommand\PScO{\pi^3_{\Scl,{\text O}}}
\newcommand\PScF{\pi^3_{\Scl,{\text F}}}
\newcommand\PScC{\pi^3_{\Scl,{\text C}}}
\newcommand\PScS{\pi^3_{\Scl,{\text S}}}
\newcommand\pScL{\pi^2_{\Scl,{\text L}}}
\newcommand\pScR{\pi^2_{\Scl,{\text R}}}
\newcommand\CLF{\CL^F}
\newcommand\CLO{\CL^O}
\newcommand\CLS{\CL^S}
\newcommand\CLC{\CL^C}
\newcommand\DeltaYb{\Delta_{\bl,Y}}
\newcommand\DeltaYsc{\Delta_{\sus-\scl}}
\newcommand\Vf{{\mathcal V}}
\newcommand\Vb{{\mathcal V}_{\bl}}
\newcommand\Vsc{{\mathcal V}_{\scl}}
\newcommand\VSc{{\mathcal V}_{\Scl}}
\newcommand\VfI{\Vf_{\text{I}}}
\newcommand\VfIq{\Vf_{\text{I},q}}
\newcommand\scH{{}^\scl H}
\newcommand\scHg{\scH_g}
\newcommand\xh{\hat x}
\newcommand\Yh{\hat Y}
\newcommand\Zh{\hat Z}
\newcommand\Yb{\bar Y}
\newcommand\hb{\bar h}
\newcommand\xih{\hat\xi}
\newcommand\etah{\hat\eta}
\newcommand\muh{\hat\mu}
\newcommand\mub{\bar\mu}
\newcommand\nub{\bar\nu}
\newcommand\mubh{\widehat{\bar\mu}}
\newcommand\yb{\bar y}
\newcommand\ub{\bar u}
\newcommand\Qb{\bar Q}
\newcommand\Wbp{{\bar W}^\perp}
\newcommand\Wp{W^\perp}
\newcommand\Kt{\tilde K}
\newcommand\Wt{\tilde W}
\newcommand\Ut{\tilde U}
\newcommand\yt{\tilde y}
\newcommand\ft{\tilde f}
\newcommand\htil{\tilde h}
\newcommand\St{\tilde S}
\newcommand\Pt{\tilde P}
\newcommand\Rt{\tilde R}
\newcommand\qt{\tilde q}
\newcommand\Qt{\tilde Q}
\newcommand\Xb{\bar X}
\newcommand\lambdat{\tilde\lambda}
\newcommand\betat{\tilde\beta}
\newcommand\epst{\tilde\epsilon}
\newcommand\ep{\epsilon}
\newcommand\bt{\tilde b}
\newcommand\Xt{\tilde X}
\newcommand\At{\tilde A}
\newcommand\at{\tilde a}
\newcommand\Ct{\tilde C}
\newcommand\pih{\hat\pi}
\newcommand\Rh{\hat R}
\newcommand\Ah{\hat A}
\newcommand\Bh{\hat B}
\newcommand\Ch{\hat C}
\newcommand\Gh{\hat G}
\newcommand\Hh{\hat H}
\newcommand\Qh{\hat Q}
\newcommand\Ph{\hat P}
\newcommand\Nh{\hat N}
\newcommand\Sh{\hat S}
\newcommand\Gcal{{\mathcal G}}
\newcommand\GcalC{{\mathcal G}_C}
\newcommand\Jcal{{\mathcal J}}
\newcommand\JcalC{{\mathcal J}_C}
\setcounter{secnumdepth}{3}
\newtheorem{lemma}{Lemma}[section]
\newtheorem{prop}[lemma]{Proposition}
\newtheorem{thm}[lemma]{Theorem}
\newtheorem{cor}[lemma]{Corollary}
\newtheorem{result}[lemma]{Result}
\newtheorem*{thm*}{Theorem}
\newtheorem*{conj*}{Conjecture}
\numberwithin{equation}{section}
\theoremstyle{remark}
\newtheorem{rem}[lemma]{Remark}
\theoremstyle{definition}
\newtheorem{Def}[lemma]{Definition}
\newtheorem*{Def*}{Definition}
\def\signature#1#2{\par\noindent#1\dotfill\null\\*
{\raggedleft #2\par}}

\title[Many-body scattering]
{Propagation of singularities in many-body
scattering}
\author[Andras Vasy]{Andr\'as Vasy}
\date{March 15, 1999; references added June 24, 1999}
\thanks{{\em Address:} Department of Mathematics, University of California, 
Berkeley, CA 94720-3840. {\em E-mail:} \texttt{andras@math.berkeley.edu}.}
\begin{abstract}
In this paper we describe the propagation of singularities of
tempered distributional solutions $u\in\temp$ of $(H-\lambda)u=0$,
$\lambda>0$, where $H$ is
a many-body Hamiltonian $H=\Delta+V$, $\Delta\geq 0$, $V=\sum_a V_a$,
under the assumption
that no subsystem has a bound state and that the two-body interactions $V_a$
are real-valued
polyhomogeneous symbols of order $-1$ (e.g.\ Coulomb-type with the singularity
at the origin removed). Here the term `singularity' provides a
microlocal description of the lack of decay at infinity.
We use this result to prove that
the wave front relation of the free-to-free S-matrix (which,
under our assumptions, is all of the S-matrix)
is given by the broken geodesic
flow, broken at the `singular directions', on $\Sn$ at time $\pi$.
We also present a natural geometric generalization to asymptotically
Euclidean spaces.

\ \\

{\sc{Propagation des singularit\'es dans la probl\`eme de diffusion
\`a $N$ corps}}

\vspace{3 mm}

{\sc{R\'esum\'e.}} Dans cet article on d\'ecrit la propagation des
singularit\'es des
solutions temp\'er\'ees $u\in\temp$ de $(H-\lambda)u=0$, $\lambda>0$,
o\`u $H$ est un
Hamiltonien \`a $N$ corps $H=\Delta+V$, $\Delta\geq 0$, $V=\sum_a V_a$,
en supposant que les Hamiltoniens des
sous-syst\`emes n'ont 
pas de vecteurs propres (dans $L^2$),
et que les potentiels \`a deux corps $V_a$
sont 
des symboles
polyhomog\`enes r\'eels d'ordre $-1$ (par exemple, de type Coulomb, mais
sans la singularit\'e \`a l'origine). Ici le terme ``singularit\'e''
fournit une description microlocale
de la croissance des fonctions \`a l'infini.
On emploie ce r\'esultat
pour montrer que la relation de front d'onde de la matrice
de diffusion, $N$-amas $N$-amas (qui est la seule partie de la matrice
de diffusion sous nos hypoth\`eses), est donn\'ee par le flot
g\'eodesique bris\'e dans les ``directions singuli\`eres'',
sur $\Sn$ \`a temps $\pi$. On pr\'esente aussi une
g\'en\'eralisation g\'eometrique naturelle au cas des vari\'et\'es
asymptotiquement euclidiennes.
\end{abstract}

\maketitle

\section{Introduction}
In this paper we describe the propagation of singularities of generalized
eigenfunctions of a many-body Hamiltonian $H=\Delta+V$ under the assumption
that no subsystem has a bound state. We use this result to prove that
the wave front relation of the free-to-free S-matrix (which is the only
part of the S-matrix under our assumptions) is given by the broken geodesic
flow, broken at the `singular directions',
on $\Sn$ at distance $\pi$.
We remark that these
results have been proved in three-body scattering, without the assumption
on the absence of bound states, in \cite{Vasy:Propagation}.
Also, Bommier \cite{Bommier:Proprietes} and Skibsted \cite{Skibsted:Smoothness}
have shown that the kernels of the 2-cluster to free cluster and 2-cluster
to 2-cluster S-matrices are smooth, and previously Isozaki had showed this
in the three-body setting \cite{Isozaki:Structures}. However, as is clear
from the smoothness statement, the microlocal propagation picture that is
crucial, for instance, in the discussion of free-to-free scattering, does not
emerge in the previous examples when the initial state is a 2-cluster.

In this section we discuss
the setup in the Euclidean setting,
but in the following ones we move to a natural geometric
generalization introduced by Melrose in \cite{RBMSpec}.
Namely, suppose that $X$ is a manifold with boundary
equipped with a scattering metric $g$ and a cleanly intersecting family
$\calC$ of closed embedded submanifolds of $\bX$ with $C_0=\bX\in\calC$.
Thus, $g$ is
a Riemannian metric in $\interior(X)$ which is
of the form $g=x^{-4}\,dx^2+x^{-2}h$ near $\bX$; here $x\in\Cinf(X)$
is a boundary defining function.
We also assume that near every $p\in \bX$, $\calC$ is locally linearizable
(i.e.\ in suitable coordinates near $p$, every element of $\calC$ is linear);
this holds if every element of $\calC$ is totally geodesic with respect
to some metric (not necessarily $h$) on $\bX$.
Let $\Delta$ be the Laplacian of
$g$ and suppose that
$V\in\Cinf([X;\calC];\Real)$ vanishes at $\bX\setminus\bigcup\{C\in\calC:
\ C\neq\bX\}$, and
$H=\Delta+V$ -- we refer to Sections~\ref{sec:geometry} and
\ref{sec:Hamiltonian}
for a more detailed discussion of the geometric and analytic aspects
of the setup. We prove under the assumption that there are no bound states
for each of the subsystems (we describe the assumption more
precisely in Section~\ref{sec:Hamiltonian}, but it holds
for example if $V\geq 0$)
that singularities of solutions $u\in\dist(X)$ of $(H-\lambda)u\in\dCinf(X)$
propagate along generalized broken bicharacteristics of $\Delta$ which
are broken at $\calC$. We also
show that this implies a bound on the singularities of the kernel of the
free-to-free S-matrix. In effect, we show that many-body scattering is
in many respects a hyperbolic problem, much like the wave equation in domains
with corners, for which the propagation of analytic singularities
was proved by Lebeau~\cite{Lebeau:Propagation}. The geometrically
simpler setting, where the elements of $\calC$ (except $C_0=\bX$)
are disjoint, corresponds
to three-body scattering in the Euclidean setting, and then the analogy
is with the wave equation in smoothly
bounded domains, where the results for
$\Cinf$ singularities were proved by
Melrose and Sj\"ostrand \cite{Melrose-Sjostrand:I, Melrose-Sjostrand:II}
and Taylor \cite{Taylor:Grazing},
and for analytic singularities by Sj\"ostrand \cite{Sjostrand:Propagation-I}.

Here however we caution that another important aspect of typical many-body
systems {\em is} the presence of bound states of subsystems. While propagation
theorems indicate that geometry plays a central role in scattering, bound
states afford a similar role to spectral theory. Thus, in general, the
two interact, even changing the characteristic set of the Hamiltonian.
The generalized broken bicharacteristics are also more complicated in
this setting, and, as a quick argument shows, the `time $\pi$' part of our
result will not hold if bound states are present. In addition,
the Hamiltonian must possess additional structure (as the Euclidean ones
do) so that
propagation in bound states can be analyzed. Hence, in this paper,
it is natural to
impose our assumption that there are no bound states in the subsystems.

We now return to the Euclidean setting. Before we can state the
precise definitions, we need to introduce some basic (and mostly
standard) notation.
We consider the Euclidean space $\Rn$, and we assume that we
are given a (finite) family $\calX$ of linear subspaces
$X_a$, $a\in I$,
of $\Rn$ which is closed under intersections and includes the subspace
$X_1=\{0\}$ consisting of the origin, and the whole space $X_0=\Rn$.
Let $X^a$ be the
orthocomplement of $X_a$, and let $\pi^a$ be the orthogonal
projection to $X^a$, $\pi_a$ to $X_a$. A many-body Hamiltonian is an operator
of the form
\begin{equation}
H=\Delta+\sum_{a\in I} (\pi^a)^*V_a;
\end{equation}
here $\Delta$ is the positive Laplacian, $V_0=0$,
and the $V_a$ are real-valued
functions in an appropriate class which we take here to be
polyhomogeneous symbols of
order $-1$ on the vector space $X_a$
to simplify the problem:
\begin{equation}
V_a\in S^{-1}_\phg(X^a).
\end{equation}
In particular, smooth potentials $V_a$
which behave at infinity like the Coulomb
potential are allowed.
Since $(\pi^a)^*V_a$ is bounded and self-adjoint and $\Delta$ is self-adjoint
with domain $H^2(\Rn)$ on $L^2=L^2(\Rn)$, $H$ is also a self-adjoint operator
on $L^2$ with domain $H^2(\Rn)$.
We let $R(\lambda)=(H-\lambda)^{-1}$ for $\lambda
\in\Cx\setminus\Real$ be the resolvent of $H$.

There is a natural partial
ordering on $I$ induced by the ordering of $X^a$
by inclusion. (Though the ordering based on inclusion of the $X_a$
would be sometimes more natural, and we use that for the geometric
generalization of many-body scattering starting from the next section, here
we use the conventional ordering.) 
Let $I_1=\{1\}$ (recall that
$X_1=\{0\}$); $1$ is the maximal element of $I$. A maximal
element of $I\setminus I_1$ is called a 2-cluster; $I_2$ denotes the set of
2-clusters. In general, once $I_k$ has been defined for
$k=1,\ldots,m-1$, we let $I_m$ (the set of $m$-clusters)
be the set of maximal elements of
$I'_{m}=I\setminus\cup_{k=1}^{m-1} I_k$, if $I'_{m}$ is not empty.
If $I'_m=\{0\}$ (so $I'_{m+1}$ is empty), we call $H$ an $m$-body
Hamiltonian. For example, if $I\neq\{0,1\}$, and for all $a,b\nin\{0,1\}$
with $a\neq b$
we have
$X_a\cap X_b=\{0\}$, then $H$ is a 3-body Hamiltonian.
The $N$-cluster of an $N$-body Hamiltonian
is also called the free cluster, since it corresponds to the particles
which are asymptotically free.

It is convenient to compactify these spaces as in \cite{RBMSpec}. Thus, we let
$\Snp$ to be the radial compactification of $\Rn$ to a closed
hemisphere, i.e.\ a ball,
(using the standard map $\SP$ given here in \eqref{eq:SP-def}),
and $\Sn=\partial\Snp$. We write $w=r\omega$,
$\omega\in\Sn$, for polar coordinates on $\Rn$, and we let $x\in\Cinf(\Snp)$
be such that
$x=(\SP^{-1})^*(r^{-1})$ for $r>1$. Hence,
$x$ is a smoothed version of $r^{-1}$
(smoothed at the origin of $\Rn$), and it is a boundary defining function
of $\Snp$. We usually identify (the interior of)
$\Snp$ with $\Rn$. Thus, we write $S^{m}_\phg(\Snp)$ and $S^{m}
_\phg(\Rn)$
interchangeably and we drop the explicit pull-back notation in the future
and simply write $x=r^{-1}$ (for $r>1$). We also remark that we have
\begin{equation}
S^m_{\phg}(\Snp)=x^m\Cinf(\Snp).
\end{equation}
We recall that under $\SP$, $\dCinf(\Snp)$, the space of smooth functions
on $\Snp$ vanishing to infinite order at the boundary corresponds
to the space of Schwartz functions $\Sch(\Rn)$, and its dual, $\dist(\Snp)$,
to tempered distributions $\temp(\Rn)$.
We also have the following correspondence of weighted
Sobolev spaces
\begin{equation}\label{eq:Sob-def}
\Hsc^{k,l}(\Snp)=H^{k,l}=H^{k,l}(\Rn)=\langle w\rangle^{-l}H^k(\Rn)
\end{equation}
where $\langle w\rangle=(1+|w|^2)^{1/2}$. Thus, for $\lambda
\in\Cx\setminus\Real$ the resolvent extends to a map
\begin{equation}\label{eq:res-weight}
R(\lambda):\Hsc^{k,l}(\Snp)\to \Hsc^{k+2,l}(\Snp).
\end{equation}

Similarly, we let
\begin{equation}
\Xb_a=\cl(\SP(X_a)),\quad C_a=\Xb_a\cap\partial\Snp.
\end{equation}
Hence, $C_a$ is a sphere of dimension $n_a-1$ where
$n_a=\dim X_a$. We also let
\begin{equation}
\calC=\{C_a:\ a\in I\}.
\end{equation}
Again, we write the polar coordinates on $X_a$ (with respect
to the induced metric) as $w_a=r_a\omega_a$, $\omega_a\in C_a$, and let
$x_a=r_a^{-1}$ (for $r_a>1$).
We note that if $a$ is a 2-cluster then
$C_a\cap C_b=\emptyset$ unless $b\leq a$.
We also define the `singular part' of $C_a$ as the set
\begin{equation}
C_{a,\sing}=\cup_{b\not\leq a}(C_b\cap C_a),
\end{equation}
and its `regular part' as the set
\begin{equation}
C'_a=C_a\setminus\cup_{b\not\leq a} C_b=C_a\setminus C_{a,\sing}.
\end{equation}
For example, if $a$ is a 2-cluster then $C_{a,\sing}=\emptyset$ and
$C'_a=C_a$.
We sometimes write the
coordinates on $X_a\oplus X^a$ as $(w_a,w^a)$.

Corresponding to each cluster $a$ we introduce the cluster
Hamiltonian $h_a$ as an operator on $L^2(X^a)$ given by
\begin{equation}
h_a=\Delta+\sum_{b\leq a} V_b,
\end{equation}
$\Delta$ being the Laplacian of the induced metric on $X^a$.
Thus, if $H$ is a $N$-body Hamiltonian and $a$ is a $k$-cluster,
then $h_a$ is a $(N+1-k)$-body Hamiltonian. The $L^2$
eigenfunctions
of $h_a$ play an important role in many-body scattering; we
remark that by Froese's and Herbst's result, \cite{FroExp},
$\pspec(h_a)\subset(-\infty,0]$ (there are no positive
eigenvalues). Moreover, $\pspec(h_a)$ is bounded below since
$h_a$ differs from $\Delta$ by a bounded operator. Note that
$X^0=\{0\}$, $h_0=0$, so the unique eigenvalue of $h_0$ is $0$.

The eigenvalues of $h_a$ can be used to define the set of
thresholds of $h_b$. Namely, we let
\begin{equation}
\Lambda_a=\cup_{b<a}\pspec(h_b)
\end{equation}
be the set of thresholds of $h_a$, and we also let
\begin{equation}
\Lambda'_a=\Lambda_a\cup\pspec(h_a)
=\cup_{b\leq a}\pspec(h_b).
\end{equation}
Thus, $0\in\Lambda_a$ for $a\neq 0$ and $\Lambda_a\subset(-\infty,0]$.
It follows from the Mourre
theory (see e.g.\ \cite{FroMourre, Perry-Sigal-Simon:Spectral})
that $\Lambda_a$ is closed,
countable, and $\pspec(h_a)$ can only accumulate at $\Lambda
_a$. Moreover, $R(\lambda)$, considered as an operator on weighted
Sobolev spaces, has a limit
\begin{equation}
R(\lambda\pm i0):\Hsc^{k,l}(\Snp)\to \Hsc^{k+2,l'}(\Snp)
\end{equation}
for $l>1/2$, $l'<-1/2$,
from either half of the complex plane away from
\begin{equation}
\Lambda=\Lambda_1\cup\pspec(H).
\end{equation}
In addition, $L^2$ eigenfunctions
of $h_a$ with eigenvalues which are not thresholds are necessarily
Schwartz functions on $X^a$ (see \cite{FroExp}). We also label the
eigenvalues of $h_a$, counted with multiplicities, by integers $m$,
and we call the pairs $\alpha=(a,m)$ channels. We denote the eigenvalue
of the channel $\alpha$ by $\epsilon_\alpha$, write $\psi_\alpha$ for
a corresponding normalized eigenfunction, and let $e_\alpha$ be the
orthogonal projection to $\psi_\alpha$ in $L^2(X^a)$.

The definition of the free-to-free S-matrix we consider
comes from the stationary theory, more precisely
from the asymptotic behavior of generalized eigenfunctions, see
\cite{Vasy:Asymptotic}, and cf.\ \cite{RBMSpec, Vasy:Propagation-2}. Apart from
the difference in normalization, it is the same as the S-matrix
given by the wave operators, see \cite{Vasy:Scattering}.
For simplicity, we state the asymptotic expansion under the assumption
that $V_a$ is polyhomogeneous of order $-2$ (so it decays as $|w^a|^{-2}$).
Namely,
for $\lambda\in(0,\infty)$
and $g\in\Cinf_c(C_0')$,
there is a unique $u\in\dist(\Snp)$
(i.e.\ $u\in\temp(\Rn)$)
such that $(H-\lambda)u=0$, and $u$ has the form
\begin{equation}\label{eq:intro-21}
u=e^{-i\sqrt{\lambda} r}r^{-(n-1)/2}v_-
+R(\lambda+i0)f,
\end{equation}
where $v_-\in\Cinf(\Snp)$, $v_-|_{\Sn}=g$, and
$f\in\dCinf(\Snp)$.
In addition, this $u$ is of the form
\begin{equation}\label{eq:intro-23}
u=e^{-i\sqrt{\lambda} r}r^{-(n-1)/2}v_-+e^{i\sqrt{\lambda} r}r^{-(n-1)/2}v_+,
\quad v_+\in\Cinf(\Snp\setminus C_{0,\sing}).
\end{equation}
The Poisson operator with free initial data is the operator
\begin{equation}
P_{0,+}(\lambda):\Cinf_c(C'_0)\to\dist(\Snp),\quad P_{0,+}(\lambda)g=u.
\end{equation}
Following \cite{Vasy:Asymptotic}, we define the free-to-free scattering
matrix, $S_{00}(\lambda)$ as the map
\begin{equation}
S_{00}(\lambda):\Cinf_c(C'_0)\to\Cinf(C'_0),
\end{equation}
\begin{equation}
S_{00}(\lambda)g=v_+|_{C'_0},
\end{equation}
so it relates the incoming amplitude $v_-|_{\Sn}$ to the outgoing one,
$v_+|_{\Sn}$.
We recall from \cite{Vasy:Scattering}
that the wave operator free-to-free S-matrix is then given by
$i^{n-1}S_{00}(\lambda)R$
(as maps $\Cinf_c(C'_0)\to\dist(C'_0)$)
where $R$ is pull back by the antipodal map on $C_0$.

There are only minor changes if $V_a$ is polyhomogeneous of order $-1$.
Namely, the asymptotic expansions in \eqref{eq:intro-21} and
\eqref{eq:intro-23} must be replaced by
\begin{equation}\label{eq:intro-31}
e^{\pm i\sqrt{\lambda}r}r^{-i
\alpha_\pm-(n-1)/2}v_\pm,\quad \alpha_\pm=\alpha_{\pm,\lambda}=\pm V'|_{C'_0}/
2\sqrt{\lambda}\in\Cinf(C'_0),\ V=xV',
\end{equation}
\begin{equation}\label{eq:intro-32}
v_\pm\sim\sum_{j=0}^\infty\sum_{s\leq 2j}
a_{j,s,\pm}(\omega)r^{-j}(\log r)^s,\quad a_{j,s,-}\in\Cinf_c(C'_0),
\ a_{j,s,+}\in\Cinf(C'_0).
\end{equation}
Note that $\alpha_\pm$ are not defined at $C_{0,\sing}$, but that does not
cause any problems even in the uniqueness statement,
\eqref{eq:intro-21}, since $v_-$ vanishes at $\Sn$ near $C_{0,\sing}$ to
infinite order.

\begin{figure}[ht]
\begin{center}
\mbox{\epsfig{file=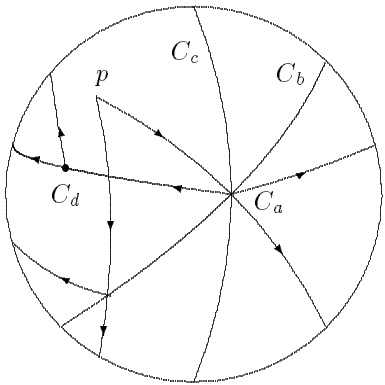}}
\end{center}
\caption{Broken geodesics on $\sphere^2$ starting at $p$. Here
$C_a=C_b\cap C_c$.}
\label{fig:broken-geodesics}
\end{figure}

Our main theorem describes the structure of $S_{00}(\lambda)$.
We first introduce
the broken geodesic flow (of the standard Riemannian metric $h$)
on $\Sn$, broken at $\calC$.
We denote by $S\Sn$ the
sphere bundle of $\Sn$ identified as the unit-length subbundle of $T\Sn$
with respect to $h$.
Let $I=[\alpha,\beta]\subset\Real$ be
an interval.
We say that a curve
$\gamma:I\to\Sn$ is a broken geodesic of $h$, broken at $\calC$,
if two conditions
are satisfied.
First, there exists a finite set of points $t_j\in I$,
$\alpha=t_0<t_1<\ldots<t_{k-1}<t_k=\beta$ such that for each $j$,
$\gamma|_{[t_j,t_{j+1}]}$ is a geodesic of $h$, and for all $t\in
(t_j,t_{j+1})$,
$\gamma'(t)\in S\Sn$.
Second, for all $j$, if $\gamma(t_j)\in C'_a$ then
the limits $\gamma'(t_j-0)$ and
$\gamma'(t_j+0)$ both exist and differ by a vector in $T_{\gamma(t_j)}\Sn$
which is orthogonal to $T_{\gamma(t_j)}C_a$ (i.e.\ the usual law of
reflection is satisfied; see Figure~\ref{fig:broken-geodesics}).
We say that $p,q\in S\Sn$ are related by the
broken geodesic flow at time $\pi$ if there is a broken geodesic $\gamma$
defined on $[0,\pi]$, such that $\gamma'(0)=p$, $\gamma'(\pi)=q$. Using
the metric $h$ to identify $S\bX$ and $S^*\bX$, this defines the broken
geodesic `flow' at time $\pi$ on $S^*\bX$. We refer to
Definition~\ref{Def:gen-br-geod} and Section~\ref{sec:tot-geod-bich}
for a more complete discussion. We then have the following result:

\begin{thm*}
Suppose that no subsystem of $H$ has bound states, i.e.\ for $a\neq 0$,
$\pspec(h_a)=\emptyset$.
Then the free-to-free scattering matrix, $S_{00}(\lambda)$,
extends to a continuous linear map
$\dist_c(C'_0)\to\dist(C'_0)$. The wave front relation of $S_{00}(\lambda)$ is
given by the broken geodesic flow at time $\pi$.
\end{thm*}

In the actual many-body problem, $w\in X_a$ means that several particles
are close to each other, namely the ones corresponding to the cluster
decomposition $a$. Thus, $\omega\in C_a$ is a statement that the particles
corresponding to cluster $a$ collide. Hence, the Theorem describes how
many-body scattering can be understood, modulo smoothing (hence in the
$\Cinf$ sense trivial) terms, as a sequence of a finite number of
collisions involving the particles. Namely, each `break' $t_j$ in the
broken geodesic describes a collision involving the cluster
decomposition $a$. In the three-body setting with Schwartz potentials
it was shown in \cite{Vasy:Structure}
that the amplitude of the reflected wave is given, to top
order, by the corresponding 2-body S-matrix; an analogous statement also
holds for short-range potentials. In particular, this shows that the Theorem
is sharp as far as the location of singularities is concerned.

We also remark that in the Euclidean setting,
unbroken geodesic flow to distance $\pi$
amounts to pull-back by the antipodal map on $\Sn=\partial\Snp$,
so it corresponds to free
propagation: particles leave in the direction opposite to the one
from which they entered.

Our approach to proving this theorem is via the analysis of generalized
eigenfunctions of $H$, i.e.\ of $u\in\dist(\Snp)$ satisfying
$(H-\lambda)u=0$. We prove that `singularities' of generalized eigenfunctions
of $H$ propagate along broken bicharacteristics in the characteristic
set of $H$, similarly to singularities of the solutions of the wave
equation. Here `singularities' are not understood as the lack of smoothness:
indeed $H$ is elliptic in the usual sense, so every generalized eigenfunction
is $\Cinf$ in the interior of $\Snp$, i.e.\ on $\Rn$. Instead, in this
situation singularity
means the lack of rapid decay $u$. Correspondingly, we define a wave
front set, $\WFSc(u)$, at infinity, i.e.\ at $\partial\Snp$, and
we will prove its invariance under the broken bicharacteristic flow.

The two notions of singularities are very closely related
via the Fourier transform. Here for simplicity consider $\Delta-\lambda$
in place of $H-\lambda$. If $(\Delta-\lambda)u=0$, then the Fourier
transform of $u$, $\Fr u$, satisfies $(|\xi|^2-\lambda)\Fr u=0$ where
$\xi$ is the dual variable of $w$. Now, the multiplication operator
$P=|\xi|^2-\lambda$ can be regarded as a 0th order differential operator.
Hence, by H\"ormander's theorem, see e.g.\ \cite{Hor}, $\WF(\Fr u)$
is invariant under the bicharacteristic flow in the characteristic
variety of $P$, i.e.\ in the set $\{(\xi,\xi^*):\ |\xi|^2-\lambda=0\}$
where we have written $\xi^*$ for the dual variable of $\xi$, so $\xi^*$ is
in fact $w$. Moreover, in the two-body problem, i.e.\ if $V$ is a symbol
(of say order $-1$) on $\Rn$, $H=\Delta+V$, and if $(H-\lambda)u=0$,
we still have $P\Fr u=0$ where now $P=|\xi|^2-\lambda+\Fr V\Frinv$.
Since $V$ is a symbol of order $-1$,
$\Fr V\Frinv$ is a pseudo-differential operator
of order $-1$, hence lower order than $|\xi|^2-\lambda$. Thus, the
principal symbol of $P$ is still $|\xi|^2-\lambda$ (recall that $\xi^*$
is the cotangent variable, so this is indeed homogeneous of order $0$
in $\xi^*$ -- it is independent of $\xi^*$). Hence, H\"ormander's
theorem is applicable and we have the same propagation statement as before.

In the two-body setting the relevant wave front set measuring lack of
decay at infinity is the scattering one, $\WFsc$. For $u\in\temp(\Rn)$,
$\WFsc(u)$ is essentially given by the usual wave front set of the Fourier
transform of $u$, i.e.\ by $\WF(\Fr u)$, after interchanging the role
of the base and dual variables. Since the Fourier transform interchanges
decay at infinity and smoothness, $\WF(\Fr u)$ indeed measures the
decay of $u$ at infinity in a microlocal sense. Hence, H\"ormander's
propagation theorem translated directly into a propagation theorem for
$\WFsc(u)$. This result
was described by Melrose in \cite{RBMSpec} where
he introduced the notion of $\WFsc$.

In the many-body setting conjugation by the Fourier transform is much less
convenient. Hence, we will design an appropriate microlocal way of measuring
the lack of decay at infinity without resorting to the Fourier transform.
Instead, we introduce an algebra of many-body pseudo-differential operators
$\PsiSc(\Snp,\calC)$ which reflects the geometry, and use it
to define the wave front
set at infinity. We then prove a propagation of singularities theorem
for generalized eigenfunctions of many-body Hamiltonians $H$; here
`singularities' are understood in the sense of the new wave front set
at infinity. The proof of this theorem is via a microlocal positive commutator
estimate, similarly to the proof of H\"ormander's theorem, or indeed
to the proof of the propagation theorems for $\Cinf$ singularities
of solutions of the wave
equation with domains with boundaries \cite{Melrose-Sjostrand:I}.
Finally, we relate such a
result to the structure of the S-matrix. This step is comparatively easy as
indicated by our description of the S-matrix in terms of generalized
eigenfunctions of $H$.

Positive commutator estimates have also played a major role in
many-body scattering
starting with the work of Mourre \cite{Mourre-Absence}, Perry, Sigal
and Simon \cite{Perry-Sigal-Simon:Spectral}, Froese and Herbst
\cite{FroMourre}, Jensen \cite{Jensen:Propagation},
G\'erard, Isozaki and Skibsted \cite{GerComm, GIS:N-body}
and Wang \cite{XPWang}.
In particular,
the Mourre estimate is one of them; it estimates $i[H,w\cdot D_w+D_w\cdot w]$.
This and some other {\em global} positive commutator results
have been used to prove the global results mentioned in the first
paragraph about some of the S-matrices
with initial state
in a two-cluster. They also give the basis for the existence, uniqueness
and equivalence statements in our definition of the S-matrix by
asymptotic expansions; these statements are discussed in
\cite{Vasy:Scattering} in more detail. Correspondingly, these
global estimates will appear in
Sections~\ref{sec:resolvent}-\ref{sec:S-matrix} of this paper
where we turn the propagation results for generalized eigenfunctions
into statements about the S-matrix.

We remark that the wave-operator approach defines the S-matrix
as a bounded
operator $L^2(C_0)\to L^2(C_0)$. Since $C_{0,\sing}$ has measure $0$,
$L^2(C_0)$ and $L^2(C'_0)$ can be identified. As $\Cinf_c(C'_0)$ is dense
in $L^2(C'_0)$, the asymptotic expansion S-matrix $S_{00}(\lambda)$ indeed
determines the wave-operator one.

The propagation of singularities of generalized eigenfunctions of $H$
is determined by the principal part of $H$; terms decaying at the
boundary do not change the analysis. As opposed to this, the precise
structure of incoming and outgoing functions, $R(\lambda\pm i0)f$,
$f\in\dCinf(\Snp)$, depends on lower order terms; a relatively trivial
example is given by the appearence of $r^{-i\alpha_\pm}$ in
\eqref{eq:intro-31} for long-range potentials.
Since we consider
$S_{00}(\lambda)$ and $P_{0,+}(\lambda)$
as operators on distributions
supported away from $C_{0,\sing}$, we do not need to analyze the
precise structure of incoming/outgoing functions at $C_{0,\sing}$, which is not
`principal type', although we certainly analyze the propagation
of singularities there. Thus, we do not discuss what happens
when the support of the incoming scattering
data increases to $C'_0$, even if the data are
$L^2$. But the behavior of $P_{0,+}(\lambda)$,
as the support of the data increases to $C'_0$, plays an important
part in asymptotic completeness, which states that all possible
outcomes of a scattering experiment are indeed described by a combination
of bound states of the cluster Hamiltonians, with asymptotically
free motion in the intercluster variables. Thus, our results
cannot be used directly to supply a
proof of asymptotic completeness. This completeness property of many-body
Hamiltonians was proved by Sigal and Soffer, Graf, Derezi\'nski and Yafaev
under different assumptions on the potentials and by different techniques
\cite{Sigal-Soffer:N, Sigal-Soffer:Long-range, Sigal-Soffer:Asymptotic,
Sigal-Soffer:Asymptotic-4, Graf:Asymptotic, Derezinski:Asymptotic, Yaf}. In
particular, Yafaev's paper \cite{Yaf}
shows quite explicitly the importance of the
special structure of the Euclidean Hamiltonian. This
structure enables him to obtain
a positive commutator estimate, which would not follow from our
indicial operator arguments in Section~\ref{sec:commutators}, and which is then
used to prove asymptotic completeness.

Finally we comment on the requirement that the collection $\calC$ be
locally linearizable. We show in the next section that
it is equivalent to the existence of a neighborhood
of every point $p\in\bX$ and a metric on it, in terms of which
all elements of $\calC$ are totally geodesic.
The importance of this assumption is closely
related to the existence of a sufficient number of {\em smooth}
vector fields on $\bX$ which are tangent to every element of $\calC$.
Such smooth vector fields always exists once we {\em resolve} the geometry
of $\calC$, i.e.\ on the blown-up space $[\bX;\calC]$, but in general,
without our assumption, there are not enough such smooth vector fields
on $\bX$. In the first part of the paper, we discuss the pseudo-differential
algebra associated to many-body scattering. For this purpose we need to
blow up $\calC$, in part for analyzing the indicial operators (see the
following paragraph). Thus, in this part of the paper, the issue of
local linearizability is irrelevant, and we do not assume it. However,
in the second part of the paper, both the discussion of generalized
broken bicharacteristics and the construction of the positive
commutators would be more complicated without it,
so from Section~\ref{sec:Hamiltonian}
on, we assume the local linearizability of $\calC$.

This paper is organized as follows. In the next section we describe the
geometric generalization of the many-body problem which was outlined
above. This includes a discussion of many-body geometry and the
definition of many-body differential operators. In Section~\ref{sec:calculus}
we proceed to define and analyze
the corresponding algebra of pseudo-differential
operators, $\PsiSc(X,\calC)$, which reflects this geometry. It
includes many-body
Hamiltonians, as well as their resolvent away from the real axis. It
extends the definition of the three-body calculus presented in
\cite{Vasy:Propagation-2}, though here we emphasize the definition
of the calculus via localization and quantization as opposed to the
conormal description of the kernels on an appropriate resolved space.
In Section~\ref{sec:indicial} we construct the indicial operators in this
calculus. They provide a non-commutative analog of the principal symbol
in standard microlocal analysis. Our proof of positivity in
commutator estimates is based on replacing the argument of
Froese and Herbst \cite{FroMourre} by indicial operator techniques.
In Section~\ref{sec:WF} we define the wave front set at infinity,
$\WFSc(u)$, corresponding to the many-body geometry and pseudo-differential
operators. The definition given here differs from the one in
\cite{Vasy:Propagation-2}; it follows the fibred cusp definition of
Mazzeo and Melrose \cite{Mazzeo-Melrose:Fibred}. These definitions, however,
give the same result for approximate generalized eigenfunctions of $H$.

In Section~\ref{sec:Hamiltonian} we discuss many-body type Hamiltonians
and their generalized broken bicharacteristics. This section is, to
a significant degree, based on Lebeau's paper \cite{Lebeau:Propagation}.
In Section~\ref{sec:tot-geod-bich} we give a much more detailed description
of the generalized broken bicharacteristics in the case when all
elements $C\in\calC$ are totally geodesic. Of course, this is true
in the Euclidean setting.
In Sections~\ref{sec:positive}-\ref{sec:commutators}
we build the technical tools for turning a symbolic positive commutator
calculation into an operator estimate.
In Section~\ref{sec:propagation} we prove
that singularities of generalized eigenfunctions of many-body type
Hamiltonians propagate along generalized broken bicharacteristics.
This is the main new result of the paper.
In Sections~\ref{sec:resolvent}-\ref{sec:S-matrix} we use this
and adaptations of the global estimates, in particular those of
G\'erard, Isozaki and Skibsted \cite{GerComm, GIS:N-body},
to analyze the structure of the resolvent and
that of the scattering matrix.
Finally, in the Appendix we prove some of the results quoted from
Lebeau's paper, using slightly different methods.

The propagation estimates of
Section~\ref{sec:propagation} lie at the heart of this paper.
The reader may want to skip some of the
technical sections when reading the paper for the first time. It may
be useful to
keep Mourre-type estimates and especially their microlocalized versions
as in \cite{GerComm, GIS:N-body} in mind while reading
Section~\ref{sec:propagation}.

I would like to thank Richard Melrose for suggesting this problem to me
(in the three-body setting) as my PhD thesis problem and for our
very fruitful discussions. His firm belief that scattering theory can
be understood in microlocal terms similar to the well-known theory
of hyperbolic operators motivated me both during my PhD work
\cite{Vasy:Propagation} and while working on its extension that appears
in this paper. I am grateful to Maciej Zworski for introducing
me to the work of Gilles Lebeau \cite{Lebeau:Propagation}, for many
helpful discussions and for his encouragement. It was Lebeau's paper
that convinced me that the results presented here were within reach,
and it plays a particularly central role in Section~\ref{sec:Hamiltonian}
where the generalized broken bicharacteristics are described.
I would also like to thank
Andrew Hassell, Rafe Mazzeo, Erik Skibsted
and Jared Wunsch for helpful discussions,
their encouragement and
for their interest in this research.

\section{Many-body geometry and differential operators}\label{sec:geometry}
It is convenient to carry out the construction in the general
geometric setting. We first describe the many-body geometry.

Thus, let $X$ be a compact manifold with boundary, and let
\begin{equation}
\calC=\{C_a:\ a\in I\}
\end{equation}
be a finite set of closed embedded submanifolds of $\bX$ such
that $\bX=C_0\in\calC$ and
for all $a,b\in I$ either $C_a$ and $C_b$
are disjoint, or they
intersect cleanly and $C_a\cap C_b=C_c$ for some
$c\in I$. We introduce a partial order on $\calC$ given by
inclusion on $\calC$, namely
\begin{equation}
C_a\leq C_b\ \text{if and only if}\ C_a\subset C_b.
\end{equation}
This partial order is the opposite of the partial order used 
traditionally in many-body scattering, discussed in the introduction,
but it will be more convenient
for us since it simply corresponds to inclusion.
A chain is defined as usual as a set on which $<$ gives a linear order.

\begin{Def}
Let $X$ and $\calC$ be as above.
We say that $(X,\calC)$ is a space with $N$-body geometry
(or an $N$-body space), $N\geq 2$,
if the maximal length of chains is $N-1$. Similarly, we say that
$C_a$ is a $k$-cluster if the maximal length of chains whose
maximal element is $C_a$ is $k-1$. We also say that $(X,\calC)$ is
a many-body space if we do not wish to specify $N$.
\end{Def}

Thus, if $C_a$ is minimal,
it is a 2-cluster, and if $(X,\calC)$ is a space with $N$-body geometry
then $\bX$ is an $N$-cluster. The numerology is chosen here so that
we conform to the usual definitions in Euclidean many-body scattering,
described in the Introduction.

Before defining the algebra
of many-body scattering differential operators on $(X,\calC)$, we
discuss the simultaneous local linearizability of the
collection $\calC$. As we have mentioned in the Introduction, the analysis
of generalized broken geodesics as well as the commutator constructions
of this paper become simpler if $\calC$ is locally linearizable.
To make this notion precise, we make the following definition.

\begin{Def}
We say that a many-body space $(X,\calC)$ is locally linearizable
(or is locally trivial) if
for every $p\in\bX$ there exists a diffeomorphism $\phi$ from a neighborhood
$U$ of $p$ in $\bX$
to a neighborhood $U'$ of the origin of a vector space $V$ such
that for each $C\in\calC$, the image of $C\cap U$ under $\phi$ is
the intersection of a linear subspace of $V$ with $U'$.
\end{Def}

\begin{rem}
In three-body type geometry, where the elements of $\calC$ except
$C_0$ are disjoint, $(X,\calC)$ is automatically locally linearizable.
The same holds, essentially by definition, if $\calC$ is a normal
collection, see \cite[Chapter~V]{RBMDiff}.
\end{rem}

Local triviality of $\calC$ is closely related to the question whether
every element of $\calC$ is locally totally geodesic with respect to some
metric. In fact,

\begin{lemma}
A many-body space $(X,\calC)$ is locally linearizable if and only if
every $p\in\bX$ has a neighborhood $U$ in $\bX$ and a Riemannian metric
$h_U$ on $U$ such that for each element $C$ of $\calC$,
$C\cap U$ is totally geodesic with respect to $h_U$.
\end{lemma}

\begin{proof}
Suppose first that $p\in\bX$ and $U$, $h_U$ are as above. By shrinking
$U$ if necessary, we can make sure that $p\nin C$ implies
$C\cap U=\emptyset$ for
every $C\in\calC$ . By shrinking $U$ further if necessary, we can arrange
that the
exponential map of $h_U$ at $p\in\bX$ identifies a neighborhood $U'$ of the
origin in $V=T_p\bX$ and $U$. Moreover,
the elements $C\in\calC$ for which $p\in C$, are
identified with $T_p C\cap U'$, since these $C$ are totally geodesic.
This proves that $(X,\calC)$ is locally linearizable.

Conversely, if $(X,\calC)$ is locally linearizable, then the choice
of an inner product on $V$ induces a metric on $TV$, hence on $U$ via the
diffeomorphism $\phi$, and as linear subspaces of $V$ are totally
geodesic with respect to this metric on $TV$, the same holds for
$\calC$ over $U$.
\end{proof}

After this brief discussion on the local
linearizability of $\calC$, we turn to the
setting of most interest, namely to Euclidean many-body geometry.
Suppose that
$X=\Snp$ is the radial compactification of $\Rn$ and $\calX$ is a
family of linear subspaces of $\Rn$ as discussed in the introduction.
Recall from
\cite{RBMSpec} that $\SP:\Rn\to\Snp$ is given by
\begin{equation}\label{eq:SP-def}
\SP(w)=(1/(1+|w|^2)^{1/2},w/(1+|w|^2)^{1/2})\in\Snp\subset\Real^{n+1},
\quad w\in\Rn.
\end{equation}
Here we use the notation $\SP$ instead of $\operatorname{SP}$, used
in \cite{RBMSpec}, to avoid confusion with the standard stereographic
projection giving a one-point compactification of $\Rn$.
We write the coordinates on $\Rn=X_a\oplus X^a$ as $(w_a,w^a)$.
Let $m=\dim X_a$.
We again let
\begin{equation}
\Xb_a=\cl(\SP(X_a)),\quad C_a=\Xb_a\cap\partial\Snp.
\end{equation}

We next show that polyhomogeneous symbols on $X^a$, pulled back to $\Rn$
by $\pi^a$, are smooth on the blown-up space $[X;C_a]$. Recall that
the blow-up process is simply an invariant way of introducing
polar coordinates about a submanifold. A full description appears
in \cite{RBMDiff} and a more concise one in \cite[Appendix~A]{RBMSpec},
but we give a brief summary here.
Thus, suppose that $X$ is a manifold
with corners and
$C$ is a p-submanifold (i.e.\ product submanifold) of $\bX$.
Thus, near any $p\in C$ we have local
coordinates $x_i$ $(i=1,\ldots,r)$, $y_j$ $(j=1,\ldots,n-r)$, $n=\dim X$,
such that the boundary hypersurfaces of $X$
through $p$ are defined by $x_i=0$, and $X$ is given by $x_i\geq 0$,
$i=1,\ldots,r$. A tangent vector $V\in T_q X$, $q$ near $p$, is inward-pointing
if $Vx_i(q)\geq 0$ for all $i$. The normal bundle of $C$ is the
quotient bundle
\begin{equation}
N C=T_C X/ T C.
\end{equation}
The inward pointing normal bundle of $C$, $N^+ C$, is the image of $T^+ X$,
consisting of inward pointing tangent vectors, in $N C$. Thus,
near $p$, $X$ is diffeomorphic
to the inward-pointing normal bundle of $C$. The blow-up of $X$ along $C$
is locally defined as the blow up of the $0$ section of $N^+ C$, i.e.\ by
introducing the new $\Cinf$ structure in $N^+ C$ given by polar coordinates
in the fibers of the bundle and by the base coordinates pulled back
from $C$. While this construction depends on some choices, the resulting
$\Cinf$ structure does not. The blow-up of $X$ along $C$ is denoted
by $[X;C]$. The blow-down map $[X;C]\to X$ is the smooth map corresponding
to expressing standard coordinates on a vector space, $N_q^+C$, in terms
of polar coordinates. It is denoted by $\beta[X;C]$. The front face
of the blow-up is the inverse image of $C$ (i.e.\ of the zero section
of $N^+ C$) under $\beta[X;C]$. Hence, it is a bundle over $C$
whose fibers are the intersection of a
sphere with a `quadrant' corresponding to the inward-pointing
condition, i.e.\ to $x_i\geq 0$. In fact, it is the inward pointing
sphere bundle $S^+ NC$ which is the quotient of $N^+ C\setminus o$,
$o$ denoting the zero section, by the natural
$\Real^+$ actions in its fibers.

\begin{figure}[ht]
\begin{center}
\mbox{\epsfig{file=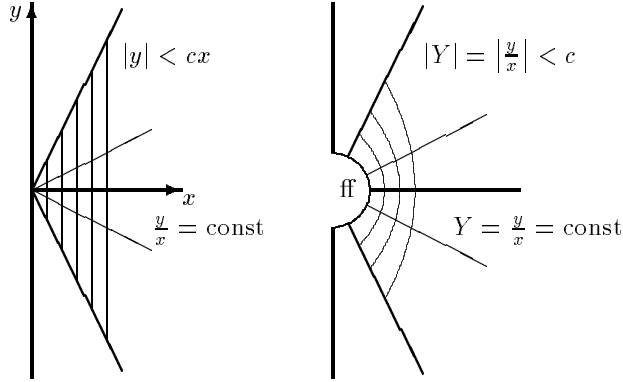}}
\end{center}
\caption{The blowup of $C_a=\{x=0,\ y=0\}$; the $z$ coordinates are
normal to the page and are not shown.
The thin lines are the coordinate curves
$Y=\text{const}$ and $x=\text{const}$ in the region $|Y|<c$
(which is disjoint from $\beta[X;C_a]^*\bX$), and their images under the
blow-down map $\beta[X;C_a]$.}
\label{fig:ff-coords}
\end{figure}

We again return to the Euclidean setting. In particular $X=\Snp$.
We denote the blow-down map
by $\beta[X;C_a]:[X;C_a]\to X$. Now $S^+ NC_a$ is a hemisphere bundle
over $C_a$, which can be identified with the
radial compactification of the normal bundle of $C_a$ in $\bX$
whose fibers can in turn be identified with  $X^a$. To see this
in more concrete terms,
we proceed by finding local coordinates on $[X;C_a]$ explicitly.
It is convenient to do so by using projective coordinates rather
than the standard polar coordinates.
Near $C_a$ in $\Snp$ we have $|w_a|>c|w^a|$ for some $c>0$.
Hence, near any point $p\in C_a$ one of the coordinate functions $(w_a)_j$
which we may take to be $(w_a)_m$, satisfies $|(w_a)_m|>c'|(w_a)_j|$,
$|(w_a)_m|>c'|w^a|$ for some $c'>0$. Taking into account the
coordinate form of $\SP$ we see that
\begin{equation}\label{eq:Snp-C_a-coords}
x=|(w_a)_m|^{-1},\ z_j=\frac{(w_a)_j}{|(w_a)_m|}\ (j=1,\ldots,m-1),
\ y_j=\frac{(w^a)_j}{|(w_a)_m|}\ (j=1,\ldots,n-m)
\end{equation}
give coordinates on $\Snp$ near $p$. In these coordinates $C_a$ is
defined by $x=0$, $y=0$.
Correspondingly, we have coordinates
\begin{equation}\label{eq:ff-coords}
x,\ z_j\ (j=1,\ldots,m-1),\ Y_j=y_j/x\ (j=1,\ldots,n-m),
\end{equation}
i.e.
\begin{equation}\label{eq:ff-coords-2}
x=|(w_a)_m|^{-1},\ z_j=\frac{(w_a)_j}{|(w_a)_m|}\ (j=1,\ldots,m-1),
\ Y_j=(w^a)_j\ (j=1,\ldots,n-m)
\end{equation}
near the interior of the front face $\ff$
of the blow-up $[X;C_a]$, i.e.\ near the
interior of $\ff=\beta[X;C_a]^*C_a$; see Figure~\ref{fig:ff-coords}.

Near the corner $\partial\beta[X;C_a]^*C_a=\beta[X;C_a]^*C_a\cap
\beta[X;C_a]^*\bX$,
in the lift of the region defined for some $k$ by
$|y_k|\geq c|y_j|$ for some $c>0$ and all $j\neq k$,
\begin{equation}\label{eq:mf-coords}
\xh=x/y_k,\ \Yh_j=y_j/y_k\ (j\neq k),\ y_k,\ z
\end{equation}
give coordinates. In terms of the original Euclidean variables these are
\begin{equation}\begin{split}
\xh=|(w^a)_k|^{-1},&\ z_j=\frac{(w_a)_j}{|(w_a)_m|}\ (j=1,\ldots,m-1),\\
&\ \Yh_j=\frac{(w^a)_j}{(w^a)_k}\ (j=1,\ldots,n-m,\ j\neq k),\ y_k=
\frac{(w^a)_k}{|(w_a)_m|}.
\end{split}\end{equation}
Since in every region near the lift $\beta[X;C_a]^*C_a$ of $C_a$ we can use
one of these coordinate systems, and since away from there we can use
coordinates as in \eqref{eq:Snp-C_a-coords} but with $w_a$ and $w^a$
interchanged, we have proved the following lemma.

\begin{lemma}\label{lemma:blowup-1}
Suppose that $X=\Snp$ and let $\beta=\beta[X;C_a]$ be the
blow-down map. Then the pull-back $\beta^*(\SP^{-1})^*\pi^a$ of
$\pi^a:\Rn\to X^a$ extends to a $\Cinf$ map, which we also
denote by $\pi^a$,
\begin{equation}\label{eq:geom-11}
\pi^a:[X;C_a]\to \Xb^a.
\end{equation}
Moreover, if $x^a$ is a boundary defining function on $\Xb^a$
(e.g.\ $x^a=|w^a|^{-1}$ for $|w^a|>1$), then $\rho_{\bX}=
(\pi^a)^* x^a$ is a defining function for
the lift of $\partial X$ to $[X;C_a]$, i.e.\ for $\beta^*\partial X$.
\end{lemma}

\begin{figure}[ht]
\begin{center}
\mbox{\epsfig{file=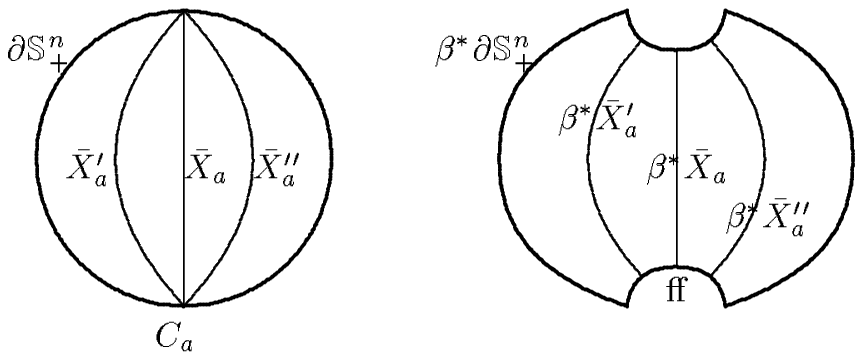}}
\end{center}
\caption{The blowup of $C_a$ in $\Snp$;
$\beta=\beta[\Snp;C_a]$ is the blow-down
map and $\ff=\beta^* C_a$. $X'_a$ and $X''_a$ denote translates
of $X_a$ in $\Rn$, $\Xb'_a=\cl(\SP(X'_a))$, etc. Note that the lifts
of $\Xb_a$, $\Xb'_a$ and $\Xb''_a$ become disjoint on $[\Snp;C_a]$.}
\label{fig:ressp}
\end{figure}

\begin{cor}
Suppose that $X=\Snp$, $f\in S^r_\phg(X^a)$. Then
\begin{equation}
(\pi^a)^*f\in\rho_{\bX}^{-r}\Cinf([X;C_a]).
\end{equation}
Here, following the previous lemma, we regard $\pi^a$ as the map
in \eqref{eq:geom-11},
and $\rho_{\bX}$ is
the defining function of $\beta[X;C_a]^*\bX$, i.e.\ of the lift
of $\bX$, and the
subscript $\phg$ refers to classical (one-step polyhomogeneous) symbols
(see Figure~\ref{fig:ressp}).
\end{cor}

This corollary shows that for a Euclidean many-body Hamiltonian, $H=\Delta+
\sum_a V_a$, $V_a$ becomes a nice function on the compact resolved
space $[\Snp;C_a]$. Thus, to understand $H$, we need to blow up
{\em all} the $C_a$. In order to analyze this iterated blow-up
procedure, it is convenient to generalize
the clean intersection properties to
manifolds with corners $X$.

Let $X$ be a
manifold with corners. Given a finite
family $\calF$ of closed p-submanifolds
$F_i$ of $X$ we say that $\calF$ is a cleanly intersecting family if it is
closed under intersection (in the sense that any two members are either
disjoint, or their intersection is in the family) and
for any $i$ and $j$, $\{F_i,F_j\}$ form a normal collection in the sense
of Melrose \cite[Chapter~V]{RBMDiff}.
Thus, for any point $p\in F_i\cap F_j$ there are local coordinates
on a neighborhood $U$ of $p$ such that with some index sets $I'$, $I''$,
\begin{equation}
F_i\cap U=\{x_r=0,\ r\in I'_i,\ y_s=0,\ s\in I''_i\},
\end{equation}
\begin{equation}
F_j\cap U=\{x_r=0,\ r\in I'_j,\ y_s=0,\ s\in I''_j\};
\end{equation}
here the $x_k$ are defining functions of the
boundary hypersurfaces through $p$.
This simply means that there is a common product decomposition for any
pair of elements of $\calF$. In particular, if $X$ is a manifold without
boundary, then this simply means that the $F_i$ pairwise intersect cleanly.
Hence, $(X,\calC)$ is a many-body space if and only if $\calC$ is a
cleanly intersecting family in $\bX$ which includes $\bX$.

Just as in the case of a many-body space, inclusions gives a partial order
on a cleanly intersecting family $\calF$. Thus, $F\in\calF$ is minimal with
respect to inclusion if there is no $F'\in\calF$ such that $F'\neq F$,
$F'\subset F$. Since $\calF$ is closed under intersection, this means
exactly that for all $F'\in\calF$ either $F'$ and $F$ are disjoint, or
$F\subset F'$.

\begin{lemma}
Let $\calF$ be a cleanly intersecting family of p-submanifolds of $\bX$.
Suppose that $F\in\calF$ is minimal with respect to inclusion.
Then the lifted family, $\calF'$, consisting of the lifts of $F_j$,
distinct from $F$, to $[X;F]$, is also a cleanly intersecting family.
\end{lemma}

\begin{proof}
We claim that for any $F_i, F_j\in\calF$ the
4-tuple $\{F,F_k,F_i,F_j\}$, $F_k=F_i\cap F_j$,
is a normal collection in the sense of Melrose.
Indeed, this is clear if $F_k$ is disjoint from $F$; otherwise
$F\subset F_k$ by our assumption.

So assume that $F\subset F_k$.
By the normality of $\{F, F_k\}$, near any point $p$ in $F$
there are local coordinates $x_r$, $y_s$, on
$X$ such that
\begin{equation}
F_k=\{x_r=0,\ r\in I'_k,\ y_s=0,\ s\in I''_k\},
\end{equation}
\begin{equation}
F=\{x_r=0,\ r\in I',\ y_s=0,\ s\in I''\},
\end{equation}
and $I'_k\subset I'$, $I''_k\subset I''$.
Similarly, by the normality of $\{F_i,F_j\}$ there are local coordinates
$x'_r$, $y'_s$ near $p$ on $X$ such that
\begin{equation}\label{eq:pr-5a}
F_i=\{x'_r=0,\ r\in I'_i,\ y'_s=0,\ s\in I''_i\},
\end{equation}
\begin{equation}\label{eq:pr-5b}
F_j=\{x'_r=0,\ r\in I'_j,\ y'_s=0,\ s\in I''_j\},
\end{equation}
Thus,
\begin{equation}\label{eq:pr-5c}
F_k=\{x'_r=0,\ r\in I'_i\cup I'_j,\ y'_s=0,\ s\in I''_i\cup I''_j\}.
\end{equation}
Thus, the differentials of the coordinates $x'_r$, $r\in I'_i\cup I'_j$,
and $y'_s$, $s\in I''_i\cup I''_j$, span the conormal bundle of $F_k$.
The same holds for the differentials of
$x_r$, $r\in I'_k$, $y_s$, $s\in I''_k$. It follows that the differentials
of $x'_r$, $r\in I'_i\cup I'_j$, $x_r$, $r\nin I'_k$, $y'_s$, $s\in I''_i
\cup I''_j$, $y_s$, $s\nin I''_k$ are independent at $F_k$ in a
coordinate neighborhood of $p$, so these functions give local coordinates
on $X$ near $p$ in terms of which $F$, $F_k$, $F_i$ and $F_j$ have
common product decomposition: $F_i$, $F_j$ and $F_k$ given by
\eqref{eq:pr-5a}-\eqref{eq:pr-5c}, and $F$ by
\begin{equation}\begin{split}
F=\{x'_r=0,\ r\in I'_i\cup I'_j,&\ x_r=0,\ r\in I'\setminus I'_k,\\
&\ y'_s=0,\ s\in I''_i\cup I''_j,\ y_s=0,\ s\in I''\setminus I''_k\}.
\end{split}\end{equation}
This proves that $\{F, F_k, F_i, F_j\}$ is indeed a normal collection.
Hence, by \cite[Lemma V.11.2]{RBMDiff}, it lifts to a normal
collection of p-submanifolds on $[X;F]$. Writing $\beta$ for the blow-down
map, and $\beta^*F_k$ for the lift of $F_k$, etc., we see in particular
that $\{\beta^*F_i, \beta^*F_j\}$ is
a normal collection whose intersection is $\beta^* F_k$ if $F_k\neq F$,
and is empty otherwise. Putting together these facts we see that we have
proved the lemma.
\end{proof}

This lemma allows us to define $[X;\calF]$ if $\calF$ is a cleanly
intersecting family of p-submanifolds of $\bX$. We do this by putting a
total order on $\calF$ which is compatible with the partial order
given by inclusion. This can always be accomplished: pick a minimal
element with respect to inclusion, and make it the minimal element of
the total order. Proceeding inductively, if we already placed a total
order on $\calF'\subset\calF$, we choose any $F\in\calF\setminus\calF'$
which is minimal with respect to inclusion in $\calF\setminus\calF'$,
and extend the total order to $\calF'\cup\{F\}$ by making $F$ the maximal
element with respect to it. Having imposed a total order on $\calF$
which is compatible with inclusion, we define $[X;\calF]$ to be the
blow up $[X;F_1,F_2,\ldots,F_n]$ where $\calF=\{F_1,F_2,\ldots,F_n\}$
and $F_1<F_2<\ldots<F_n$, $<$ being
the total order. Of course, a priori $[X;\calF]$ depends on the total
order. The following lemma shows that this is not the case.

\begin{lemma}
If $\calF$ is a cleanly intersecting family and $<$, $<'$ are total
orders on it which are compatible with inclusion, then the blow ups
\begin{equation}
[X;F_1,F_2,\ldots,F_n],\qquad F_1<F_2<\ldots<F_n,
\end{equation}
\begin{equation}
[X;F'_1,F'_2,\ldots,F'_n],\qquad F_1'<'F_2'<'\ldots<'F_n'
\end{equation}
are canonically diffeomorphic.
\end{lemma}

\begin{proof}
Since any total order compatible with inclusion can be obtained from any
other one by repeatedly interchanging the order of adjacent elements, but
keeping the order compatible with inclusion, it suffices to show that
\begin{equation}
[X;F_1,\ldots,F_k,F_{k+1},\ldots,F_n]\ \text{and}
\ [X;F_1,\ldots,F_{k+1},F_k,\ldots,F_n]
\end{equation} are naturally isomorphic if both of these total orders
respect inclusion. Now, either $F_k\cap F_{k+1}=\emptyset$, in which
case the statement is clearly true, or $F_k\cap F_{k+1}=F_j$ for some $j$.
Since inclusion is respected, we must have $j<k$. But upon the blow up
of their intersection, any two closed p-submanifolds with normal
intersection lift to be disjoint. Hence, on $[X;F_1,\ldots,F_{k-1}]$ the
lifts $\beta^* F_k$ and $\beta^* F_{k+1}$ are disjoint, and thus they
can be blown up in either order. This proves the lemma.
\end{proof}

Correspondingly, $[X;\calF]$ is defined independently of the total order
used in the definition of the blown up space, assuming that it
respects inclusion, so we can speak about $[X;\calF]$ without
specifying such a total order.

If $F_i\in \calF$, we can always specify the total order so that every
$F_j\in\calF$ with $F_j<F_i$ satisfies $F_j\subset F_i$. Then the blow-up
of $F_i$ commutes with all the ones preceeding it. Hence, any function
that is smooth on $[X;F_i]$ pulls back to be smooth on $[X;\calF]$.
Applying this in the Euclidean many-body setting we conclude that

\begin{lemma}
Suppose that $X=\Snp$ and $\calX$ is a linear family of subspaces of $\Rn$
as in the introduction. Then $V=\sum_a V_a$, $V_a\in S^{-m}_\phg(X^a)$,
lifts to be an element of $\rho_{\bX}^m\Cinf([X;\calC])$ where
$\rho_{\bX}$ is the defining function of the lift of $C_0=\bX$ under the
blow-down map
\begin{equation}
\betaSc=\beta[X;\calC]:[X;\calC]\to X.
\end{equation}
\end{lemma}

Our main interest is the study of differential operators, in particular
the analysis of many-body Hamiltonians $H$.
For this purpose we next investigate
how vector fields lift under the blow up.
First,
we define $\Vb(X;\calF)$ as the Lie algebra of smooth vector fields on $X$
which are tangent to the boundary faces of $X$ and to each element of $\calF$.

\begin{lemma}\label{lemma:tgt-vec-blow-up}
Each element of $\Vb(X;\calF)$ lifts to an element of $\Vb([X;\calF])$.
\end{lemma}

\begin{proof}
It suffices to show that $V\in\Vb(X;\calF)$ lifts to be an element of
$\Vb([X;F];\calF')$ where $F$ is minimal with respect to inclusion and
\begin{equation}
\calF'=\{\beta^*F':\ F'\in\calF\setminus\{F\}\}.
\end{equation}
Taking into account that for any $F'\neq F$, $\{F,F'\}$ is a normal
collection of p-submanifolds of $X$, this claim follows from
\cite[Proposition~V.11.1]{RBMDiff}, or it can be checked directly
by using projective coordinates on $[X;F]$.
\end{proof}

\begin{rem}
It is {\em not} the case in general that $\Vb(X;\calF)$ lifts to span
$\Vb([X;\calF])$ over $\Cinf([X;\calF])$. This statement is true,
however, if $\calF$ is a normal collection (i.e.\ all elements of $\calF$
have product decomposition in the same coordinate system, not just pairs
of elements), see \cite[Proposition~V.11.1]{RBMDiff}.
\end{rem}

We can now introduce the appropriate class of differential and
pseudo-differential operators on many-body spaces $(X,\calC)$.
These will include many-body Hamiltonians
in the Euclidean setting as well as their resolvents (in the resolvent set).

First, we recall from \cite{RBMSpec} Melrose's definition of
the Lie algebra of `scattering vector fields' $\Vsc(X)$, defined for
every manifold with boundary $X$. Thus,
\begin{equation}
\Vsc(X)=x\Vb(X)
\end{equation}
where $\Vb(X)$ is the set of smooth vector fields
on $X$ which are tangent to $\bX$. If $(x,y_1,\ldots, y_{n-1})$
are coordinates on $X$ where
$x$ is a boundary defining function, then locally a basis of $\Vsc(X)$
is given by
\begin{equation}
x^2\partial_x,\ x\partial_{y_j},\ j=1,\ldots,n-1.
\end{equation}
Correspondingly, there is a vector bundle $\Tsc X$ over $X$, called
the scattering tangent bundle of $X$, such that $\Vsc(X)$ is the set of
all smooth sections of $\Tsc X$:
\begin{equation}
\Vsc(X)=\Cinf(X;\Tsc X).
\end{equation}
The
dual bundle of $\Tsc X$ (called the scattering cotangent bundle) is
denoted by $\sct X$.
Thus, covectors $v\in\sct_p X$, $p$ near
$\bX$, can be written as $v=\tau\,\frac{dx}{x^2}+\mu\cdot\frac{dy}{x}$.
Hence, we have local coordinates $(x,y,\tau,\mu)$ on $\sct X$ near $\bX$.
The scattering density bundle $\Osc X$
is the density bundle associated to
$\sct X$, so locally near $\bX$ it is spanned by $x^{-n-1}\,dx\,dy$ over
$\Cinf(X)$.
Finally, $\Diffsc(X)$ is the algebra of differential operators generated
by the vector fields in $\Vsc(X)$; $\Diffsc^m(X)$ stands for scattering
differential operators of order (at most) $m$.

To establish the
relationship between the scattering structure and
the Euclidean scattering theory, we
introduce local coordinates on $X$ near $p\in\bX$ as above, and
use these to identify the coordinate neighborhood $U$ of $p$ with
a coordinate patch $U'$ on the closed upper hemisphere $\Snp$ (which is
just a closed ball) near its boundary.
Such an identification preserves the scattering structure since this
structure is
completely natural. We
further identify $\Snp$ with $\Rn$ via the
radial compactification $\SP$ as in \eqref{eq:SP-def}.
The constant coefficent vector fields $\partial_{w_j}$ on $\Rn$ lift
under $\SP$
to give a basis of $\Tsc\Snp$. Thus, $V\in\Vsc(\Snp)$ can
be expressed as (ignoring the lifting in the notation)
\begin{equation}
V=\sum_{j=1}^n a_j\partial_{w_j},\quad a_j\in\Cinf(\Snp).
\end{equation}
As mentioned in the introduction, $a_j\in\Cinf(\Snp)$ is equivalent to
requiring that $\SP^* a_j$ is a classical (i.e.\ one-step polyhomogeneous)
symbol of order $0$ on $\Rn$. This description also shows that
the positive Euclidean Laplacian, $\Delta$, is an element of
$\Diffsc^2(\Snp)$, and that $\Osc\Snp$ is spanned by the pull-back of
the standard Euclidean density $|dw|$.

If $X$ is a manifold with boundary then any element of $\Vsc(X)=x\Vb(X)$ is
automatically tangent to any submanifold $C$ of $\bX$. Hence, due
to Lemma~\ref{lemma:tgt-vec-blow-up}, we can
define the algebra of many-body differential operators as shown by
the following proposition.

\begin{prop}
If $(X,\calC)$ is a many-body space, then $\Vsc(X)$ lifts to a subalgebra
of $\Vb([X;\calC])$. Correspondingly,
\begin{equation}
\DiffSc(X,\calC)=\Cinf([X;\calC])\otimes_{\Cinf(X)}\Diffsc(X)
\end{equation}
is an algebra.
\end{prop}

\begin{proof}
By the first part of the statement, for any $V\in\Vsc(X)$, $f\in\Cinf
([X;\calC])$, the commutator $[V,f]=Vf$ is in $\Cinf([X;C])$.
\end{proof}

Since in Euclidean many-body scattering
$\Delta\in\Diffsc^2(\Snp)$ and $V=\sum_a V_a\in\Cinf([\Snp;\calC])$,
it follows immediately that $H=\Delta+V\in\DiffSc^2(\Snp,\calC)$.

\section{Many-body pseudo-differential calculus}\label{sec:calculus}
Let $(X;\calC)$ be a many-body space, and $\betaSc:[X;\calC]\to X$ the
blow-down map. There are two equivalent way of defining many-body
pseudo-differential operators. We can either specify their kernels as
conormal distributions on an appropriately resolved space, or we can
define them as the quantization of certain symbols. Here we give
both definitions and show their equivalence.

The b-double space, $\XXb$,
has been defined by Melrose as $[X^2;(\bX)^2]$. The front face of the
blow-up is called the b-front face and is denoted by $\ffb$, while the
lifts of the left and right boundary hypersurfaces of $X^2$, i.e.\ of
$\bX\times X$ and $X\times\bX$ are denoted by $\lf$ and $\rf$ respectively.
The diagonal $\Delta$
of $X^2$ lifts to a p-submanifold $\Deltab$ of $\XXb$ which intersects
$\partial\XXb$ in the interior of the b-front face, $\ffb$. (The definition
of p-submanifolds and the blow-up process were
discussed at the beginning of the
previous section.) Moreover,
$\Deltab$ is naturally diffeomorphic to $X$. Hence, $\calC$ can be
regarded as a collection $\calC'$ of submanifolds of $\Deltab$, and, since
$\Deltab$ is a p-submanifold of $\XXb$, these
submanifolds form a cleanly intersecting family in $\XXb$. Therefore,
the blow up
\begin{equation}\label{eq:XXSc-def}
\XXSc=[\XXb;\calC']
\end{equation}
is well-defined by our previous results. Note that $\bX\in\calC$ by our
assumption, so the definition includes the blow up of the lift of
$\partial\Deltab$. It is easy to see that this space coincides with the
$\XXSc$ defined in \cite{Vasy:Propagation} if $(X,\calC)$ is a 3-body
space.

Noting that even $\calC'\cup\{\Deltab\}$ is a cleanly intersecting family,
we conclude that $\Deltab$ lifts to a a p-submanifold, $\DeltaSc$, of
$\XXSc$. Correspondingly, we define the set of many-body
pseudo-differential operators by
\begin{equation}\label{eq:PsiScc-def}
\PsiScc^{m,l}(X,\calC)=\{\kappa\in\bcon^{m,l}(\XXSc,\DeltaSc;\Osc_R):
\ \kappa\equiv 0\ \text{at}\ \beta^*\ffb\cup\beta^*\lf\cup\beta^*\rf\};
\end{equation}
here $\Osc_R$ is the pull-back of the scattering density bundle from
the right factor and $\beta:\XXSc\to\XXb$ is the
blow-down map. Similarly we define the corresponding polyhomogeneous operators
\begin{equation}\label{eq:PsiSc-def}
\PsiSc^{m,l}(X,\calC)=\{\kappa\in \rho^l I^m(\XXSc,\DeltaSc;\Osc_R):
\ \kappa\equiv 0\ \text{at}\ \beta^*\ffb\cup\beta^*\lf\cup\beta^*\rf\}
\end{equation}
where $\rho$ is the total boundary defining function of $\XXSc$.
In particular, conormal distributions of order $-\infty$ are smooth
functions, so
\begin{equation}
\PsiSc^{-\infty,l}(X,\calC)=\{\kappa\in \rho^l \Cinf(\XXSc,\DeltaSc;\Osc_R):
\ \kappa\equiv 0\ \text{at}\ \beta^*\ffb\cup\beta^*\lf\cup\beta^*\rf\},
\end{equation}
i.e.\ the kernels of operators in $\PsiSc^{-\infty,l}(X,\calC)$ are smooth up
to all boundary hypersurfaces of $\XXSc$ (at least if $l$ is a non-negative
integer), and vanish to infinite order at the lift of evey boundary
hypersurface of $\XXb$.
Tensoring with vector bundles defines $\PsiScc^{m,l}(X,\calC;E,F)$
and $\PsiSc^{m,l}(X,\calC;E,F)$
for vector bundles $E$ and $F$ over $X$ as usual.

Since for all $F\in\calC'$ we have $F\subset\partial\Deltab$, we can
do the blow up of $\partial\Deltab\in\calC'$ first, before blowing up
other elements of $\calC'$ (normally we would do this blow up last by
our total order construction). It follows that $\XXSc$ is a blow up of
the space $\XXsc=[\XXb;\partial\Deltab]$. Hence, conormal
distributions on $\XXsc$ pull back to be conormal on $\XXSc$. Since
the kernels of scattering pseudo-differential operators are conormal
to $\Deltasc$ and to the boundary of $\XXsc$ with infinite order vanishing
at every boundary face except the scattering front face, we conclude that
these kernels pull back to $\XXSc$ to be elements of the kernel space defined
in \eqref{eq:PsiScc-def}, so $\Psiscc^{m,l}(X)\subset\PsiScc^{m,l}
(X,\calC)$.

Suppose now that $X=\Snp$ and $\calC$ is a cleanly intersecting family
of submanifolds of $\bX=\partial\Snp=\Sn$. Here we {\em do not} assume that
$\calC$ arises from a family $\calX$ of linear subspaces of $\Rn$.
An equivalent definition of $\PsiScc^{m,l}(\Snp,\calC)$ is the following.
Suppose that
\begin{equation}\label{eq:PsiScc-5}
a\in\bcon^{-m,l}([\Snp;\calC]\times\Snp).
\end{equation}
That is,
identifying $\interior(\Snp)$ and $\interior([\Snp;\calC])$ with $\Rn$
as usual (via $\SP^{-1}$), suppose that
$a\in\Cinf(\Rn_w\times\Rn_\xi)$ has the following property.
For every $P\in\Diffb^k(\Snp)$, acting on the second factor of
$\Snp$ (i.e.\ in the $\xi$ variable), and
$Q\in\Diffb^{k'}([\Snp;\calC])$, acting
on the first factor of $\Snp$ (i.e.\ in the $w$ variable),
$k,k'\in\Nat$,
\begin{equation}\label{eq:symb-est-3}
PQa\in\rho_{\infty}^{-m}\rho_{\partial}^l L^{\infty}(\Snp\times\Snp)
\end{equation}
where $\rho_\infty$ and $\rho_\partial$ are defining functions of the first
and second factors of $\Snp$ respectively. Let $A=q_L(a)$ denote the
left quantization of $a$:
\begin{equation}\label{eq:q_L-def}
Au(w)=(2\pi)^{-n}\int e^{i(w-w')\cdot\xi}a(w,\xi)u(w')\,dw'\,d\xi,
\end{equation}
understood as an oscillatory integral. Then $A\in\PsiScc^{m,l}(\Snp,\calC)$.
Indeed, the kernel of $A$ is
\begin{equation}\label{eq:PsiSc-23}
K(w,w')=\tilde a(w,w-w')
\end{equation}
where
$\tilde a$ is the inverse Fourier transform of $a$ in the $\xi$ variable,
i.e.\ $\tilde a=\Frinv_\xi a$. Thus, $\tilde a(w,W)$ is smooth away from
$W=0$, is conormal to $W=0$, and it is rapidly decreasing with all
derivatives in $W$. More precisely, the rapid decay means that for all $k$ and
$Q\in\Diffb([\Snp;\calC])$ and for all $\alpha$,
\begin{equation}
\sup_{|W|\geq 1,\ w\in\Rn}(|w|^l|W|^k |Q_w D^\alpha_W\tilde a(w,W)|)<\infty.
\end{equation}
Taking into account the geometry of $\XXSc$, in particular that $|w-w'|^{-1}$
vanishes at all faces of the blow-up \eqref{eq:XXSc-def}
but the front faces (i.e.\ it vanishes at $\beta^*\lf$, $\beta^*\rf$ and
$\beta^*\ffb$), we see
that $K$ vanishes to infinite order at these faces. Similar arguments
describe the behavior of $K$ near $\DeltaSc$, proving that
$A\in\PsiScc^{m,l}(\Snp,\calC)$.

Conversely, if $A\in\PsiScc^{m,l}(\Snp,\calC)$ then there exists $a$ satisfying
\eqref{eq:symb-est-3} such that $A=q_L(a)$. Namely, we let
$\tilde a(w,W)=K(w,w-W)$ and let $a$ be the Fourier transform of $\tilde a$
in $W$. The conormal estimates for $K$ (hence for $\tilde a$) give the
symbolic estimates \eqref{eq:symb-est-3} for $a$.

Similar conclusions hold
for the right quantization $B=q_R(b)$ of a symbol $b$:
\begin{equation}
Bu(w)=(2\pi)^{-n}\int e^{i(w-w')\cdot\xi}b(w',\xi)u(w')\,dw'\,d\xi.
\end{equation}
In addition, the polyhomogeneous class $\PsiSc^{m,l}(\Snp,\calC)$
is given by the quantization of symbols
\begin{equation}\label{eq:cl-symb-def}
a\in\rho_{\infty}^{-m}\rho_{\partial}^l\Cinf([\Snp;\calC]\times\Snp).
\end{equation}
Since
differential operators $\sum a_\alpha(w) D^\alpha$ are just the
left quantization of the symbols $a(w,\xi)=\sum a_\alpha(w)\xi^\alpha$, it
follows immediately that
\begin{equation}
\DiffSc^m(X,\calC)\subset\PsiSc^m(X,\calC).
\end{equation}
This conclusion also follows directly from the description of the kernels
since the kernel of a differential operator is a differentiated
delta-distribution associated to the diagonal.

Note that, as usual, one can allow symbols $a$ depending on $w$, $w'$ and
$\xi$, so e.g.\ if $a\in\rho_\infty^{-m}\rho_{\partial,L}^l
\rho_{\partial,R}^{l'}\Cinf
([\Snp;\calC]\times[\Snp;\calC]\times\Snp)$, $\rho_{\partial,L}$ and
$\rho_{\partial,R}$ denoting total boundary defining functions of the
first and second factor of $[\Snp;\calC]$ respectively (i.e.\ they are
pull-backs of a boundary defining function of $\Snp$), then
\begin{equation}\label{eq:Psop-37}
Au(w)=(2\pi)^{-n}\int e^{i(w-w')\cdot\xi}a(w,w',\xi)\,u(w')\,dw'\,d\xi
\end{equation}
defines an operator $A\in\PsiSc^{m,l+l'}(\Snp,\calC)$.

This characterization allows the application of the standard tools of
the theory of pseudo-differential operators. In particular,
if $A\in\PsiSc^{m,l}(X,\calC)$ is written as the left quantization
of a symbol $a$ and $B\in\PsiSc^{m',l'}(X,\calC)$ is written as
the right quantization of a symbol $b$, so
\begin{equation}
a\in
\rho_{\infty}^{-m}\rho_{\partial}^l\Cinf([\Snp;\calC]\times\Snp),
\ b\in
\rho_{\infty}^{-m'}\rho_{\partial}^{l'}\Cinf([\Snp;\calC]\times\Snp),
\end{equation}
then the operator $AB$ is given by
\begin{equation}
ABu(w)=(2\pi)^{-n}\int e^{i(w-w')\cdot\xi}a(w,\xi)b(w',\xi)\,u(w')\,dw'\,d\xi.
\end{equation}
Here $c(w,w',\xi)=a(w,\xi)b(w',\xi)$ is in
$\rho_\infty^{-m-m'}\rho_{\partial,L}^l
\rho_{\partial,R}^{l'}\Cinf
([\Snp;\calC]\times[\Snp;\calC]\times\Snp)$, so we conclude that
$AB\in\PsiSc^{m+m',l+l'}(X,\calC)$. In addition, the adjoint $A^*$ of $A$ is
the right quantization of $\bar a$, so $A^*\in\PsiSc^{m,l}(X,\calC)$.
Analogously, $\PsiScc(\Snp,\calC)$ is also closed under composition and
adjoints. These statements can be seen also from the standard
more explicit formulae. For example, if $B$ is the left quantization
of a symbol $b'$, the
composition formula, including the remainder terms,
only involves derivatives of the form $D_\xi^\alpha
D_w^\alpha b'$, and $D_w^\alpha\in\Diffsc^{|\alpha|}(\Snp)\subset
\Diffb^{|\alpha|}([\Snp;\calC])$, so we see that $\PsiScc(\Snp,\calC)$ is
closed under composition.

This discussion can be carried over to arbitrary manifolds with boundary
$X$ by locally identifying $X$ with $\Snp$ and using that our arguments
are local in $\Snp$. More precisely,
suppose that $\{U_1,\ldots, U_k\}$ is an open cover of $X$ by coordinate
patches, and identify each $U_i$ with a coordinate patch $U'_i$ of
$\Snp$. We write $\phi_i:U_i\to U'_i$ for the identification.
Let $\calC'_i$ denote the family given by the image of elements
of $\calC$ in $U'_i$. Then $A\in\PsiSc^{m,l}(X,\calC)$ if and only if
there exists operators $A'_i\in\PsiSc^{m,l}(\Snp;\calC')$ with kernel
supported in the inverse image of $U'_i\times U'_i$ in $(\Snp)^2_\Scl$
and $R\in\dCinf(X\times X;\Osc_R)$ such that
\begin{equation}
A=\sum_i (\phi_i^*A'_i
(\phi_i^{-1})^*)+R.
\end{equation}
Note that the support condition on $A'_i$ ensures
that this expression makes sense. To see this, just introduce a
partition of unity $\rho_i\in\Cinf(X)$ subordinate to the cover, and
let $\psi_i\in\Cinf(X)$ be identically $1$ in a neighborhood of $\supp
\rho_i$. Then
\begin{equation}\label{eq:Psop-17}
A=\sum_i A\rho_i=\sum_i\psi_i A\rho_i+\sum_i(1-\psi_i)A\rho_i.
\end{equation}
It is straightforward to check directly from the definition of
$\PsiSc^{m,l}(X,\calC)$ that the last terms is given by a kernel in
$\dCinf(X\times X;\Osc_R)$, while $A'_i=
(\phi_i^{-1})^*\psi_i A\rho_i\phi_i^*\in\PsiSc^{m,l}(\Snp,\calC'_i)$ with the
claimed support properties.
Thus, our results for $\PsiSc^{m,l}(\Snp,\calC)$ immediately show
the following theorem.

\begin{thm}
Both $\PsiSc(X,\calC)$ and $\PsiScc(X,\calC)$
are $*$-algebras (with respect to
composition and taking adjoints).
\end{thm}

Since $\Psiscc^{m,0}(\Snp,\calC)\subset\Psi_\infty^m(\Rn)$,
where $\Psi_\infty^m
(\Rn)$ is the class of pseudo-differential operators defined by
H\"ormander \cite[Section~18.1]{Hor},
arising by a quantization of symbols $a\in\Cinf(\Rn\times\Rn)$
satisfying
\begin{equation}
|D_w^\alpha D_\xi^\beta a(w,\xi)|\leq C_{\alpha\beta}
\langle\xi\rangle^{m-|\beta|},
\end{equation}
and
\begin{equation}
\Psi^m_\infty(\Rn):\langle w\rangle^{-s}H^r(\Rn)\to\langle w\rangle^{-s}
H^{r-m}(\Rn),
\end{equation}
we immediately deduce the boundedness of elements of $\Psiscc^{m,l}(X,\calC)$
between
the appropriate weighted Sobolev spaces.

\begin{thm}
If $A\in\Psiscc^{m,l}(X,\calC)$ then $A:\Hsc^{r,s}(X)\to\Hsc^{r-m,s+l}(X)$ is
bounded.
\end{thm}

There is another way of characterizing the calculus $\PsiScc^{m,l}(\Snp,
\calC)$ via H\"ormander's Weyl calculus (see \cite[Section~18.5]{Hor}).
We describe it briefly here, only considering the Euclidean setting where
the $C_a$ arise from linear subspaces $X_a$;
it is straightforward to check that it
agrees with the definition we have given above
in terms of quantization of symbols
as in \eqref{eq:PsiScc-5}.
Namely, $\Psiscc^{\infty,-\infty}(\Snp)$ is just the calculus on $\Rn$ arising
from the metric
\begin{equation}
g^{(0)}=\frac{dw^2}{\langle w\rangle^2}+\frac{d\xi^2}{\langle \xi\rangle^2}.
\end{equation}
Similarly, if we take $\calC'$ to consist of a single element $C_a$, $a\neq 0$,
and if $(w_a,w^a)$ is the usual splitting of the coordinates, then
$\PsiScc^{\infty,-\infty}(\Snp,\calC')$ arises from the metric
\begin{equation}
g^{(a)}=\frac{dw_a^2}{\langle w\rangle^2}+\frac{(dw^a)^2}{\langle w^a\rangle^2}
+\frac{d\xi^2}{\langle \xi\rangle^2}.
\end{equation}
In the three-body problem, $C_a\cap C_b=\emptyset$ if $a,b\neq 0$,
we define the metric by localizing the $g^{(a)}$, i.e.\ we consider a
partition of unity $\phi_a\in\Cinf(\Snp)$,
$a\in I$, $\supp \phi_a\cap C_b=\emptyset$
unless $b=0$, and define the metric
\begin{equation}
g=\sum_a \phi_a\, g^{(a)}.
\end{equation}
(Here the $\phi_a$ are pulled back to the cotangent bundle by the bundle
projection.)
Since the $g^{(a)}$ are equivalent near $C'_0$, it follows that $g$ is
indeed slowly varying. Note that if $\phi_a$ is supported close to $C_a$,
which we can arrange by enlarging the support of $\phi_0$,
$dw_a^2/\langle w\rangle^2$ above can be replaced by
$dw_a^2/\langle w_a\rangle^2$.

In general we simply repeat this procedure. Thus,
to define the appropriate
metric on $T^*X^c$ if it has been defined on $T^*X^a$ for every
$a$ with $X^a\subset X^c$, we define a partition of unity $\phi_a\in
\Cinf(\Xb^c)$ with $\supp\phi_a\cap C_b^c=\emptyset$ unless $C_a^c
\subsetneq C_b^c$. Here $X^c=X^a\oplus X^c_a$ and
$C_a^c=\partial\Xb^c\cap \cl(X_a^c)$). We extend the
metric $g^{a}$ on $T^*X^a$ to a symmetric 2-cotensor on $T^* X^c$
using the orthogonal decomposition $X^c=X^a\oplus X^c_a$, and let
\begin{equation}
g^{(a)}=g^a+\frac{(dw^c_a)^2}{\langle w^c_a\rangle^2}+
\frac{(d\xi^c_a)^2}{\langle \xi^c\rangle^2}.
\end{equation}
Then
\begin{equation}
g^c=\sum_{a: X^a\subset X^c} \phi_a\, g^{(a)}
\end{equation}
gives the desired metric on $T^*X^c$.

After this brief discussion of the relationship of $\PsiSc^{*,*}(\Snp,\calC)$
with H\"ormander's Weyl calculus, we return to the general setting to
describe the principal symbol map and its analog at $\bX$.

\section{The principal symbol and the indicial operators}\label{sec:indicial}

Since the inclusion of $\Hsc^{r',s'}(X)$ to $\Hsc^{r,s}(X)$ is compact for
$r'>r$, $s'>s$, it suffices to understand $A\in\PsiScc^{m,l}(X,\calC)$ modulo
$\PsiScc^{m-1,l+1}(X,\calC)$ to analyze its spectral properties.
Now, H\"ormander's principal symbol map on $\Psi_\infty^m(\Rn)$ restricts to
a principal symbol map
\begin{equation}\label{eq:calc-23}
\prs_{\Scl,m}:\PsiSc^{m,0}(\Snp,\calC)\to\Shom^m(\ScT[\Snp;\calC]),
\end{equation}
$\Shom(\ScT[\Snp;\calC])$
denoting the space of smooth symbols which are homogeneous
of degree $m$. Due to its invariance and its local nature,
it immediately extends to a map
\begin{equation}
\prs_{\Scl,m}:\PsiSc^{m,0}(X,\calC)\to\Shom^m(\ScT [X;\calC]).
\end{equation}
We radially compactify the fibers of $\ScT [X;\calC]$ (i.e.\ replace
the vector spaces by balls) and let
$\ScS [X;\calC]$ be the
new boundary face (i.e.\ the boundary of $\ScT[X;\calC]$ at fiber-infinity).
This allows us to write $\prs_{\Scl,m}$ as a map
\begin{equation}\label{eq:calc-25}
\prs_{\Scl,m}:\PsiSc^{m,0}(X,\calC)
\to\Cinf(\ScS [X;\calC];(N^*\ScS[X;\calC])^{-m}).
\end{equation}
The line bundle $N^*\ScS[X;\calC]$ is locally spanned by the
pull-back of $d(|\xi|^{-1})$ from $\ScT[X;\calC]$, so \eqref{eq:calc-25}
is obtained from \eqref{eq:calc-23} by writing homogeneous functions
$a(w,\xi)$ of degree $m$ as $a_0(w,\xih)|\xi|^m$, $\xih=\xi/|\xi|$,
considering $a_0$ as a function on the cosphere bundle, and using
$N^*\ScS[X;\calC]$ to take care of the factor $|\xi|^m$ invariantly.
We then have a short exact sequence
\begin{equation}
0\to\PsiSc^{m-1,0}(X,\calC)\to\PsiSc^{m,0}(X,\calC)
\to\Cinf(\ScS [X;\calC];(N^*\ScS[X;\calC])^{-m})\to 0
\end{equation}
as usual.

Now, an operator $A\in\PsiSc^{m,0}(X,\calC)$ is certainly determined modulo
$\PsiSc^{m,1}(X,\calC)$ by the restriction of its kernel to the front faces of
the blow up \eqref{eq:XXSc-def}. Our next task is to construct a multiplicative
indicial operator from this restriction. We use oscillatory testing for this
purpose as was done in \cite{Vasy:Propagation-2}. We start by discussing the
effect of conjugation of $A$ by oscillatory functions.

\begin{lemma}
Suppose that $A\in\PsiSc^{m,l}(X,\calC)$ and $\ft\in\Cinf(X;\Real)$. Then
\begin{equation}
\At=e^{-i\ft/x}A e^{i\ft/x}\in\PsiSc^{m,l}(X,\calC).
\end{equation}
\end{lemma}

\begin{proof}
It is convenient to use the explicit description of $\PsiSc(X,\calC)$ in terms
of localization and quantization \eqref{eq:q_L-def}. Thus, we may assume
that $X=\Snp$. Note that the
pull-back of $\ft/x$ to $\Rn$ is a polyhomogeneous symbol of order $1$ which
we denote by $F$. Then the kernel of
$\At$ is $\Kt(w,w')=e^{i(F(w')-F(w))}K(w,w')$ where $K$ is the kernel of $A$.
But by the fundamental theorem of calculus
\begin{equation}
F(w')-F(w)=\sum_{j=1}^n(w'_j-w_j)\int_0^1\partial_j F(w+t(w'-w))\,dt,
\end{equation}
and $\partial_j F$ is a polyhomogeneous symbol of order $0$.
Taking into account the rapid decay of $K$ in $W=w-w'$ we immediately
conclude that $\Kt\in\bcon^{m,l}((\Snp)^2_\Scl,\DeltaSc;\KDSc)$ vanishing
with all derivatives at $\beta^*\ffb\cup\beta^*\lf\cup\beta^*\rf$, so,
returning to the global setting, $\At\in\PsiSc^{m,l}(X,\calC)$.
\end{proof}

We next discuss mapping properties on $\Cinf([X;\calC])$.

\begin{lemma}
If $A\in\PsiSc^{m,l}(X,\calC)$, $u\in x^r\Cinf([X;\calC])$, then $Au\in x^{r+l}
\Cinf([X;\calC])$.
\end{lemma}

\begin{proof}
This result essentially reduces to the fact that $\PsiSc(X,\calC)$
is an algebra.
Indeed, write $u=u\cdot 1$, and note that $Au=(AU)1$ where $B=AU$ denotes
the composite of $A$ with the multiplication operator $U$ by $u$. Since the
latter is in $x^r\DiffSc^0(X,\calC)$, hence in $\PsiSc^{0,r}(X,\calC)$,
we conclude that
$B\in\PsiSc^{m,l+r}(X,\calC)$. Thus, we only have to analyze $B1$. Again, we
can reduce the discussion to a local one. But writing
$B$ as the left quantization of a symbol $b(w,\xi)$ as in \eqref{eq:q_L-def},
$b$ satisfying \eqref{eq:cl-symb-def} with $l$ replaced by $l+r$,
and writing the oscillatory integral
explicitly as a convergent integral, we see that
\begin{equation}
B1(w)=(2\pi)^{-n}\int e^{i(w-w')\cdot\xi}\langle w-w'\rangle^{-2r}
\langle\xi\rangle^{-2s}(1+\Delta_\xi)^s b(w,\xi) (1+\Delta_{w'})^r 1\,dw'
\,d\xi
\end{equation}
for $2r>n$, $2s>n+m$. Changing the variables:
\begin{equation}
B1(w)=(2\pi)^{-n}\int e^{iW\cdot\xi}\langle W\rangle^{-2r}
\langle\xi\rangle^{-2s}(1+\Delta_\xi)^s b(w,\xi)\,dW
\,d\xi.
\end{equation}
This is a convergent integral with $w$ dependence only in $b$. Since
\begin{equation}
b\in\rho_{\infty}^{-m}\rho_{\partial}^{l+r}\Cinf([\Snp;\calC]\times\Snp),
\end{equation}
we conclude that $B1\in x^{l+r}\Cinf([\Snp;\calC])$. Hence, returning to
the global setting, $Au\in x^{l+r}\Cinf([X;\calC])$ as claimed.
\end{proof}

The previous two lemmas show that if $u=e^{i\ft/x}v$, $v\in\Cinf([X;\calC])$,
$A\in\PsiSc^{m,0}(X,\calC)$
then $Au=e^{i\ft/x}v'$ with $v'\in\Cinf([X;\calC])$.
Moreover, $v'$ restricted to the boundary of $[X;\calC]$ only depends on the
restriction of $v$ to $\partial[X;\calC]$. It also follows from the above
proof that if $p\in
\bX$ and $v\in\Cinf([X;\calC])$ vanishes at $\betaSc^{-1}(p)$ then $v'$
also vanishes there, i.e.\ composition is local in $X$ (though not in the
resolved space $[X;\calC]$). Similarly, if $\ft(p)=\ft'(p)$ and
$d_y\ft(p)=d_y\ft'(p)$ (which really just mean that the scattering covectors
$d(\ft/x)$ and $d(\ft'/x)$ agree at $p$) then $e^{-i\ft/x}Ae^{i\ft/x}v$
and $e^{-i\ft'/x}Ae^{i\ft'/x}v$ agree at $p$. This allows us to define
the indicial operators of $A$ at the boundary hypersurfaces of $[X;\calC]$.
These depend on certain choices in general (though the dependence is via
unitary equivalence), but if we have a scattering metric on $X$ they
can be constructed canonically, so we assume this in what follows.

Recall first that a scattering metric $g$ on $X$ is a
metric in the interior of $X$ (smooth symmetric positive definite 2-cotensor)
which is of the form
\begin{equation}
g=\frac{dx^2}{x^4}+\frac{h'}{x^2}
\end{equation}
near $\bX$,
where $x$ is a boundary defining function of $X$ and $h'$ is a smooth symmetric
2-cotensor on $X$ whose restriction to the boundary, $h$, is
positive definite.
Thus, $g$ gives a positive definite pairing on $\Tsc X$, so it is (a somewhat
special) smooth section of $\sct X\otimes\sct X$. We remark that the
choice of such a $g$ fixes $x$ up to the addition of
functions in $x^2\Cinf(X)$.

Next, we recall the definition of the relative scattering
tangent bundle $\Tsc(C;X)$ of a closed embedded submanifold $C$ of $\bX$
from \cite{Vasy:Propagation-2}.

\begin{Def} For a closed embedded submanifold $C$ of $\bX$,
the relative scattering tangent bundle $\Tsc(C;X)$ of $C$ in $X$ is
the subbundle
of $\Tsc_C X$ consisting of
$v\in\Tsc_p X$, $p\in C$, for which there exists
\begin{equation}
V\in\Vsc(X;C)\subset\Vsc(X)
\end{equation}
with $V_p=v$. Here
\begin{equation}
\Vsc(X;C)=x\Vb(X;C)=x\{V\in\Vb(X):\ V\ \text{is tangent to}\ C\}
\end{equation}
and tangency is defined using the (non-injective) inclusion map
$\Tb X\to TX$.
\end{Def}

Thus, in local coordinates $(x,y,z)$ near $p\in C$ such that
$C$ is defined by $x=0$, $y=0$, $\Tsc(C;X)$ is spanned by
$x^2\partial_x$ and $x\partial_{ z_j}$, $j=1,\ldots, m-1$ where
$n-m$ is the codimension on $C$ in $\bX$.
In the case of
Euclidean scattering, $X=\Snp$, $C=\partial\Xb_a$, $g$ the Euclidean
metric, $\Tsc(C;X)$ is naturally isomorphic to $\Tsc_C \Xb_a$, i.e.\ it
should be regarded as the bundle of scattering
tangent vectors of the collision
plane at infinity, spanned by $\partial_{(w_a)_j}$, $j=1,\ldots,m$, $m
=\dim X_a$.

For $C=C_a\in\calC$,
the metric $g$ defines the orthocomplement $(\Tsc(C;X))^\perp$
of $\Tsc (C;X)$ in $\Tsc_C X$.

\begin{Def}\label{Def:rel-sc-cot}
Given $g$, a scattering metric on $X$, the subbundle of $\sct_C X$
consisting of covectors that
annihilate $(\Tsc(C;X))^\perp$,
is denoted by $\sct(C;X)$; we say that it is
the relative scattering cotangent bundle of
$C$ in $X$.
\end{Def}

This bundle of course depends on $g$.
In the case of
Euclidean scattering, $\sct(C;X)$ is naturally isomorphic to $\sct_C \Xb_a$
and is spanned by $d(w_a)_j$, $j=1,\ldots,m$.

We now choose local coordinates $(x,y,z)$ near $p\in C$ such that
$C$ is defined by $x=0$, $y=0$, and such that $x\partial_{y_j}$ give
an orthonormal basis of $(\Tsc (C;X))^\perp$. Note that a basis
of $\Tsc (C;X)$ is given by $x^2\partial_x$ and $x\partial_{z_j}$,
while a basis of $\sct (C;X)$ is given by $x^{-2}\,dx$, $x^{-1}\,dz_j$.
A covector in $\sct X$ can be written in these local coordinates as
\begin{equation}
\tau\frac{dx}{x^2}+\mu\cdot\frac{dy}{x}+\nu\cdot\frac{dz}{x}.
\end{equation}
We will write this as
\begin{equation}
\tau_a\frac{dx}{x^2}+\mu_a\cdot\frac{dy_a}{x}+\nu_a\cdot\frac{dz_a}{x}
\end{equation}
to emphasize the element $C=C_a$ of $\calC$ around which the local
coordinates are centered.
Thus, local coordinates on $\sct_{\bX} X$ are given by $(y,z,\tau,\mu,\nu)$,
while on $\sct (C;X)$ by $(z,\tau,\nu)=(z_a,\tau_a,\nu_a)$.
Note also that
at $C$ the metric function of $h$ is of the form $|\mu|^2+\htil(z,\nu)$
with $|\mu|$ denoting the Euclidean length of $\mu$ and $\htil$ is the
metric function of the restriction of $h$ to $TC$; the metric function
of $g$ (also denoted by $g$) is thus
\begin{equation}\label{eq:g-at-C}
g=\tau^2+\htil+|\mu|^2
\end{equation}
there.

Now if $C=C_a$, $C_b\in\calC$ with $C_a\subset C_b$, we can further adjust
our coordinates so that $C_b$ is defined by $x=0$, $y'=0$, for some
splitting $y=(y',y'')$. With the corresponding splitting of the dual
variable, $\mu=(\mu',\mu'')$, we obtain a well-defined projection
\begin{equation}
\pi_{ba}:\sct_{C_a} (C_b;X)\to\sct (C_a;X),
\end{equation}
\begin{equation}
\pi_{ba}(0,z,\tau,\mu'',\nu)=(z,\tau,\nu).
\end{equation}
In the Euclidean setting this is just the obvious projection
\begin{equation}
\pi_{ba}:\sct_{\partial \Xb_a}\Xb_b\to\sct_{\partial\Xb_a} \Xb_a
\end{equation}
under the inclusion $\Xb_a\subset \Xb_b$.
We write $\pi$ for the collection of these maps.

Before we define the indicial operators, we need to analyze the
structure of the lift of $C_a$ to $[X;\calC]$.
For $C_a\in\calC$ let
\begin{equation}
\calC_a=\{C_b\in\calC:\ C_b\subsetneq C_a\},
\end{equation}
\begin{equation}
\calC^a=\{C_b\in\calC:\ C_a\subsetneq C_b\}.
\end{equation}
We carry out the blow-up $[X;\calC]$ by first blowing up $\calC_a$. Since
all elements of $\calC_a$ are p-submanifolds of $C_a$, the lift
$\beta[X;\calC_a]^*C_a$ of $C_a$ to $[X;\calC_a]$ is naturally
diffeomorphic to
\begin{equation}
\Ct_a=[C_a;\calC_a].
\end{equation}
Thus, over $C'_a$, the regular part of $C_a$, $\Ct_a$
can be identified with $C_a$.
The front face of the new blow-up, i.e.\ of the blow up of
$\beta[X;\calC_a]^*C_a$ in $[X;\calC_a]$ is thus a hemisphere
(i.e.\ ball) bundle
over $\Ct_a$, namely $S^+N\Ct_a$. We write the bundle projection, which
is just the restriction of the new blow-down map to the front face,
$S^+N\Ct_a$ as
\begin{equation}
\rho_a:S^+N\Ct_a\to\Ct_a.
\end{equation}
In the Euclidean
setting, these fibers can be naturally identified with $\Xb^a$ via
the projection $\pi^a$ (extended as in Lemma~\ref{lemma:blowup-1}). Every
remaining blow up in $[X;\calC]$ concerns submanifolds that are either
disjoint from this new front face or are the lift of elements of $\calC^a$.
The former do not affect the structure near the new front face,
$S^+N\Ct_a=\beta[X;\calC_a;C_a]^*C_a$, while
the latter, which are given by the lifts of elements of $\calC^a$,
correspond to blow ups that can be performed in the fibers of $S^+N\Ct_a$.
Note that
the lift of $C_b\in\calC^a$, meets the new front face only at
its boundary since all $C_b$ are subsets of $\bX$.
In particular, the lift $\betaSc^*C_a$
of $C_a$ to $[X;\calC]$ fibers over $\Ct_a$ and the
fibers are diffeomorphic to a hemisphere (i.e.\ ball) with certain
boundary submanifolds blown up. More specifically, the intersection
of $\beta[X;\calC_a;C_a]^*C_b$, $C_b\in\calC^a$, with the
front face $S^+N\Ct_a$ is the image of $T\beta[X;\calC_a]^*C_b$
under the
quotients; $\betaSc^*C_a$ is obtained by blowing these up
in $S^+N\Ct_a$. Hence, the fiber of $\betaSc^*C_a$ over $p\in\Ct_a$
is given by $[S^+N_q C_a;T_q \calC^a]$ where $q=\beta[X;\calC_a](p)\in
C_a$.
In particular, in the Euclidean
setting, the fibers of $\betaSc^*C_a$ over $\Ct_a$
can be naturally identified with
$[\Xb^a;\calC^a]$ via $\pi^a$. 
Thus, we have the following commutative diagrams:
\begin{equation}
\xymatrix{\ar[d]_{\betaSc}\betaSc^*C_a \ar[r]^{\betat_a}
&\Ct_a\ar[dl]^{\beta[C_a;\calC_a]}\\
C_a}
\qquad\qquad
\xymatrix{\ar[d]_{\betat_a}\betaSc^*C_a \ar[r]
&S^+N \Ct_a\ar[dl]^{\rho_a}\\\Ct_a}
\end{equation}
with $\betat_a$ being the fibration to the base $\Ct_a$.

We now define $\sct (\Ct_a;X)$ denote the pull-back of
$\sct(C_a;X)$ by the blow-down
map $\beta[C_a;\calC_a]$:
\begin{equation}
\sct (\Ct_a;X)=\beta[C_a;\calC_a]^*\sct(C_a;X).
\end{equation}
If $C_a\subset C_b$ then $\pi_{ba}$ lifts to a map
\begin{equation}
\pit_{ba}:\sct_{\beta[C_b;\calC_b]^*C_a} (\Ct_b;X)\to\sct(\Ct_a;X).
\end{equation}

We recall from \cite[Section~4]{Vasy:Propagation-2}
that the interior of the fibers $S^+N_p\Ct_a=\rho_a^{-1}(p)$
of
$\rho_a:S^+N \Ct_a\to \Ct_a$, $p\in\Ct_a$, possess a natural
transitive free affine action by the quotient bundle
$(\beta[X;\calC_a]^*_p\Tsc X)/\Tsc_p(\Ct_a;X)$. Thus, the tangent space
of $S^+N_p\Ct_a$ at
every point $q\in\interior(S^+N_p\Ct_a)$ can be naturally identified with
$(\beta[X;\calC_a]^*_p\Tsc X)/\Tsc_p(\Ct_a;X)$, hence with the tangent
space at other
$q'\in\interior(S^+N_p\Ct_a)$.

For each operator $A\in\PsiSc^{m,l}(X,\calC)$,
the $C_a$-indicial operator of $A$,
denoted by $\Ah_{a,l}$, will be a collection of operators, one for each $\zeta
\in\sct_p(\Ct_a;X)$, acting on functions on the fiber $\betat_a^{-1}(p)$
of $\betat_a$. So suppose that $u\in\dCinf(\betat_a^{-1}(p))$; we
need to define $\Ah_a(\zeta)u$. For this purpose
choose $\ft\in\Cinf(X;\Real)$ such that
$d(\ft/x)$, evaluated at $\beta[C_a;\calC_a](p)$,
is equal to $\zeta$. Then let
$\At=e^{-i\ft/x}x^{-l}Ae^{i\ft/x}\in\PsiSc^{m,0}(X,\calC)$,
and choose $u'\in\Cinf([X;\calC])$ such that
$u'|_{\betat_a^{-1}(p)}=u$. Then
\begin{equation}
\Ah_{a,l}(\zeta)u=(\At u')|_{\betat_a^{-1}(p)},
\end{equation}
which is independent of all the choices we made. This can be shown by
an argument which is analogous to the proof of the preceeding lemmas, but
it will also follow from the explicit calculation we make below
leading to \eqref{eq:Psop-46}.
If $l\neq 0$, then $\Ah_{a,l}$ would a priori depend
on the choice of $x$ up to $\calO(x^2)$ terms, but the choice of the
scattering metric $g$ fixes $x$ up to such terms.
We often simplify (and thereby abuse) the
notation and drop the index $l$, i.e.\ we write $\Ah_a=\Ah_{a,l}$, when
the value of $l$ is understood.
Before discussing the $C_a$-indicial operators of
$A\in\PsiSc^{m,l}(X,\calC)$ in detail, we discuss how
we can combine them into a single object.

In the case of Euclidean many-body scattering, $C_a=\partial\Xb_a$ and
$\Ah_{a,l}$ is
a function on $\beta_a^*\sct_{C_a}\Xb_a$
with values in operators on $\Sch(X^a)$; here
\begin{equation}
\beta_a=\beta[C_a;\calC_a]:\Ct_a=[C_a;\calC_a]\to C_a
\end{equation}
is the blow-down map. Note that $\beta_a$ is simply the
restriction of $\beta[\Xb_a;\calC_a]$ to the lift $\Ct_a
=\beta[\Xb_a;\calC_a]^* C_a$. In fact,
\begin{equation}\label{eq:ind-35}
\Ah_{a,l}\in
\Cinf(\beta_a^*\sct_{\partial\Xb_a}\Xb_a,\PsiSc^{m,0}(\Xb^a,\calC^a))
\end{equation}
as we show shortly. Note that if $Z$ is a
(not necessarily compact) manifold with
corners and $(\Xt,\calCt)$ is a many-body space (in \eqref{eq:ind-35}
we take $Z=\beta_a^*\sct_{\partial\Xb_a}\Xb_a$ and $(\Xb^a,\calC^a)$
for the many-body space),
it makes perfectly good sense to talk about
$\Cinf(Z,\PsiSc^{m,l}(\Xt,\calCt))$, i.e.\  about smooth functions on $Z$
with values in $\PsiSc^{m,l}(\Xt,\calCt)$. The topology on
$\PsiSc^{m,l}(\Xt,\calCt)$ is the standard one, namely that of conormal
distributions on $\XXSct$, conormal to $\DeltaSc$, vanishing
to infinite order at $\beta^*\ffb\cup\beta^*\lf\cup\beta^*\rf$,
$\beta:\XXSct\to\XXbt$ the blow-down map. This is equivalent to the
topology arising by localizing operators $A\in\PsiSc^{m,l}(\Xt,\calCt)$
as in \eqref{eq:Psop-17}, and using the topology of the symbol spaces
on the local pieces, i.e.,
with the notation of \eqref{eq:cl-symb-def} and \eqref{eq:PsiScc-5}, of
$\rho_{\infty}^{-m}\rho_{\partial}^l\Cinf([\Snp;\calCt]\times\Snp)$ and
$\bcon^{-m,l}([\Snp;\calCt]\times\Snp)$, in the polyhomogeneous and
non-polyhomogeneous setting respectively (and that of $\dCinf(\Xt
\times\Xt;\Osc_R)$ for the remainder term).

We need to generalize this example to accommodate the geometric setting.
It should be kept in mind throughout following discussion that $Z$ is
simply a `parameter space'.
So suppose first that $\phi:E\to Z$ is a fibration of manifolds with corners
with fiber $\Xt$, a manifold with boundary,
$\calCt_E$ a cleanly intersecting family of p-submanifolds of $E$ which
is fibered over $Z$ with fiber $\calCt$,
a cleanly intersecting family of p-submanifolds of $\partial \Xt$ that gives
rise to a many-body space $(\Xt,\calCt)$.
That is, we suppose that there is an open cover
$\{U_j:\ j\in J\}$ of $Z$ such that $(\phi^{-1}(U_j),\calCt_E\cap
\phi^{-1}(U_j))$ is diffeomorphic to $U_j\times (\Xt,\calCt)$; we denote the
diffeomorphism by $\psi_j$. Let $\partial_\phi E$ denote
the fiber-boundary of $E$, i.e.\ locally it is given by
$U_j\times\partial\Xt$ (under the identification $\psi_j$).
The algebra $\PsiScph^{\infty,-\infty}(E,\calC_E)$ is then defined as the
algebra of operators $A$ acting on, say, functions $u\in\Cinf(E)$
which vanish to infinite order at
$\partial_\phi E$, with the following local characterization.
For each
$U_j$ there is an operator $A'_j\in\Cinf(U_j;\PsiSc^{\infty,-\infty}
(\Xt,\calCt))$ such that for
$u\in\Cinf(E)$ with $\supp u\subset\phi^{-1} (U_j)$ and vanishing
to infinite order at $\partial_\phi E$, $Au=\psi_j^*A'_j(\psi_j^{-1})^*u$.

This local description does not depend on any choices. Indeed, the
local definition is
equivalent to saying that the distribution
kernel $K_A$ of $A$ on the fiber-product
$E\times_Z E$ (with values in scattering densities on the fiber
$\Xt$ from the right factor, to be precise)
is conormal on the appropriate blow-up $E^2_{\Scl,Z}$
of $E\times_Z E$. Here $K_A$ gives rise to the operator $A$ by
fiber-integration
\begin{equation}
Au(w,z)=\int K_A(w,w',z)\,u(w',z)\,|dw'|,
\end{equation}
where $z$ gives coordinates on $Z$, $w$ and $w'$ are variables in the left
and right factor of the fiber $\Xt$ respectively, and we wrote
$K_A=K_A(w,w',z)\,|dw'|$.
Indeed, following the discussion at the beginning
of the previous section, we take $E^2_{\bl, Z}$ to be the blow-up of
$\partial_\phi E\times_Z \partial_\phi E$ in $E\times_Z E$,
$\Delta_{\bl,\phi}$ the
lift of the fiber-diagonal, $\partial_\phi \Delta_{\bl,\phi}$ its
fiber-boundary which we identify with $\partial_\phi E$, $\calCt'_E$
the image of $\calCt_E$ under this identification, and $E^2_{\Scl,Z}$
the blow-up $[E^2_{\bl,Z};\calCt'_E]$. Then the definition of
$\PsiScph^{\infty,-\infty}(E,\calC_E)$
is given by modifying \eqref{eq:PsiSc-def}
the natural way. Since all blowups can be done in the fibers over $Z$
(i.e.\ $Z$ can be regarded as a parameter), this description indeed
agrees with local definition given above.

This intrinsic definition of $\PsiScph^{\infty,-\infty}(E,\calC_E)$
given in the previous paragraph automatically
extends even to the setting where
the fibration $\phi$ is transversal to the collection $\calC_E$, each fiber of
$\phi$ being diffeomorphic to $\Xt$.
Note that in general there
are no diffeomorphisms $\psi_j$ even locally such that image of $\calC_E$ takes
a product form as above, though such diffeomorphisms exist, for
example, if $\calC_E$ is locally linearizable.
In particular,
we can take $Z=\sct(\Ct_a;X)$, $E$ to be the pull-back of $Z$ to
$S^+N\Ct_a$ by $\rho_a$,
$\phi:E\to Z$ the map $\rho_a^\sharp$ induced by the pull-back,
\begin{equation}
E=\rho_a^*\sct(\Ct_a;X),\quad \rho_a^\sharp:E\to\sct(\Ct_a;X).
\end{equation}
Thus, $E$ is a vector bundle over $S^+N\Ct_a$ with projection $\pi$.
Finally, we let $\calC_E$ consist of the inverse images under $\pi$
of the lifts of $C_b\in\calC^a$
to $[X;\calC_a;C_a]$ intersected with the new front face, $S^+N\Ct_a$;
in fact, we also add $\partial_\phi E$ to $\calC_E$ to play the role of
$C_0$ in $\calC$.
We are then in the setting discussed above, so we have defined
\begin{equation}
\PsiScra^{\infty,-\infty}(\rho_a^*\sct(\Ct_a;X),\calCt_a),
\quad \calCt_a=\pi^{-1}(S^+N\Ct_a\cap\beta[X;\calC_a;C_a]^*
\calC^a)\cup\{\partial_\phi E\}.
\end{equation}
Recall that for $C_b\in\calC^a$,
\begin{equation}
S^+N\Ct_a\cap\beta[X;\calC_a;C_a]^*
C_b=T\beta[X;\calC_a]^*C_b,
\end{equation}
the right hand side understood as the image of the tangent space under the
quotient map. We are now ready to prove the following proposition.

\begin{prop}
Suppose that $A\in\PsiSc^{r,l}(X,\calC)$. Then the
indicial operators of $A$ satisfy
\begin{equation}
\Ah_{a,l}\in\PsiScra^{r,0}(\rho_a^*\sct(\Ct_a;X),\calCt_a).
\end{equation}
\end{prop}

\begin{proof}
We prove this statement
by finding $\Ah_a(\zeta)$ explicitly in terms of local
coordinates. To simplify the notation we assume that $A\in\PsiSc^{r,0}
(X,\calC)$. We identify $X$ with $\Snp$ locally so that $C_a$ is
given by $x=0$, $y=0$.
In the interior of $\betaSc^* C_a$ we can use the
same coordinates as at the front face of $[X;C_a]$, i.e.\ the ones
given in \eqref{eq:ff-coords}-\eqref{eq:ff-coords-2}. So suppose that
$u'$ is supported in the region of validity of these coordinates.
Then
\begin{equation}\label{eq:Psop-40}
Au'(w)=\int K(w,w') u'(w')\,dw'=\int \at(w,W) u'(w-W)\,dW
\end{equation}
with $\at$ as in \eqref{eq:PsiSc-23}. Here the integral is understood
as a distributional pairing in general, but it actually converges
if $r<-n$.
We now consider the coordinates \eqref{eq:ff-coords} on the both factors,
i.e.\ we take
$(x',Y',z')$ corresponding to $w'=w-W$, and $(x,Y,z)$ corresponding to $w$.
Expressing $(x',Y',z')$ in terms of $(x,Y,z)$ and $W$ (using $w'=w-W$)
gives
\begin{equation}
x'=x(1-x(W_a)_m)^{-1},\ z'_j=\frac{z_j-x(W_a)_j}{1-x(W_a)_m},\ Y'_j=
Y_j-(W^a)_j,
\end{equation}
where we wrote $W=(W_a,W^a)$ and $(W_a)_j$, $(W^a)_j$ denote the components
of $W_a$ and $W^a$ respectively.
Thus, \eqref{eq:Psop-40} yields
\begin{equation}
Au'(x,Y,z)=\int \at(x,Y,z,W) u'\left(\frac{x}{1-x(W_a)_m}, Y-W^a,
\frac{z_j-x(W_a)_j}{1-x(W_a)_m}\right)\,dW.
\end{equation}
Evaluating at $x=0$ gives
\begin{equation}\begin{split}
Au'(0,Y,z)&=\int \at(0,Y,z,W) u'(0, Y-W^a,z)\,dW\\
&=\int \left(\int \at(0,Y,z,W)\,dW_a\right) u'(0, Y-W^a,z)\,dW^a.
\end{split}\end{equation}
Since $\at$ is the inverse Fourier transform in the $\xi$ variable
of the symbol $a$ whose left quantization is $A$, and since the $W_a$
integral above can be understood as the Fourier transform in $W_a$
evaluated at the origin, we deduce that
\begin{equation}
Au'(0,Y,z)=(2\pi)^{-(n-m)}
\int e^{iW^a\cdot\xi^a}a(0,Y,z,0,\xi^a) u'(0, Y-W^a,z)\,d\xi^a\,dW^a.
\end{equation}
Thus, the indicial operator $\Ah_a((p,0))$ where $(p,0)\in\sct(\tilde C_a;X)$
is the zero covector above $p=(0,0,z)\in C_a$ is given by
\begin{equation}
\Ah_a((p,0))u(Y)=(2\pi)^{-(n-m)}
\int e^{iW^a\cdot\xi^a}a(0,Y,z,0,\xi^a) u(Y-W^a)\,d\xi^a\,dW^a,
\end{equation}
i.e.\ by the left quantization in $(Y,\xi^a)=(W^a,\xi^a)$
of $a(0,Y,z,0,\xi^a)$. Similar results hold for $\Ah_a(\zeta)$ in general,
namely
\begin{equation}\label{eq:Psop-46}
\Ah_a(z,\xi_a)u(Y)=(2\pi)^{-(n-m)}
\int e^{iW^a\cdot\xi^a}a(0,Y,z,\xi_a,\xi^a) u(Y-W^a)\,d\xi^a\,dW^a,
\end{equation}
Though the local coordinates are only valid in the interior of $\betaSc^*
C_a$, hence not at $\betat_a^*\partial\Ct_a$, the continuity
of $\At u$ up to $\betat_a^*\partial\Ct_a$ shows that \eqref{eq:Psop-46}
also holds with $p\in\Ct_a$.

The explicit expression, \eqref{eq:Psop-46} shows,
in particular, that $\Ah_a(\zeta)u$ is indeed independent of the extension
$u'$ of $u$ that we chose, and also of the choice of $\ft$ with
$d(\ft/x)$ prescribed at $\betaSc(p)$. Moreover, also from
\eqref{eq:Psop-46},
for each $\zeta\in\sct_p(\tilde C_a;X)$, $p\in\Ct_a$,
\begin{equation}
\Ah_a(\zeta)\in\PsiSc^{r,0}(\rho_a^{-1}(p),T_p\calC^a);
\end{equation}
here we wrote $T_p \calC^a$ for $T_p \beta[X;\calC_a]^*\calC^a$
for simplicity.
In fact, \eqref{eq:Psop-46} shows the more precise statement
which encodes the smooth dependence of $\Ah_a(\zeta)$ on $\zeta$,
namely that
\begin{equation}
\Ah_{a,l}\in\PsiScra^{r,0}(\rho_a^*\sct(\Ct_a;X),\calCt_a).
\end{equation}
In the Euclidean setting the many-body space
$(\rho_a^{-1}(p),T_p \calC^a)$
can be identified with $(\Xb^a,\calC^a)$, and we can write
\begin{equation}
\Ah_a(\zeta)\in\PsiSc^{r,0}(\Xb^a,\calC^a),
\end{equation}
and correspondingly
\begin{equation}
\Ah_{a,l}\in
\Cinf(\beta_a^*\sct_{\partial\Xb_a}\Xb_a,\PsiSc^{r,0}(\Xb^a,\calC^a))
\end{equation}
as we have claimed.
\end{proof}

If $A\in\PsiSc^{r,0}(X,\calC)$, then
the vanishing of $\Ah_{a,0}(\zeta)$ for every $a$ and every
$\zeta\in\sct (\tilde C_a;X)$ implies, by our explicit formula, that
$a\in\Cinf([X;\calC]\times\Snp)$ vanishes at $(\partial[X;\calC])
\times\Snp$, so $A\in\PsiSc^{r,1}(X,\calC)$. Thus, the vanishing of
$\prs_{\Scl,r}(A)$ and all indicial operators together, for $A\in
\PsiSc^{r,0}(X,\calC)$, say, implies that $A\in\PsiSc^{r-1,1}(X,\calC)$.

An advantage of the oscillatory testing
definition of the indicial operators is that it makes their multiplicative
property clear.

\begin{prop}
If $A\in\PsiSc^{m,l}(X,\calC)$, $B\in\PsiSc^{m',l'}(X,\calC)$ then
\begin{equation}
\widehat{AB}_{a,l+l'}(\zeta)u=\Ah_{a,l}(\zeta)\Bh_{a,l'}(\zeta) u.
\end{equation}
\end{prop}

The indicial operators are related via the projections $\pit_{ba}$. Thus,
if $\zeta\in\sct (\Ct_a;X)$, then the indicial operators
of $\Ah_{a,l}(\zeta)$ are $\Ah_{b,l}(\zetat)$ where $C_a\subset C_b$, $C_a\neq
C_b$, and $\zetat\in\sct_{\beta[C_b;
\calC_b]^*C_a}(\Ct_b;X)$ is such that $\pit_{ba}(\zetat)=\zeta$. This
follows easily from the explicit coordinate form of the indicial
operators. In particular,
if $A\in\PsiSc^{m,l}(X,\calC)$ with $m<0$, and if
$\Ah_{b,l}(\zetat)$ vanishes for all such $b$ and $\zetat$,
then $\Ah_{a,l}(\zeta)$
is compact. We thus have the following proposition.

\begin{prop}
If $A\in\PsiSc^{r,0}(X,\calC)$ is such that
$\prs_{\Scl,r}(A)$ never vanishes and
$\Ah_a(\zeta)$ is invertible in
$\PsiSc^{r,0}(\rho_a^{-1}(p),T_p \calC^a)$
(i.e.\ in $\PsiSc^{r,0}(\Xb^a,\calC^a)$ in the Euclidean setting) for every $a$
and for every
$\zeta\in\sct (\tilde C_a; X)$, then there exists a parametrix $P$ for
$A$ such that $PA-\Id, AP-\Id\in\PsiSc^{-\infty,\infty}(X,\calC)$.
\end{prop}

\begin{proof}
Choose $P_0\in\PsiSc^{-r,0}(X,\calC)$
with principal symbol $\prs_{\Scl,r}(A)^{-1}$
and indicial operators $\Ah_{a,0}(\zeta)^{-1}$.
Note that these match
up under the projections $\pit_{ba}$ since those of $A$ match up.
Thus, they indeed specify the restriction of some $p_0\in\rho_\infty^r
\Cinf([X;\calC]\times\Snp)$ to the boundary of $[X;\calC]\times\Snp$.
Quantizing $p_0$ gives an operator $P_0$ with the required indicial
operators and principal symbol. Hence, $E=\Id-P_0A\in\PsiSc^{0,0}(X,\calC)$
has vanishing
principal symbol and indicial operators, so it is in $\PsiSc^{-1,1}(X,\calC)$.
Summing the Neumann series $\sum_{j=1}^\infty E^j$ asymptotically to some
$F\in\PsiSc^{-1,1}(X,\calC)$ and letting $P=(\Id+F)P_0$ gives the
required left parametrix. A right parametrix can be constructed similarly,
and then the usual argument shows that they can be taken to be the same.
\end{proof}

For $A\in\PsiSc^{m,0}(X,\calC)$ self-adjoint, $m>0$, with $\prs_{\Scl,m}(A)$
never vanishing, we automatically have
that $(A-\lambda)^{-1}\in\PsiSc^{-m,0}(X,\calC)$ for $\lambda\in\Cx\setminus
\Real$. Moreover, the blow-up of $(A-\lambda)^{-1}$ in $\PsiSc^{-m,0}
(X,\calC)$ can be analyzed uniformly as $\lambda$ approaches the real axis,
see e.g.\ \cite{Hassell-Vasy:Symbolic, Vasy:Propagation-2}.
Therefore,
the functional calculus for self-adjoint operators $A$ and the Cauchy
integral representation of $\phi(A)$ via almost analytic extensions,
as in the work of Helffer and Sj\"ostrand \cite{Helffer-Sjostrand:Schrodinger},
Derezi\'nski and G\'erard \cite{Derezinski-Gerard:Scattering},
see also \cite{Hassell-Vasy:Symbolic},
gives immediately

\begin{prop}
Suppose that $A\in\PsiSc^{m,0}(X,\calC)$ self-adjoint, $m>0$,
and $\prs_{\Scl,m}(A)$
never vanishes. Suppose also that $\phi\in\Cinf_c(\Real)$. Then
$\phi(A)\in\PsiSc^{-\infty,0}(X,\calC)$ and its indicial operators are
$\phi(\Ah_a(\zeta))$. If instead we assume $\phi\in S^{-r}_\phg(\Real)$
then $\phi(A)\in\PsiSc^{-rm,0}(X,\calC)$.
\end{prop}

If $m=0$, that is $A\in\PsiSc^{0,0}(X,\calC)$, then
$\phi(A)\in\PsiSc^{0,0}(X,\calC)$ without any assumption 
on the invertibility of
$\prs_{\Scl,0}(A)$.
We thus have:

\begin{prop}\label{prop:fc-3}
Suppose that $A\in\PsiSc^{0,0}(X,\calC)$ self-adjoint. If $\phi\in
\Cinf(\Real)$ then $\phi(A)\in\PsiSc^{0,0}(X,\calC)$.
\end{prop}

\begin{proof}
Since $A$ is bounded, we can replace $\phi$ by a function $\psi\in\Cinf_c
(\Real)$ such that $\phi\equiv\psi$ on the spectrum of $A$. Now
$\prs_{\Scl,0}(A-\lambda)=\prs_{\Scl,0}(A)-\lambda$ is invertible for
$\lambda\in\Cx\setminus\Real$, so $(A-\lambda)^{-1}\in
\PsiSc^{0,0}(X,\calC)$ for $\lambda\nin\Real$. Again, $(A-\lambda)^{-1}$
can be analyzed uniformly up to the real axis, and then
the
Cauchy integral representation of $\psi(A)$ now proves the
proposition.
\end{proof}

\section{The wave front set}\label{sec:WF}
The $\Scl$-wave front set $\WFSc(u)$
of a distribution $u$, and the $\Scl$-operator
wave front set $\WFScp(A)$ of $A\in\PsiSc^{m,l}(X,\calC)$,
at infinity will be defined as subsets of the compressed scattering
cotangent bundle
\begin{equation}
\scdt X=\cup_a\sct_{C'_a} (C_a;X);
\end{equation}
we have defined $\sct (C_a;X)$ in Definition~\ref{Def:rel-sc-cot}.
This is very similar to the image of the cotangent bundle in
the compressed cotangent bundle (the b-cotangent bundle) that
Melrose and Sj\"ostrand used to describe the propagation of singularities
for the wave equation in domains with smooth boundaries
\cite{Melrose-Sjostrand:I} and also to
the corresponding phase space for domains with corners $\Omega$,
$\dot T^*_\bl\Omega$, which was the setting for Lebeau's analysis
of the singularities of solutions to the wave equation on $\Omega$.
We topologize $\scdt X$ using the projection $\pi:\sct_{\bX} X\to\scdt X$.
Thus, $\WFSc(u)$ will contain less detailed information than $\WFsc(u)$,
the latter being a subset of $\sct_{\bX} X$. However, neither of these
wave front sets can be used to describe the other; they are simply
different. The fact that $\WFSc(u)$ contains fewer details about $u$
simply corresponds to the singular behavior of elements of $\PsiSc^{m,l}
(X,\calC)$, compared to those of $\Psisc^{m,l}(X)$.

The definition of $\WFSc(u)$ and $\WFScp(A)$ will be local in $X$. Thus,
we can always work on $\Snp$ instead. Just like when we defined
$\PsiSc^{m,l}(X,\calC)$, we will be able to proceed either by giving
an explicit description in $\Snp$ via the Fourier transform, or by
giving invariant definitions.

We start with the operator wave front sets.
The invariant definition proceeds by oscillatory testing.

\begin{Def}
Suppose that $A\in\PsiSc^{m,l}(X,\calC)$ and
$\zeta\in\sct_p(C_a;X)$, $p\in C'_a$. We say that
$\zeta\nin\WFScp(A)$
if and only if
there exist a neighborhood $U$ of
$\zeta$ in $\scdt X$ and a neighborhood
$V$ of $p$ in $X$ such that
$Au\in\dCinf(X)$ for every oscillatory
function $u=e^{if/x}v$, $v\in\Cinf([X;\calC])$ with
$\pi(\Graph(d(f/x)))\cap\scdt_{V\cap\bX} X\subset U$ and
$\supp v\subset \beta_\Scl^{-1}(V)$.
\end{Def}

This definition
implies immediately that $\WFScp(A)$ is closed in $\scdt X$,
\begin{equation}
\WFScp(A+B)\subset\WFScp(A)\cup\WFScp(B),\quad A,B\in\PsiSc(X,\calC),
\end{equation}
\begin{equation}
\WFScp(AB)\subset\WFScp(A)\cap\WFScp(B),\quad A,B\in\PsiSc(X,\calC).
\end{equation}

We can also formulate the definition explicitly.
We thus locally identify $X$ with $\Snp$ and consider
$A\in\PsiSc^{m,l}(\Snp,\calC)$. We also
identify $\sct\Snp$ with
$\Snp\times\Rn$.
So suppose that $A$ is the left
quantization of a symbol $a\in\rho_\infty^{-m}\rho_\partial^l\Cinf
([\Snp;\calC]\times\Snp)$. Then  $\zeta\nin\WFScp(A)$,
$\zeta\in\sct_p(C_a;X)$, $p\in C'_a$, if
and only if
there exists a neighborhood $U$ of $\zeta$ in $\scdt\Snp$ such that
$a$ vanishes at $U'\subset(\partial[\Snp;\calC])\times\Rn$ to infinite
order where $U'$ is the
inverse image of $U$
under the composite map
\begin{equation}\begin{CD}\label{eq:WFScp-20}
(\partial[\Snp;\calC])\times\Rn@>\betaSc\times\identity>>
(\partial\Snp)\times\Rn=\sct_{\Sn}\Snp
@>\pi>>\scdt\Snp.
\end{CD}\end{equation}
Note that this condition implies automatically that $a$ vanishes
to infinite order on the closure of $U'$ in
$(\partial[\Snp;\calC])\times\Snp$, in particular, it is rapidly
decreasing in some directions as $|\xi|\to\infty$.
It follows immediately from the usual formulae relating quantizations
and the effect of diffeomorphisms
that this definition is independent
of such choices. For example, we could have equally well written
$A$ as the right quantization of a symbol with similar properties.

In general, the structure of $\pi$ can be somewhat complicated
in explicit coordinates. However,
in the actual Euclidean setting, i.e.\ where $\calC$ arises
from a family of linear subspaces $\calX$, 
or indeed,
if we merely assume that $(X,\calC)$ is locally linearizable,
it is particularly easy to give a sufficient condition for
$\zeta\nin\WFScp(A)$.
Namely, if there is a neighborhood
$V$ of $\zeta=(0,z^0_a,\xi^0_a)\in\sct_p(C_a;\Snp)$, $p\in C'_a$,
in $\partial\Snp\times\Rm$ such that
$a$ vanishes to infinite order at every point $(q',\xi)\in
(\partial[\Snp;\calC])\times\Rn$ with $(\betaSc(q'),\xi_a)\in V$,
then $\zeta\nin\WFScp(A)$. Note that as $p\in C'_a$, we can always
assume, by reducing the size of $V$ if necessary, that $(q,\xi^a)\in V$
implies $q\in C'_b$ for some $b$ with $C_a\subset C_b$.
We can see that this condition is sufficient for $\zeta\nin\WFScp(A)$
since for nearby $q\in\Sn$, assuming as we may that
$q\in C'_b$, $C_a\subset C_b$, the restriction of $\pi$ to
$\sct_q\Snp$ takes the form $(q,\xi_b,\xi^b)\mapsto(q,\xi_b)$ and
$\xi_b$ splits as $(\xi'_b,\xi''_b)$ with $\xi'_b=\xi_a$. Thus,
the condition of the previous paragraph holds if we take
\begin{equation}
U=\cup_b\{(q,\xi_b):\ q\in C'_b,\ \exists\xi_a,\xi_b''\ \text{s.t.}
\ (q,\xi_a)\in V\Mand\xi_b=(\xi'_b,\xi''_b)\}.
\end{equation}

The definition of the wave front set of a distribution $u\in\dist(X)$ at
$\bX$ is more complicated. To determine whether $\zeta\in\sct_p(C_a;X)$,
$p\in C'_a$ is in $\WFSc(u)$, we would like to cut off $u$ to be supported
near $p$, i.e.\ consider $\psi u$, $\psi\in\Cinf(X)$, $\psi\equiv 1$
near $p$, identify a neighborhood of $p$ with an open set in $\Snp$
near $\partial\Snp$, and consider smoothness of
the Fourier transform of $u$,
$\Fr \psi u$. Indeed, in the two-body setting, hence in the many-body
setting if we consider
$\zeta\in \sct_{C'_0}(C_0;X)=\sct_{C'_0}X$,
written as a covector $\xi\cdot dw$ over $p\in C'_0$, we
have
\begin{equation}\label{eq:WFsc-def-1}
\zeta\nin\WFsc(u)\Miff\exists\psi\ \text{as above, s.t.}\ \Fr\psi u
\ \text{is smooth near}\ \xi.
\end{equation}
In the general many-body setting,
$\zeta\in\sct_p(C_a;X)$, $p\in C'_a$, $\zeta$
takes the form $\xi_a\cdot dw_a$, and correspondingly we would like to say
that $\Fr\psi u$ is Schwartz in a region including the subspace $S$
consisting of all points of the form
$(\xi_a,\xi^a)$ where $\xi^a$ is arbitrary. Here Schwartz takes the
place of smooth since the region is not compact in $\Rn$. However, as
shown by the example of ordinary wave front set, we
cannot expect that this wave front set behaves reasonably unless the
region $U$ is conic near infinity, i.e.\ unless it is a neighborhood of
the closure of $S$ in the radial compactification $\Snp$ of $\Rn$. This
however introduces the complication that all parallel translates of
$S$ intersect $U$, and we are exactly interested in separating from
each other the
singularities on the various translates of $S$. This problem is not
too serious, especially for generalized eigenfunctions of many-body
Hamiltonians $H$, but it introduces additional terms into the
following definition which is modelled on that of the fibred cusp wave front
set by Mazzeo and Melrose \cite{Mazzeo-Melrose:Fibred}.

\begin{Def}\label{Def:WFSc}
We say that
\begin{equation}\begin{split}
\zeta\nin\WFSc(u)\cap\sct_{C'_a} (C_a;X)\Miff
&\exists A\in\PsiSc^{0,0}(X,\calC),
\ \Ah_{a,0}(\zeta)\ \text{invertible in}\ \PsiSc^{0,0}(\Xb^a,\calC^a),\\
&\exists B_j\in\PsiSc^{-\infty,0}(X,\calC),\ \zeta\nin\WFScp(B_j),\\
&\exists u_j\in\dist(X),\ j=1,\ldots,s,\ f\in\dCinf(X),\\
&Au=\sum_{j=1}^s B_j u_j+f.
\end{split}\end{equation}
Here we used the Euclidean notation $\PsiSc^{0,0}(\Xb^a,\calC^a)$
instead of
$\PsiSc^{0,0}(\rho_a^{-1}(p),T_p\calC^a)$
for the sake
of simplicity.
Similarly, the filtered version of the $\Scl$-wave front set is given by
\begin{equation}\begin{split}
\zeta\nin\WFSc^{m,l}(u)\cap\sct_{C'_a} (C_a;X)\Miff
&\exists A\in\PsiSc^{0,0}(X,\calC),
\ \Ah_{a,0}(\zeta)\ \text{invertible in}\ \PsiSc^{0,0}(\Xb^a,\calC^a),\\
&\exists B_j\in\PsiSc^{-\infty,0}(X,\calC),\ \zeta\nin\WFScp(B_j),\\
&\exists u_j\in\dist(X),\ j=1,\ldots,s,\ f\in\Hsc^{m,l}(X),\\
&Au=\sum_{j=1}^s B_j u_j+f.
\end{split}\end{equation}
\end{Def}

Thus, if $p\in C'_a$,
then the part of
$\WFSc$ over $p$ lives in $\sct_p (C_a;X)$. If we define the
scattering wave front set, $\WFsc(u)$, in terms of operators
instead of the description of $\WFsc(u)$ given in \eqref{eq:WFsc-def-1}
then the extra terms $B_j u_j$ can be dropped. In fact, \eqref{eq:WFsc-def-1}
is equivalent to requiring that $Au\in\dCinf(\Snp)$ where $A=\Frinv
\phi\Fr\psi\in\Psisc^{-\infty,0}(\Snp)$,
$\psi$ as above, and $\phi\in\Cinf_c(\Real)$ is identically
$1$ near $\xi$. The additional terms $B_j u_j$ for $\WFSc(u)$ thus arise
because the invertibility of $\Ah_a(\zeta)$ implies that
$\prs_{\Scl,0}(\Ah_a(\zeta))$
cannot vanish which in turn means that $\prs_{\Scl,0}(\Ah_a
(\zeta'))$ is non-zero for every $\zeta'\in\sct_p(C_a;X)$ since
$\prs_{\Scl,0}(\Ah_a(\zeta))=\prs_{\Scl,0}(\Ah_a
(\zeta'))$. This simply corresponds
to the conic cutoff requirement discussed before the definition.

With the topology we put on $\scdt X$,
$\WFSc(u)$ is closed due to the relationship between the
indicial operators mentioned above. Namely, the invertibility of
$\Ah_{a,0}(\zeta)$ implies that of $\Ah_{b,0}(\zetat)$ with $\pit_{ba}(\zetat)=\zeta$,
hence of $\Ah_{b,0}(\zetat')$ for nearby $\zetat'$. As the complement
of $\WFScp(B_j)$
is open, this implies that the complement of
$\WFSc(u)$ is also open. In addition,
\begin{equation}
\WFSc(u_1+u_2)\subset\WFSc(u_1)\cup\WFSc(u_2)
\end{equation}
and the corresponding result also holds for the filtered wave front set.
Moreover,
\begin{equation}\label{eq:WF-Au-1}
A\in\PsiSc^{m,l}(X,\calC),\ u\in\dist(X)\Rightarrow \WFSc(Au)\subset
\WFScp(A)\cap\WFSc(u),
\end{equation}
and similarly
\begin{equation}\label{eq:WF-11}
A\in\PsiSc^{m,l}(X,\calC),
\ u\in\dist(X)\Rightarrow \WFSc^{m'-m,l+l'}(Au)\subset
\WFScp(A)\cap\WFSc^{m',l'}(u).
\end{equation}
We refer to \cite{Mazzeo-Melrose:Fibred} for detailed arguments; we only
need simple modifications of their proofs. We also refer to the remarks after
Proposition~\ref{prop:WF-3} for connecting this wave front set to the
one discussed in \cite{Vasy:Propagation-2} in three-body scattering.

This wave front set, $\WFSc$, gives a complete microlocal
description of distributions at $\bX$. To state it generally, we would
need to define the extension of the standard wave front set of $u$ to
give a subset of $\ScS[X;\calC]$, but for us the following extension of
\eqref{eq:WF-Au-1} suffices:
\begin{equation}\label{eq:WF-12}
A\in\PsiSc^{-\infty,l}(X,\calC),
\ u\in\dist(X),\ \WFScp(A)\cap\WFSc(u)=\emptyset\Rightarrow Au\in\dCinf(X).
\end{equation}
We remark that in \cite{RBMSpec}, $\WFsc(u)$ (or rather its part over
$\bX$) is defined as a subset of $\scct_{\bX}X$, the radial compactification
of $\sct_{\bX}X$ in the fibers. The part at fiber-infinity, i.e.\ at
the boundary arising from the radial compactification
of the fibers, extends the usual wave front set from the interior. However,
for us this extension is not important; the operator wave front set
of nearly all operators we are
interested in is contained in a compact region of $\scdt X$.

The description of the wave front set becomes simpler for generalized
eigenfunctions of many-body Hamiltonians $H$. Namely, we have the following
result.

\begin{prop}\label{prop:WF-3}
Suppose that $u\in\dist(X)$, $H\in\PsiSc^{m,0}(X,\calC)$, $m>0$ is self-adjoint
and $\prs_{\Scl,m}(H)$ never vanishes. Let $\lambda\in\Real$, and define
$W\subset\scdt X$ by
\begin{equation}\begin{split}
\zeta\nin W\cap\sct_{C'_a}(C_a;X)
\Miff&\exists\psi\in\Cinf_c(\Real),\ \psi(\lambda)=1,\\
&\exists A\in\PsiSc^{-\infty,0}(X,\calC),\ \Ah_a(\zeta)=\widehat{\psi(H)}_a,
Au\in\dCinf(X).
\end{split}\end{equation}
Then
\begin{equation}
\WFSc(u)\subset\WFSc((H-\lambda)u)\cup W.
\end{equation}
The same conclusion holds with $\WFSc$ replaced by $\WFSc^{m,l}$ and
$Au\in\dCinf(X)$ by $Au\in\Hsc^{m,l}(X)$.
\end{prop}

\begin{proof}
Suppose that $\zeta\nin\WFSc((H-\lambda)u)$ and $\zeta\nin W$.
With $\psi$ as above, let $\psit(t)=(1-\psi(t))/(t-\lambda)$, so
$\psit\in S^{-1}_\phg(\Real)$ as $\psi(\lambda)=1$. Then
$\psit(H)\in\PsiSc^{-m,0}(X,\calC)$ and $\Id=\psit(H)(H-\lambda)+\psi(H)$.
With $A$ as above, let $A'=A+(\Id-\psi(H))\in\PsiSc^{0,0}(X,\calC)$.
Then $\Ah'_a(\zeta)=\Id$ and
\begin{equation}
A'u=Au+\psit(H)(H-\lambda)u.
\end{equation}
But $Au\in\dCinf(X)$ by assumption, so by \eqref{eq:WF-Au-1}
\begin{equation}
\WFSc(A'u)=\WFSc(\psit(H)(H-\lambda)u)\subset\WFSc((H-\lambda)u).
\end{equation}
Hence, there exist $A''$ (in place of $A$), $B_j$, etc., as in
Definition~\ref{Def:WFSc}, $A''A'u=f+\sum B_j u_j$, and the
indicial operator of $A''A'$ at $\zeta$
is just the composite of those of $A''$ and
$A'$, hence invertible,
showing that $\zeta\nin\WFSc(u)$.
\end{proof}

\begin{rem}
Our definition of $\WFSc(u)$, which is in particular valid if $(X,\calC)$ is
a three-body space, is {\em different} from the wave front set
$\WFtSc(u)$ used in
\cite{Vasy:Propagation-2} in the three-body setting. Indeed,
in the definition of $\WFtSc(u)$, the terms $B_j u_j$
appearing in Definition~\ref{Def:WFSc} were not allowed. Consequently,
\eqref{eq:WF-Au-1}, more precisely $\WFtSc(Au)\subset\WFtSc(u)$, and
its filtered analogue did not hold in general. However, for the
positive commutator proofs of both \cite{Vasy:Propagation-2} and the
present paper, one only needs \eqref{eq:WF-12}, which was proved
for $\WFtSc$ under the weak additional condition that $\WFtSc'(A)$
is compact; this condition holds for all operators appearing
in such estimates in both papers.

Note that $\WFSc(u)\subset\WFtSc(u)$ directly from the definition.
Moreover, if $(\Id-P)u\in\dCinf(X)$ for some $P\in\PsiSc^{-\infty,0}(X,\calC)$
(e.g.\ $P=\psi(H)$ in the setting of the proposition)
then $\WFSc(u)=\WFtSc(u)$. In fact, suppose that
$\zeta\nin\WFSc(u)$, so $Au=\sum B_j u_j+f$ as in Definition~\ref{Def:WFSc}.
Since $\Ah$ is invertible near $\zeta$, we can arrange (by inverting
$\Ah$ nearby, i.e.\ by constructing a `microlocal parametrix')
that $A'u=\sum B'_j u_j+f'$ with $A'$ and $B'_j$ as in the definition, but
in addition with $\zeta\nin\WFScp(\Id-A')$
(cf.\ \cite{Mazzeo-Melrose:Fibred}).
Using the methods of Section~\ref{sec:commutators},
given any neighborhood $U$ of $\zeta$,
it is easy to construct an operator $C\in\PsiSc^{-\infty,0}(X,\calC)$
such that $\WFScp(C)\subset U$
and $\zeta\nin\WFScp(P-C)$ (hence the same holds for a neighborhood
of $\zeta$). Since the indicial operator of $Q=C+(\Id-P)$ at $\zeta$ is
the identity, and since $(\Id-P)u\in\dCinf(X)$, we only need
to prove that $Cu\in\dCinf(X)$ to conclude that $\zeta\nin\WFtSc(u)$.
But $Cu=\sum CB'_j u_j+Cf'+C(\Id-A')u$, so if $U$ is chosen sufficiently
small, then $CB'_j,\ C(\Id-A')\in\PsiSc^{-\infty,\infty}(X,\calC)$,
so $Cu\in\dCinf(X)$ indeed.
\end{rem}

\section{The Hamiltonian and generalized broken bicharacteristics}
\label{sec:Hamiltonian}
We next analyze the operator $H-\lambda$ where $H=\Delta+V$ and $\Delta$
is the Laplacian of a scattering metric
\begin{equation}
g=\frac{dx^2}{x^4}+\frac{h'}{x^2}.
\end{equation}
Recall that $h'$ is a smooth symmetric
2-cotensor on $X$ whose restriction to $\bX$ (i.e.\ its pull-back), $h$, is
positive definite. We
assume that
\begin{equation}\label{eq:mb-pot}
V\in\Cinf([X;\calC];\Real)\ \text{vanishes at}\ \betaSc^* C_0,
\end{equation}
i.e.\ $V$ vanishes in the free region. This implies that
\begin{equation}
H\in\DiffSc^2(X,\calC).
\end{equation}
Such a situation arises, for example,
in actual Euclidean scattering if the potentials $V_a$ (in the notation of
the introduction) are classical symbols of order $-1$ on $X^a$. Hence,
we make the following definition.

\begin{Def}
A many-body Hamiltonian is an operator $H=\Delta+V$ where $\Delta$ is
the Laplacian of a scattering metric $g$, and $V$ satisfies
\eqref{eq:mb-pot}.
\end{Def}

As indicated in the Introduction,
from this point on we also make the assumption
\begin{equation}\label{eq:hypo-loc-lin}
(X,\calC)\ \text{is locally linearizable};
\end{equation}
this will simplify the analysis. We recall that this is
equivalent to the local existence of Riemannian
metrics on $\bX$, possibly different from $h$, with respect to which all
elements of $\calC$ are totally geodesic.

Since $\prs_{\scl,2}(\Delta)$ never vanishes, the same holds for
$\prs_{\Scl,2}(H)$ which is the pull-back of the former.
A simple calculation, see \cite[Sections~4 and 11]{Vasy:Propagation-2}
for more details, shows that
the indicial operators of $H$ are given by
\begin{equation}\label{eq:ham-4}
\Hh_{a,0}(\xi)=\Hh_{a,0}((p,0))+
\tau^2+\htil(z,\nu),\ \xi=(z,\tau,\nu)\in\sct(\Ct_a;X),
\end{equation}
\begin{equation}\label{eq:ham-5}
\Hh_{a,0}(p,0)=\Delta_Y+V(p,Y)
\end{equation}
where $Y$ are `Euclidean
coordinates' on the interior of $\rho_a^{-1}(p)$, i.e.\ that of
$\betat_a^{-1}(p)$, and $\Delta_Y$ is the
Euclidean Laplacian.

More precisely, we have seen in Section~\ref{sec:indicial} that
$(\beta[X;\calC_a]^*_p\Tsc X)/\Tsc_p(\Ct_a;X)$ naturally acts transitively
and freely on
the interior of
$\rho_a^{-1}(p)=S^+N_p\Ct_a$, so it makes sense to talk about
translation invariant vector fields and differential operators
on the interior of $S^+N_p\Ct_a$.
Indeed, the restriction to $S^+N\Ct_a$ of the
lift of elements of $\Diffsc(X)$ (under $\betaSc$)
are such. We can see this since $\Vsc(X)$ is
given by sections of $\Tsc X$; the restriction of the lift of $P\in\Vsc(X)$
is then given by the identification of
$(\beta[X;\calC_a]^*_p\Tsc X)/\Tsc_p(\Ct_a;X)$ with the tangent space at
each point of the fiber $\rho_a^{-1}(p)$. Using the metric $g$ to identify
the quotient bundle with the orthocomplement of $\Tsc_p(\Ct_a;X)$,
$S^+N_p\Ct_a$ becomes an affine space with a translation-invariant
metric (i.e.\ `Euclidean')
with the metric induced by $g$; $\Delta_Y$ is the Laplacian of this metric.

Equations \eqref{eq:ham-4}-\eqref{eq:ham-5} show
that $\Hh_{a,0}(p,0)$ is uniformly bounded
below, so for any $\psi\in\Cinf_c(\Real)$ the set
\begin{equation}
\cup_a\cl(\{\xi\in\sct(\Ct_a;X):\ \psi(\Hh_{a}(\xi))\neq 0\})
\end{equation}
is compact.

The bound states of the subsystems of $H$ play an important role in
Euclidean many-body scattering. The appropriate replacement in the
general geometric setting is
given via the indicial operators of $H$. Thus, in this paper
the statement `no subsystem of $H$ has a bound state' means that
\begin{equation}\label{eq:hypo-0}
\Hh_{a,0}(\xi)\ \text{has no}\ L^2\ \text{eigenvalues for any}
\ a\neq 0\Mand \xi\in\sct(\Ct_a;X).
\end{equation}
Due to \eqref{eq:ham-4}-\eqref{eq:ham-5},
this means simply that
\begin{equation}\label{eq:h_a(p)-def}
h_a(p)=\Hh_{a,0}((p,0))\ \text{has no}\ L^2\ \text{eigenvalues for any}
\ a\neq 0\Mand p\in\Ct_a.
\end{equation}
In Euclidean scattering $h_a(p)$ is just the subsystem Hamiltonian
$h_a$ (which is then independent of $p$), so in that setting
\eqref{eq:hypo-0} indeed means that the
(proper) subsystems of $H$ have no bound states.

If no subsystem of $H$ has bound states it can be expected that
$\Delta-\lambda$ governs the propagation of singularities of distributions $u$
with $(H-\lambda)u\in\dCinf(X)$, except that the flow will break at
the places where $V$ is singular (i.e.\ where locally $V\nin\Cinf(X)$),
similarly to boundary and transmission problems for the wave equation
\cite[Chapter~XXIV]{Hor}, \cite{Melrose-Sjostrand:I, Lebeau:Propagation}.
Now, the
symbol of $\Delta-\lambda$ at $\bX$ (i.e.\ its
$\scl$-indicial operator) is $g-\lambda$. Hence, its characteristic variety is
\begin{equation}
\Sigma=\Sigma_{\Delta-\lambda}=\{\xi\in\sct_{\bX}X:\ g(\xi)-\lambda=0\}.
\end{equation}
The rescaled Hamilton vector field $\scHg=x^{-1}H_g$ of $g$ (or $g-\lambda$),
introduced in \cite{RBMSpec},
is
\begin{equation}
\scHg=2\tau(x\partial_x+\mu\cdot\partial_\mu+\nu\cdot\partial_\nu)-2h
\partial_\tau+H_h+xW',\quad W'\in\Vb(\sct X),
\end{equation}
so its restriction to $\bX$, also denoted by $\scHg$, is
\begin{equation}
\scHg=2\tau(\mu\cdot\partial_\mu+\nu\cdot\partial_\nu)-2h
\partial_\tau+H_h.
\end{equation}
Here $(y,z,\tau,\mu,\nu)$ denote coordinates about some $C=C_a$ as before,
though notice that $\mu\cdot\partial_\mu+\nu\cdot\partial_\nu$ is simply
the radial vector field in $T^*\bX$, so the above expression is indeed
invariant (as it must be).
The bicharacteristics of $\Delta-\lambda$ are just integral curves of
$\scHg$.

We divide the image $\dot\Sigma\subset\scdt X$
of $\Sigma$ under $\pi$ into a normal and a tangential
part,
\begin{equation}
\dot\Sigma=\Sigma_n(\lambda)\cup\Sigma_t(\lambda),
\end{equation}
as follows. Let $\pih$ be the restriction of $\pi$ to $\Sigma$. We let
\begin{equation}
\Sigma_n(\lambda)=\cup_a\{\xi\in\sct_{C'_a}(C_a;X)\cap\dot\Sigma:
\ \pih^{-1}(\xi)\ \text{consists of more than one point}\}
\end{equation}
and
\begin{equation}
\Sigma_t(\lambda)=\cup_a\{\xi\in\sct_{C'_a}(C_a;X)\cap\dot\Sigma:
\ \pih^{-1}(\xi)\ \text{consists of exactly one point}\}.
\end{equation}
In terms of our local coordinates around $C'_a$, in view of
\eqref{eq:g-at-C} and $|\mu_a|^2\geq 0$, this means that
\begin{equation}
\Sigma_n(\lambda)=\cup_a\{(z_a,\tau_a,\nu_a)
\in\sct_{C'_a}(C_a;X):\ \tau_a^2+\htil(z_a,\nu_a)<\lambda\}
\end{equation}
and
\begin{equation}
\Sigma_t(\lambda)=\cup_a\{(z_a,\tau_a,\nu_a)
\in\sct_{C'_a}(C_a;X):\ \tau_a^2+\htil(z_a,\nu_a)=\lambda\}.
\end{equation}
Notice that for $\xi=(z_a,\tau_a,\nu_a)\in\Sigma_t(\lambda)$
and the unique point $\xit=(0,z_a,\tau_a,\mu_a,\nu_a)\in\sct_{\bX}
X$ with $\pi(\xit)=\xi$ we have $\mu_a=0$. As the $\partial_{y_a}$
component of $\scHg$ is $2\mu_a\cdot\partial_{y_a}$ at $y_a=0$ (i.e.\ at
$C_a$), for such $\xi$ and $\xit$, $\scHg(\xit)$ is tangent
to $\sct_{C'_a} X$. On the other hand, if $\xi\in\Sigma_n(\lambda)$, $\xit\in
\pih^{-1}(\xi)$, then $\scHg(\xit)$ is normal to $\sct_{C'_a} X$, hence the
choice of our terminology. Notice also that on $\sct_{C'_0}X$, $\pi$ is
the identity map, so
\begin{equation}
\dot\Sigma\cap\sct_{C'_0}X\subset\Sigma_t(\lambda).
\end{equation}

We also define the radial sets $R_\pm(\lambda)$ as the sets
\begin{equation}
R_\pm(\lambda)=\pi(\{(y,z,\tau,\mu,\nu):\ \tau=\pm\sqrt{\lambda},
\ h(y,z,\mu,\nu)=0\}).
\end{equation}
Thus, $R_+(\lambda)\cup R_-(\lambda)$ is the image (under $\pi$)
of the set where $\scHg$ vanishes. Notice that
\begin{equation}
R_+(\lambda)\cup R_-(\lambda)\subset\Sigma_t(\lambda).
\end{equation}

Following Lebeau, we define generalized broken
bicharacteristics of $\Delta-\lambda$ as follows. First, we say that a function
$f\in\Cinf(\sct_{\bX} X)$ is $\pi$-invariant if for
$\xit,\xit'\in\sct_{\bX} X$, $\pi(\xit)=\pi(\xit')$ implies $f(\xit)=f(\xit')$.
A $\pi$-invariant function $f$ naturally defines a function $f_\pi$ on
$\scdt X$ by $f_\pi(\xi)=f(\xit)$ where $\xit\in\sct_{\bX}X$ is chosen so
that $\pi(\xit)=\xi$.

\begin{Def}\label{Def:gen-br-bichar}
Suppose that $(X,\calC)$ is locally linearizable.
A generalized broken bicharacteristic of $\Delta-\lambda$
is a continuous map
$\gamma:I\to\scdt X$,
where $I\subset\Real$ is an interval, satisfying
the following requirements:

\renewcommand{\theenumi}{\roman{enumi}}
\renewcommand{\labelenumi}{(\theenumi)}
\begin{enumerate}
\item
If $\xi_0=\gamma(t_0)\in\Sigma_t(\lambda)$
then for all $\pi$-invariant functions
$f\in\Cinf(\sct_{\bX}X)$,
\begin{equation}
\frac{d}{dt}(f_\pi\circ \gamma)(t_0)=\scHg f(\xit_0),\ \xit_0=\pih^{-1}(\xi_0).
\end{equation}

\item
If $\xi_0=\gamma(t_0)\in\Sigma_n(\lambda)\cap\sct_{C'_a} (C_a;X)$ then
there exists $\ep>0$ such that
\begin{equation}
t\in I,\ 0<|t-t_0|<\ep\Rightarrow\gamma(t)\nin \sct_{C'_a}(C_a;X).
\end{equation}
\end{enumerate}
\end{Def}

The success of this definition (so that it indeed describes what we
wish to describe) depends on a plentiful supply of $\pi$-invariant functions
on $\sct_{\bX} X$. Under our local linearizability hypothesis,
\eqref{eq:hypo-loc-lin}, there are always many such functions. Indeed,
by \eqref{eq:hypo-loc-lin},
if $p\in C'=C'_a$,
we can choose local coordinates $(y,z)$ on $\bX$ in terms
of which {\em all} the $C_b$ satisfying
$p\in C_b$ are linear, i.e.\ they are
are given by $A_b y=0$ where $A_b$ is a (constant) matrix, and
$C_a$ is given by $y=0$. Let $(\tau,\mu,\nu)$ denote the sc-dual
variables of $(x,y,z)$ as above.
Choosing
such coordinates, $y,z,\tau,\nu$ are $\pi$-invariant near $\sct_p X$.
In general, without the assumption \eqref{eq:hypo-loc-lin}, $\nu$
would not be $\pi$-invariant, and we
would not be able to modify it to make it such, so the
definition would be inadequate.

We can also arrange that the metric function is of the form
$h=\htil(z,\nu)+h_{nn}(z,\mu)$ at $C'_a$ by a further change of coordinates
$z_j'=z_j+\sum_{jk}Z_{jk}(z)y_k$, $y'=y$, which preserves the linear
structure of the $C_b$. In general we cannot arrange that
$h_{nn}(z,\mu)=|\mu|^2$
everywhere along $C'_a$ without destroying the product-linear structure of the
$C_b$. However, by a linear change in the $y$
coordinates we can make sure that {\em at a fixed $p\in C'_a$},
$h=\htil(z,\nu)+|\mu|^2$.
The continuity of a generalized broken bicharacteristic
$\gamma$ means that if $\gamma(t_0)\in\sct_{C'}(C;X)$,
then for $t$ near $t_0$,
$t\mapsto(y(\gamma(t)),z(\gamma(t)),\tau(\gamma(t)),\nu(\gamma(t)))$
is continuous, but $\mu(\gamma(t))$ may be discontinuous.
In terms of
Euclidean scattering this means that at $C_a$ the external momentum is
conserved, but not necessarily the internal one, while $\image(\gamma)\subset
\dot\Sigma$ corresponds to the conservation of kinetic energy. The
latter cannot be expected to hold if the subsystems of the Hamiltonian
have bound states; the relevant broken bicharacteristics in that case
exhibit more complex behavior. Another example of a $\pi$-invariant function
in this situation is $y\cdot\mu$; this will play a rather important role
in the propagation estimates. In fact, $\scHg (y\cdot \mu)(\xit_0)=2|\mu_0|^2$
if $\xit_0\in\sct_p X$ is of the form $(0,z_0,\tau_0,\mu_0,\nu_0)$,
so if $\pi(\xit_0)\in\Sigma_n(\lambda)$ and $\xit_0\in\Sigma$ then $y\cdot\mu$
is a parameter along generalized broken bicharacteristics near $\xit_0$ --
see also the following proposition.

A stronger characterization of generalized broken bicharacteristics
at $\Sigma_n(\lambda)$ follows as in Lebeau's
paper. Notice that if $\gamma:I\to\dot\Sigma$ is continuous then the
conclusion of the following proposition certainly implies (i) and (ii)
((ii) follows as $y_j=(y_a)_j$ are $\pi$-invariant), so the proposition
indeed provides an alternative to our definition.

\begin{prop}(Lebeau, \cite[Proposition~1]{Lebeau:Propagation})
\label{prop:Lebeau-2-sided}
If $\gamma$ is a generalized broken bicharacteristic as above, $t_0\in I$,
$\xi_0=\gamma(t_0)$, then there exist unique
$\xit_+, \xit_-\in\Sigma(\Delta-\lambda)$ satisfying
$\pi(\xit_\pm)=\xi_0$ and having
the property that if $f\in\Cinf(\sct_{\bX}X)$ is $\pi$-invariant
then $t\mapsto f_\pi(\gamma(t))$ is differentiable both from the left and
from the right at $t_0$ and
\begin{equation}\label{eq:Ham-23}
\left(\frac{d}{dt}\right)(f_\pi\circ\gamma)|_{t_0\pm}=\scHg f(\xit_{\pm}).
\end{equation}
\end{prop}

We refer to Lebeau's paper for the proof in the general setting,
but in the Appendix we give the proof under the assumption that the
elements of $\calC$ are totally geodesic. In fact, we prove slightly
more by giving a H\"older-type remainder estimate.
We present the proof in the Appendix,
but we emphasize that it is simply a minor modification of Lebeau's proof.
We remark that the most delicate part of the conclusion (under the
totally geodesic assumption) is the
differentiability of the `normal' coordinate functions
$y_j$ along $\gamma$, i.e.\ that of $y_j\circ\gamma$. Here we dropped
the projection $\pi$ from the notation (i.e.\ we did not write
$(y_j)_\pi\circ\gamma$) to simplify it; we will often do this in
the future for the other $\pi$-invariant coordinate functions $\tau$, $z_j$,
$\nu_j$. The proof proceeds by induction using the order on $\calC$. Thus,
we have to understand what happens near $t_0$ if $\gamma(t_0)=\xi_0\in
\Sigma_n(\lambda)\cap\sct_{C'_a}(C_a;X)$. The inductive hypothesis is that
we have already proved the proposition for $b$ with $C_a\subsetneq C_b$. Thus,
by Definition~\ref{Def:gen-br-bichar}, part (ii), it is true for $t_0$
replaced by $t\neq t_0$, assuming $|t-t_0|<\ep$. Hence, we need to
analyze the behavior of the coordinate functions using the Hamilton
equation, \eqref{eq:Ham-23} which is a little
more delicate than the positive commutator construction in
Proposition~\ref{prop:normal-prop-fine}, but the two proofs are very closely
related via the use of
same function $\phi$ to localize near (and along) the generalized broken
bicharacteristics.
A rather similar analogy arises in our tangential estimates in the totally
geodesic setting; see Propositions~\ref{prop:tgt-bichar-tot-geod} and
\ref{prop:tot-geod-tgt-prop} respectively.

We now describe some corollaries of this proposition. First, we remark that
the role of the globally defined $\pi$-invariant function $\tau$ is
somewhat analogous to the role played by the time variable in the
wave equation in Lebeau's paper. In particular,
$\tau$ gives a parameter along generalized broken bicharacteristics
with the exception of some trivial ones (namely the constant ones in
$R_+(\lambda)\cup R_-(\lambda)$). To see this, we show the following
corollary of the above proposition.

\begin{cor}\label{cor:tau-on-bichar}
Suppose that $\gamma:I\to\dot\Sigma$ is a generalized broken bicharacteristic.
Then $T=\tau_\pi\circ\gamma:I\to\Real$ is a $\Cinf$ function. In
addition, $T$ has one of the following forms. Either

\renewcommand{\theenumi}{\roman{enumi}}
\renewcommand{\labelenumi}{(\theenumi)}
\begin{enumerate}
\item
$T(t)=\sqrt{\lambda}$ for all $t\in I$, or

\item
$T(t)=-\sqrt{\lambda}$ for all $t\in I$, or

\item
$T'(t)<0$ for all $t$ and if $I=\Real$ then $T(t)\to\mp\sqrt{\lambda}$ as
$t\to\pm\infty$.
\end{enumerate}
\end{cor}

\begin{proof}
As $\lambda=\tau^2+h$ in $\Sigma_{\Delta-\lambda}$, we have for all
$\xit\in\Sigma_{\Delta-\lambda}$ that
\begin{equation}
\scHg\tau(\xit)=-2h(\xit)=2(\tau(\xit)^2-\lambda).
\end{equation}
Thus, with $T=\tau_\pi\circ\gamma$, the previous proposition implies
that for any $t\in I$, $T$ is differentiable
from both the left and the right at $t$, and both of these derivatives
are equal to $2(T(t)^2-\lambda)$. (We remark that this is proved directly
in the Appendix as a first step to the proof of the proposition.)
Thus, $T$ is $\calC^1$ and it satisfies
the ODE $dT/dt=2(T^2-\lambda)$. But, given say $T(t_0)=\tau_0$,
this ODE has a unique solution which is $\Cinf$. The last statement
follows by writing down the solution of the ODE explicitly, which, if
$T(t_0)\in(-\lambda,\lambda)$ for some $t_0\in I$, takes the form
$T(t)=-\sqrt{\lambda}\tanh(4\sqrt{\lambda}(t-c))$, $t\in I$,
for an appropriate constant $c$.
\end{proof}

Since for $\xi\in\dot\Sigma$ with $\tau(\xi)^2=\lambda$ we
automatically have $\xi\in R_+(\lambda)\cup R_-(\lambda)$, in (iii)
we see
that (if $I=\Real$)
as $t\to\pm\infty$, $\gamma(t)$ approaches $R_\mp(\lambda)$. In
addition, in the same case, as $T'$ never vanishes, $T\in(-\sqrt{\lambda},
\sqrt{\lambda})$ can be used to reparameterize $\gamma$ (reversing its
direction).

We proceed to examine generalized broken bicharacteristics in more detail,
starting with cases (i) and (ii).
Namely, we prove that generalized broken bicharacteristics through
$R_+(\lambda)\cup R_-(\lambda)$ are constant maps:

\begin{prop}\label{prop:rad-bichar}
If $\gamma:I\to\dot\Sigma$ is a generalized broken bicharacteristic,
$\gamma(t_0)=\xi_0\in R_+(\lambda)\cup R_-(\lambda)$, then $\gamma(t)=\xi_0$
for $t\in I$. Hence, $\hat\pi^{-1}\circ\gamma$
is a bicharacteristic of $\scHg$.
\end{prop}

\begin{proof}
The previous corollary and the above remarks show that
for all $t\in I$, $\gamma(t)\in R_+(\lambda)\cup
R_-(\lambda)$. Let $\xit(t)=\pih^{-1}(\gamma(t))$.
Thus, $\scHg$ vanishes at $\xit(t)\in
\hat\pi^{-1}(R_+(\lambda)\cup R_-(\lambda))$ for all $t$.
Since the base variables $y$ and $z$ are $\pi$-invariant, we conclude
that $d((y_j)_\pi\circ\gamma)/dt$ vanishes identically, hence $y$ is
constant, and similarly for $z$, proving that $\gamma(t)=\xi_0$ for all
$t$. The last statement of the proposition follows since $\scHg$ vanishes at
$\hat\pi^{-1}(R_+(\lambda)\cup R_-(\lambda))$.
\end{proof}

Now, we consider case (iii) of Corollary~\ref{cor:tau-on-bichar}.
Namely, we show that if we rescale and
reparametrize $\gamma$ and project off its $\tau$ component, we obtain
a generalized broken geodesic (of $h$) in $\bX$, broken at $\calC$.
This is a notion completely analogous to that of our
generalized broken bicharacteristics, and we proceed to define it.
Again, we need to introduce a `compressed' cotangent bundle. The metric
$h$ on $\bX$ naturally identifies the cotangent bundle $T^* C$
of $C\in\calC$ as a subset of $T^*\bX$. The compressed cotangent bundle
of $\bX$ is then
\begin{equation}
\dT \bX=\cup_a T^*_{C'_a}C_a.
\end{equation}
It is topologized by the projection $\pipa:T^*\bX\to\dT\bX$.
We also define the compressed cosphere bundle as the image of $S^*\bX$
under $\pi_\partial$; here $S^*\bX$ is the set of covcectors of unit length:
\begin{equation}
\dS\bX=\pi_\partial(S^*\bX).
\end{equation}
The restriction of $\pi_\partial$ to $\dS\bX$ is denoted by $\hat\pipa$.
This plays a role analogous to that of $\dot\Sigma$. We also define
its tangential and normal parts:
\begin{equation}
\dS_n\bX=\cup_a\{\zeta\in T^*_{C'_a}C_a\cap\dS\bX:
\ \pipah^{-1}(\zeta)\ \text{consists of more than one point}\}
\end{equation}
and
\begin{equation}
\dS_t\bX=\cup_a\{\zeta\in T^*_{C'_a}C_a\cap\dS\bX:
\ \pipah^{-1}(\zeta)\ \text{consists of exactly one point}\}.
\end{equation}
Generalized broken
geodesics are then defined as follows.

\begin{Def}\label{Def:gen-br-geod}
A generalized broken geodesic of $h$
is a continuous map
$\gammapa:I\to\dS \bX$,
where $I\subset\Real$ is an interval, satisfying
the following requirements:

\renewcommand{\theenumi}{\roman{enumi}}
\renewcommand{\labelenumi}{(\theenumi)}
\begin{enumerate}
\item
If $\zeta_0=\gammapa(t_0)\in\dS_t\bX$
then for all $\pipa$-invariant functions
$f\in\Cinf(T^*\bX)$,
\begin{equation}
\frac{d}{dt}(f_{\pipa}\circ \gammapa)(t_0)=H_{\half h} f(\zetat_0),
\ \zetat_0=\hat\pipa^{-1}(\zeta_0).
\end{equation}

\item
If $\zeta_0=\gammapa(t_0)\in\dS_n\bX\cap T^*_{C'_a}C_a$ then
there exists $\ep>0$ such that
\begin{equation}
t\in I,\ 0<|t-t_0|<\ep\Rightarrow\gammapa(t)\nin T^*_{C'_a}C_a.
\end{equation}
\end{enumerate}
\end{Def}

\begin{rem}
Sometimes, with an abuse of terminology, we also say that the projection
of a generalized broken geodesic to $\bX$ (via the projection
$\dS\bX\to\bX$ inherited from $T\bX$) is a generalized broken geodesic.
Indeed, this was the terminology used in the introduction.
\end{rem}

The metric $g$ gives rise to a product decomposition
\begin{equation}
\sct_{\bX} X=\Real_\tau\times T^*\bX.
\end{equation}
The compressed scattering cotangent bundle is thus also naturally a product:
\begin{equation}
\scdt X=\Real_\tau\times\dT\bX.
\end{equation}
We sometimes write the product variables as $\xi=(\tau,\xi'')$.
We write
\begin{equation}
p:\scdt X\to\dT\bX
\end{equation}
for the projection to the second factor.
Note that $\scdt X$ inherits a natural $\Real$-action from $\sct_{\bX} X$,
and if $\xi\in\dot\Sigma$, $\tau(\xi)^2\neq \lambda$, then
$\zeta=p((\lambda-\tau(\xi)^2)^{-1/2}\xi)\in\dS\bX$ since $h=\lambda-\tau^2$
on $\dot\Sigma$.

We also reparametrize generalized broken bicharacteristics $\gamma$
satisfying (iii) of Corollary~\ref{cor:tau-on-bichar} by letting $s=S(t)$
where $S$ satisfies $dS/dt=2(\lambda-\tau(\gamma(t))^2)^{1/2}$, with
$S(t_0)=s_0$
picked arbitrarily. We have the
following result.

\begin{prop}\label{prop:param-bichar}
Suppose that $\gamma:I\to\dot\Sigma$
is a generalized broken bicharacteristic which is
disjoint from $R_+(\lambda)\cup R_-(\lambda)$.
Then $\gamma\circ S^{-1}:J\to\dot\Sigma$, $S$ defined above, is given by
\begin{equation}\label{eq:br-bich-br-geod}
\tau=\sqrt{\lambda}\cos (s-s_1),
\ \xi''=\sqrt{\lambda}\sin (s-s_1)\gammapa(s)
\end{equation}
where $s_1$ is an appropriate constant and
$\gammapa:J\to\dS\bX$ is a generalized
broken geodesic, broken at $\calC$. If $I=\Real$, then
$J=(s_1,s_1+\pi)$, in particular
$J$ has length $\pi$,
and correspondingly the projection of $\gammapa$ to $\bX$ is a
curve of length $\pi$.
\end{prop}

\begin{proof}
Let
\begin{equation}
\gammapa(s)=p\left(\frac{\gamma(S^{-1}(s))}
{\sqrt{\lambda-\tau(\gamma(S^{-1}(s)))^2}}\right).
\end{equation}
Condition (ii) of Definition~\ref{Def:gen-br-bichar} implies (ii) of
Definition~\ref{Def:gen-br-geod} immediately.
Let $f\in\Cinf(T^*\bX)$ be a $\pipa$-invariant function.
Let
\begin{equation}
F(\xi)=f(p((\lambda-\tau(\xi)^2)^{-1/2}\xi));
\end{equation}
here we slightly abuse the notation and write $p:\sct_{\bX} X\to T^*\bX$.
Then $F$ is $\pi$-invariant, so (i) of
Definition~\ref{Def:gen-br-bichar} applies and gives $d(F_\pi\circ\gamma)/dt
(t_0)$. Since
$(F_\pi\circ\gamma)\circ S^{-1}=f_{\pipa}\circ\gammapa$, the chain rule
and a short calculation of $\scHg F$ gives (i) of
Definition~\ref{Def:gen-br-geod}.
The first equation in \eqref{eq:br-bich-br-geod}
follows since along $\gamma$, $ds/d\tau=(ds/dt)
(d\tau/dt)^{-1}=-(\lambda-\tau^2)^{-1/2}$. As $(\lambda-\tau^2)^{1/2}
=\sqrt{\lambda}\sin(s-s_1)$, the second equation follows as well.
Since $\tau\to\mp
\sqrt{\lambda}$ along $\gamma$ as $t\to\pm\infty$ and $\tau$ is decreasing,
we deduce the last statement.
\end{proof}

It is useful to introduce a relation on $\dS\times\dot\Sigma(\lambda)$
using the structure of the generalized broken bicharacteristics given
in this proposition.

\begin{Def}\label{Def:bichar-rel}
Suppose $\xi\in\dot\Sigma(\lambda)\setminus(R_-(\lambda)\cup
R_+(\lambda))$, $\zeta\in \dS\bX$.
We say that $\xi\sim_-\zeta$ if there is a generalized
broken bicharacteristic $\gamma:\Real\to\dot\Sigma(\lambda)$ with
$\gamma(t_0)=\xi$ such that $\gammapa:(a,a+\pi)\to S^*\bX$,
as in the above Proposition, satisfies $\lim_{s\to a+} \gammapa(s)
=\zeta$. We define $\xi\sim_+\zeta$ similarly by replacing $a+$ in the limit
by $(a+\pi)-$.
\end{Def}

We also need to analyze the uniform behavior of generalized broken
bicharacteristics. Here we quote Lebeau's results; they can also be proved
completely analogously to the proof of Proposition~\ref{prop:Lebeau-2-sided}
given here in the Appendix.

\begin{prop}(Lebeau, \cite[Proposition~5]{Lebeau:Propagation})
Suppose that $K$ is a compact subset of $\dot\Sigma$, $\gamma_n:[a,b]\to K$ is
a sequence of generalized broken bicharacteristics which converge uniformly
to $\gamma$. Then $\gamma$ is a generalized broken bicharacteristic.
\end{prop}

\begin{prop}(Lebeau, \cite[Proposition~6]{Lebeau:Propagation})
\label{prop:Lebeau-compactness}
Suppose that $K$ is a compact subset of $\dot\Sigma$, $[a,b]\subset\Real$
and
\begin{equation}
\calR=\{\text{generalized broken bicharacteristics}\ \gamma:[a,b]\to K\}.
\end{equation}
If $\calR$ is not empty then it is compact in the topology of uniform
convergence.
\end{prop}

\begin{cor}(Lebeau, \cite[Corollaire~7]{Lebeau:Propagation})
\label{cor:Lebeau-bichar-ext}
If $\gamma:(a,b)\to\Real$ is a generalized broken bicharacteristic then
$\gamma$ extends to $[a,b]$.
\end{cor}

\section{Generalized broken bicharacteristics for totally geodesic $\calC$}
\label{sec:tot-geod-bich}
We next examine the generalized broken bicharacteristics if all
elements of $\calC$
are totally geodesic with respect to $h$.
First we prove that generalized broken bicharacteristics
$\gamma:I\to\dot\Sigma$ with $\gamma(t_0)=\xi_0$,
$\xi_0\in\Sigma_t(\lambda)\cap\sct_{C'_a}(C_a;X)$
are actually bicharacteristics of $\scHg$
(and hence stay in $\sct_{C'_a}(C_a;X)$) for $t$ near $t_0$.

\begin{prop}\label{prop:tgt-bichar-tot-geod}
Suppose that all elements of $\calC$ are totally geodesic with respect
to $h$. Let
$\gamma:I\to\dot\Sigma$ be a generalized broken bicharacteristic,
\begin{equation}
\gamma(t_0)=\xi_0\in(\Sigma_t(\lambda)\cap\sct_{C'_a}(C_a;X))\setminus
(R_+(\lambda)\cup R_-(\lambda)).
\end{equation}
Then for $t\in J$, $J$ a neighborhood of
$t_0$, we have $\gamma(t)\in\Sigma_t(\lambda)\cap
\sct_{C'_a}(C_a;X)$, and $\gamma|_J$ is a bicharacteristic of $\scHg$.
\end{prop}

\begin{proof}
Our strategy consists of constructing a $\pi$-invariant function $\phi$
with $\scHg\phi\geq c>0$ in a neighborhood of $\hat\pi^{-1}(\xi_0)$.
Thus, by Proposition~\ref{prop:Lebeau-2-sided}, $d/dt(\phi_\pi(\gamma))|_{t\pm}
\geq c>0$ for $t\in J$, $J$ sufficiently small, so $\phi_\pi\circ\gamma$ is
increasing there. This will allow us to draw the desired conclusion for
the correct choice of $\phi$. We remark that this $\phi$ will reappear
in the proof of the propagation estimate in
Proposition~\ref{prop:tot-geod-tgt-prop}.
Moreover,
it is essentially the same as the corresponding function in the three-body
propagation estimate \cite[Proposition~15.4]{Vasy:Propagation-2}, though
we will use slightly different methods to estimate $\scHg\phi$.

In fact, first we find a $\pi$-invariant function $\omega$ such that
$\scHg$ will be appropriately small near $\hat\pi^{-1}(\xi_0)$.
So introduce coordinates centered at $C'_a$ as after
Definition~\ref{Def:gen-br-bichar}. Then the
metric function takes the form
\begin{equation}
h=\sum h^{ij}_{nn}(y,z)\mu_i\mu_j
+2\sum h^{ij}_{nt}(y,z)\mu_i\nu_j
+\sum h^{ij}_{tt}(y,z)\nu_i\nu_j
\end{equation}
with
\begin{equation}\label{eq:tot-geod-metric-2a}
h^{ij}_{nn}(0,0)=1,\ h^{ij}_{nt}(0,z)=0,
\end{equation}
and, due to the assumption that $C_a$ is totally geodesic,
\begin{equation}\label{eq:tot-geod-metric-2b}
\partial_y h^{ij}_{tt}(0,z)=0.
\end{equation}
We write
\begin{equation}
\htil(z,\nu)=\sum h^{ij}_{tt}(0,z)\nu_i\nu_j
\end{equation}
for the restriction of the tangential part of the metric function to $C_a$,
so
\begin{equation}
h|_{y=0}=\htil+\sum h^{ij}_{nn}(0,z)\mu_i\mu_j.
\end{equation}

Now, the Hamilton vector field of $h$ is given by
\begin{equation}\begin{split}
H_h=2&\sum_{i,j}h^{ij}_{nn}\mu_j\partial_{y_i}
+2\sum_{i,j}h^{ij}_{nt}\mu_i\partial_{z_j}
+2\sum_{ij}h^{ij}_{nt}\nu_j\partial_{y_i}
+2\sum_{i,j}h^{ij}_{tt}\nu_j\partial_{z_i}\\
&+\sum_{i,j,k}(\partial_{z_k}h^{ij}_{nn})\mu_i\mu_j\partial_{\nu_k}
+2\sum_{i,j,k}(\partial_{z_k}h^{ij}_{nt})\mu_i\nu_j\partial_{\nu_k}
+\sum_{i,j,k}(\partial_{z_k}h^{ij}_{tt})\nu_i\nu_j\partial_{\nu_k}+W'
\end{split}\end{equation}
with $W'=\sum\alpha_j\partial_{\mu_j}$. Hence, if
$\omega\in\Cinf(\Real^{m-1}_z\times\Rm_{\tau,\nu})$ then
\begin{equation}
H_h\omega|_{y=0}=H_{\htil}\omega
+\sum_{k}(\partial_{z_k}(h-\htil))\partial_{\nu_k}\omega.
\end{equation}
Now, $\mu$, hence $h-\htil$,
is small near $\hat\pi^{-1}(\xi_0)$, so to model
\begin{equation}
\scHg=2\tau(\mu\cdot\partial_\mu+\nu\cdot\partial_\nu)-2h
\partial_\tau+H_h,
\end{equation}
we introduce the
vector field
\begin{equation}
W=2\tau(\nu\cdot\partial_\nu)-2\htil\partial_\tau+H_{\htil}
\end{equation}
locally (near $\xi_0$) on $\sct (C_a;X)$.
Thus, we have
\begin{equation}
\scHg\omega|_{y=0}=W\omega-2(h-\htil)\partial_\tau\omega
+\sum_{k}(\partial_{z_k}(h-\htil))\partial_{\nu_k}\omega
\end{equation}
which is small if $W\omega$ is small.

We define $\omega$ as follows.
First, $W\tau=-2\htil$, and $\htil_{z_0}(\nu_0)\neq 0$ since
$\xi_0\nin R_+(\lambda)\cup R_-(\lambda)$, so near $\xi_0$,
$W\tau\neq 0$, i.e.\ $W$
is transversal to the hypersurface $\tau=\tau_0$. Thus, near $\xi_0$ in
$\sct(C_a;X)$ we can solve the Cauchy problem
\begin{equation}\label{eq:tot-geod-Cauchy}
W\omega=0,\qquad \omega|_{\tau=\tau_0}=(z-z_0)^2+(\nu-\nu_0)^2.
\end{equation}
Since $\omega$ and $d\omega$ vanish at $\xi_0$, the same holds on the
bicharacteristic of $W$ through $\xi_0$, but $\omega\geq 0$ and
the Hessian is still positive
in directions transversal to the bicharacteristics as these hold at $\xi_0$.
Moreover, by \cite[Lemma 7.7.2]{Hor},
\begin{equation}
|d\omega|\leq C\omega^{1/2}.
\end{equation}

Let
\begin{equation}
r_0=\tau^2+\htil_z(\nu)-\lambda,
\end{equation}
so $Wr_0=0$. At $\tau=\tau_0$ we have
$r_0=\htil_z(\nu)-\htil_{z_0}(\nu_0)$, so
\begin{equation}
|r_0|\leq C'|d\omega|\leq C''\omega^{1/2}
\end{equation}
when $\tau=\tau_0$, and then $W\omega=0=Wr_0$ implies that this
inequality holds everywhere.
Therefore,
\begin{equation}\label{eq:htil-h}
|\htil-h|\leq|\lambda-\tau^2-h|+|\lambda-\tau^2-\htil|\leq|\lambda-\tau^2
-h|+C\omega^{1/2}.
\end{equation}

Now,
\begin{equation}\begin{split}\label{eq:scHg-omega-est}
\scHg\omega=\scHg\omega-W\omega=&-2(h-\htil)\partial_\tau\omega\\
&+2\sum_{i,j}h^{ij}_{nt}(y,z)\mu_i\partial_{z_j}\omega
+2\sum_{i,j}(h^{ij}_{tt}(y,z)-h^{ij}_{tt}(0,z))\nu_j\partial_{z_i}\omega\\
&+\sum_{i,j,k}\partial_{z_k}h^{ij}_{nn}(y,z)\mu_i\mu_j\partial_{\nu_k}\omega
+2\sum_{i,j,k}\partial_{z_k}h^{ij}_{nt}(y,z)\mu_i\nu_j\partial_{\nu_k}\omega\\
&+\sum_{i,j,k}\partial_{z_k}(h^{ij}_{tt}(y,z)-h^{ij}_{tt}(0,z))
\nu_i\nu_j\partial_{\nu_k}\omega.
\end{split}\end{equation}
Thus, using  \eqref{eq:tot-geod-metric-2a}-\eqref{eq:tot-geod-metric-2b},
for some $C, C'>0$ we have
\begin{equation}\begin{split}\label{eq:scHg-omega-est-2}
|\scHg\omega-W\omega|&\leq C'(|\tau^2+h-\lambda|+\omega^{1/2}+|y|^2+|\mu|^2
+|\mu||y|) |d\omega|\\
&\leq C(|\tau^2+h-\lambda|+\omega^{1/2}+|y|^2+|\mu|^2) \omega^{1/2}.
\end{split}\end{equation}

Next, note that
\begin{equation}
\scHg |y|^2=4\sum_{i,j}h^{ij}_{nn}\mu_j y_i+4\sum_{i,j} h^{ij}_{nt}\nu_j y_i,
\end{equation}
so by \eqref{eq:tot-geod-metric-2a},
\begin{equation}
|\scHg |y|^2|\leq C|y|(|y|+|\mu|).
\end{equation}

For $\epsilon>0$ let
\begin{equation}\label{eq:def-phi-tot-geod}
\phi^{(\ep)}=\phi=\tau_0-\tau+\epsilon^{-1}|y|^2+\epsilon^{-2}\omega.
\end{equation}
Thus,
\begin{equation}\begin{split}\label{eq:scHg-phi-est-0'}
|\scHg\phi-2h|\leq C(\epsilon^{-1}|y|(|y|+|\mu|)+\epsilon^{-2}\omega^{1/2}
(|y|^2+|\mu|^2+|\tau^2+h-\lambda|+\omega^{1/2})).
\end{split}\end{equation}

We next estimate $\mu$. First, as $h^{ij}_{nn}(0,0)=\delta_{ij}$, $h_{nn}$
is positive definite in a small neighborhood of $(0,0)$ and
\begin{equation}
|\mu|^2\leq 2\sum_{i,j}h^{ij}_{nn}(y,z)\mu_i\mu_j
\end{equation}
there.
On the other hand,
\begin{equation}
\sum_{i,j}h^{ij}_{nn}(y,z)\mu_i\mu_j=h-\htil
-\sum_{i,j}h^{ij}_{nt}(y,z)\mu_i\nu_j
-\sum_{i,j}(h^{ij}_{tt}(y,z)-h^{ij}_{tt}(0,z))\nu_i\nu_j,
\end{equation}
so
\begin{equation}
|\sum_{i,j}h^{ij}_{nn}(y,z)\mu_i\mu_j|\leq|h-\htil|+C_1|y||\mu|+C_2|y|^2.
\end{equation}
Moving $C_1|y| |\mu|$ to the right hand side and completing the square gives
\begin{equation}
(|\mu|-C_3|y|)^2\leq |h-\htil|+C_4|y|^2,
\end{equation}
so
\begin{equation}
|\mu|\leq C(|h-\htil|^{1/2}+|y|),
\ \text{i.e.}\ |\mu|^2\leq C'(|h-\htil|+|y|^2).
\end{equation}

We can finally estimate $\scHg\phi$, using \eqref{eq:htil-h} as well:
\begin{equation}\begin{split}\label{eq:scHg-phi-est-0''}
|\scHg\phi-2h|\leq C(\epsilon^{-1}&|y|(|y|+\omega^{1/4}
+|\tau^2+h-\lambda|^{1/2})\\
&+\epsilon^{-2}\omega^{1/2}
(|y|^2+|\tau^2+h-\lambda|+\omega^{1/2})).
\end{split}\end{equation}
Note that $\phi_\pi(\xi_0)=0$, so near $\pi^{-1}(\xi_0)$, $\phi$ is small.
So now suppose that $0<\delta<1$ and
\begin{equation}\label{eq:near-xit_0-1}
\phi\leq 2\delta\Mand\tau-\tau_0\leq 2\delta.
\end{equation}
Then
\begin{equation}
\epsilon^{-1}|y|^2+\epsilon^{-2}\omega\leq 4\delta,
\end{equation}
so $|y|\leq (4\ep\delta)^{1/2}$, $\omega\leq 4\ep^2\delta$.
Hence, under the additional assumption
\begin{equation}\label{eq:near-xit_0-2}
|\tau^2+h-\lambda|<\epsilon\delta,
\end{equation}
i.e.\ that $\xit=(y,z,\tau,\mu,\nu)$ sufficiently close to $\Sigma_{\Delta
-\lambda}$, we have
\begin{equation}
|\scHg\phi-2h|\leq C(\ep^{-1}(\ep\delta)^{1/2}(\ep^2\delta)^{1/4}
+\ep^{-2}\ep^2\delta)\leq C'\delta^{3/4}.
\end{equation}
Since $h(\hat\pi^{-1}(\xi_0))>0$, we have $h(\xit)\geq 2c>0$ in a neighborhood
of $\hat\pi^{-1}(\xi_0)$. Now choose $\delta_0>0$ sufficiently small,
so that $C'\delta_0^{3/4}<c$. Note that this requirement is independent of
$\ep$. We thus conclude that for $\delta\in (0,\delta_0)$, $\xit$ satisfying
\eqref{eq:near-xit_0-1} and \eqref{eq:near-xit_0-2}, we have
\begin{equation}\label{eq:scHg-phi-est-tot-geod}
\scHg\phi(\xi)\geq c>0.
\end{equation}

Now, using the result of Proposition~\ref{prop:Lebeau-2-sided},
let $\xit_{\pm}(t)\in\Sigma(\Delta-\lambda)$ be the unique points
such that $\pi(\xit_{\pm}(t))=\gamma(t)$ and for all $\pi$-invariant $f$
\begin{equation}\label{eq:2-ders}
\left(\frac{d}{dt}\right)(f_\pi\circ\gamma)|_{t\pm}=\scHg f(\xit_{\pm}(t)).
\end{equation}
Choosing a sufficiently small open interval $J$ around $t_0$,
$\tau(\gamma(t))$, hence $\tau(\xit_\pm(t))$, automatically satisfies
\eqref{eq:near-xit_0-1} for $t\in J$, while \eqref{eq:near-xit_0-2}
holds automatically as $\xit_\pm(t)\in\Sigma(\Delta-\lambda)$. Thus,
applying \eqref{eq:2-ders} with $\phi$ in place of $f$, we see that, with
\begin{equation}\label{eq:tot-geod-36}
g(t)=\phi_\pi\circ\gamma(t),
\end{equation}
we have
\begin{equation}\label{eq:2-ders-phi}
t\in J\Mand g(t)\leq 2\delta\Rightarrow
\left(\frac{dg}{dt}\right)|_{t\pm}\geq c>0.
\end{equation}
As $g$ is continuous and $g(t_0)=0$,
this shows that $g$ is increasing on $J\cap(-\infty,t_0]$.
To see this, first note that $g(t)< 2\delta$ on $J\cap
(-\infty,t_0]$, for otherwise $g^{-1}(\{2\delta\})\cap(-\infty,t_0]\cap J$
is not empty,
$g^{-1}(\{2\delta\})\cap(-\infty,t_0]$ is closed, so
taking $t_1=\sup (g^{-1}(\{2\delta\})\cap(-\infty,t_0])<t_0$ and $t_1\in J$. 
Thus,
for $t\in[t_1,t_0]$, $g$ is differentiable from either side at $t$ and the
derivatives are both positive, so $g$ is increasing on $[t_1,t_0]$,
hence $g(t_1)\leq g(t_0)=0$ contradicting $g(t_1)=2\delta$. Thus,
$g<2\delta$ on $J\cap(-\infty,0]$, so $g$ is increasing here, so
$g(t)\leq 0$ for $t\in J\cap(-\infty,t_0)$.
Taking into account the definition of $\phi$ we immediately deduce that
\begin{equation}
|y(\gamma(t))|\leq C\ep^{1/2},\ t\in J\cap(-\infty,t_0).
\end{equation}
Since $\ep\in(0,1)$ is arbitrary, we conclude that $y(\gamma(t))=0$ for
$t\in J\cap(-\infty,t_0]$, so $\gamma(t)\in\sct(C_a;X)$ for such $t$.
Similarly,
$\omega(\gamma(t))=0$ for such $t$, so by the construction of $\omega$,
$\gamma(t)$ is the integral curve of $W$ through $\xi_0$ (for $t\in J$,
$t\leq t_0$). Of course, a similar argument (with a change of sign in
$\tau_0-\tau$ in \eqref{eq:def-phi-tot-geod}) works for $J\cap[0,\infty)$,
so we conclude that $\gamma|_J\subset\sct_{C'_a}(C_a;X)$ and $\gamma|_J$
is an integral curve of $W$. As $W$ preserves $\tau^2+\htil$ (being
essentially its rescaled Hamilton vector field), $\tau^2(\gamma(t))+
\htil(\gamma(t))=\lambda$, $t\in J$, so $\gamma|_J\subset\Sigma_t(\lambda)$,
and hence at $\hat\pi^{-1}(\gamma|_J)$, $\scHg$ and $W$ agree and $\gamma|_J$
is a bicharacteristic of $\scHg$ as claimed.
\end{proof}

Next, we prove that if $\xi_0\in\Sigma_n(\lambda)\cap\sct_{C'_a}(C_a;X)$,
$\gamma(t_0)=\xi_0$, $\gamma$ is a generalized broken bicharacteristic,
then for a sufficiently small $\delta>0$, $\gamma|_{[0,\delta]}$
is a generalized
broken bicharacteristic of $\Delta-\lambda$, broken at $\calC'\subset\calC$,
where $\calC'$ is cleanly intersecting and $C_a\nin\calC'$. This will
{\em not} use that $\calC$ is totally geodesic.

\begin{prop}\label{prop:normal-bichar-tot-geod}
Suppose that $\xi_0\in\Sigma_n(\lambda)\cap\sct_{C'_a}(C_a;X)$, $\gamma$
is a generalized broken geodesic with $\gamma(t_0)=\xi_0$ and
$\xit_+$ is as in Proposition~\ref{prop:Lebeau-2-sided}.
Suppose that $\xit_+\in\sct(C_b;X))$ and $b$ is minimal with this property
(i.e.\ $C_c\subset C_b$ and $\xit_+\in\sct(C_c;X)$ imply $c=b$).
Let
\begin{equation}
\calC'=\calC\setminus\{C_c:\ C_c\cap C_b\subset C_a\}.
\end{equation}
Then for sufficiently small $\delta>0$, $\gamma|_{[0,\delta]}$ is
a generalized broken bicharacteristic of $\Delta-\lambda$, broken at
$\calC'$, and $\gamma((0,\delta])$ is disjoint from $\sct(C_c;X)$ if
$C_c\nin\calC'$.
\end{prop}

\begin{proof}
Let $b$ be as above and introduce local coordinates centered at $C'_a$.
We may assume that $C_b$ is given by $y'=0$ for a suitable splitting
$y=(y',y'')$. Thus, $\xit_+$ is of the form
$\xit_+=(0,0,\tau_0,0,\mu''_0,\nu_0)$, and as $\xit_+\in\Sigma_n(\lambda)$,
$\mu''_0\neq 0$. By
Proposition~\ref{prop:Lebeau-2-sided}, taking into account that $y$ is
$\pi$-invariant,
\begin{equation}
d(y'_j\circ\gamma)/dt|_{t_0+}=0,\ d(y''_j\circ\gamma)/dt|_{t_0+}=(\mu''_0)_j.
\end{equation}
Since $\mu''_0\neq 0$, there exist $c>0$, $\delta_0>0$,
such that $|y''(\gamma(t))|\geq c(t-t_0)$ for $t\in(t_0,t_0+\delta_0)$, while
for any $\ep>0$ there exists $\delta_1>0$ such that
$|y'(\gamma(t))|\leq\ep (t-t_0)$ for $t\in(t_0,t_0+\delta_1)$.
In particular, for any $\ep>0$ there exists $\delta>0$ such that
for $t\in(t_0,t_0+\delta)$ we have $|y'(\gamma(t))|/|y''(\gamma(t))|\leq
\ep$. By choosing $\ep>0$ sufficiently small we can thus make sure that
$\gamma(t)\nin C_c$ for $t\in(t_0,t_0+\delta]$ if $C_c\nin\calC'$.
Hence, $\gamma|_{[t_0,t_0+\delta]}$ can be regarded as a curve in $\cup_{C_c\in
\calC'}\sct_{C'_c}(C_c;X)$, $C'_c$ taken with respect to $\calC'$, if we let
$\gamma(t_0)=\pi_{0b}(\xit_0)\in\sct(C_b;X)$. Of course, $\gamma|_{(t_0,
t_0+\delta]}$ {\em is} a generalized broken bicharacteristic, broken at
$\calC'$ (since it has no points above $\calC\setminus\calC'$). Thus,
by Corollary~\ref{cor:Lebeau-bichar-ext}, $\gamma|_{(t_0,t_0+\delta]}$
extends to a generalized broken bicharacteristic, broken at $\calC'$,
defined on $[t_0,t_0+\delta]$; by continuity of $\gamma$ this must
coincide with $\gamma$, so $\gamma$ is a generalized broken bicharacteristic,
broken at $\calC'$, as claimed.
\end{proof}

We can combine the previous results to deduce the structure of the
generalized broken bicharacteristics if $\calC$ is totally geodesic.

\begin{prop}
Suppose that $\calC$ is totally geodesic with respect to $h$
and $\gamma$ is a generalized
broken bicharacteristic, broken at $\calC$ with
$\xi_0=\gamma(t_0)\in\sct_{C'_a}(C_a;X)$.
Then there exists $\delta>0$ such that both
$\gamma|_{[t_0,t_0+\delta)}$ and $\gamma|_{(t_0-\delta,t_0]}$ are
bicharacteristics of $\scHg$.
\end{prop}

\begin{proof}
If $\xi_0\in R_+(\lambda)\cup R_-(\lambda)$ then $\gamma(t)=\xi_0$ for
$t$ near $t_0$ by Proposition~\ref{prop:rad-bichar}, hence near $t_0$,
$\gamma$ is a ($\pi$-projected)
bicharacteristic of $\scHg$ (as $\scHg$ vanishes at $R_+
(\lambda)\cup R_-(\lambda)$).
If $\xi_0\in\Sigma_t
(\lambda)\setminus(R_+(\lambda)\cup R_-(\lambda))$
then Proposition~\ref{prop:tgt-bichar-tot-geod} applies and
proves the result. If $\xi_0\in\Sigma_n(\lambda)$, then
with $\calC'$ as in Proposition~\ref{prop:normal-bichar-tot-geod},
$\gamma|_{[0,\delta)}$ is a generalized broken bicharacteristic,
broken at $\calC'$, with $\gamma(t_0)\in\sct_{C'_b}(C_b;X)\cap\Sigma_t
(\lambda)$ (prime taken
with respect to $\calC'$). Thus, Proposition~\ref{prop:tgt-bichar-tot-geod}
applies again and proves the result.
\end{proof}

A compactness argument gives at once

\begin{cor}
If $\gamma:[a,b]\to\dot\Sigma$ is a generalized broken bicharacteristic,
broken at $\calC$, and $\calC$ is totally geodesic, then there exist
$t_0=a<t_1<t_2<\ldots<t_m=b$ such that $\gamma|_[t_j,t_{j+1}]$ is
a bicharacteristic of $\Delta-\lambda$ (i.e.\ it is {\em not} broken).
\end{cor}

\section{Positive operators}\label{sec:positive}
In the following two sections we discuss technical points of the
microlocal positive commutators constructions. In this section we show
roughly speaking
that the positivity of the indicial operators of $A\in\PsiSc^{-\infty,0}
(X,\calC)$ implies the positivity of $A$ modulo compact operators.
We prove this by constructing an approximate square root of $A$.
In the next section we examine commutators $[A,H]$ in more detail.

Throughout this section we assume that $H$
is a many-body Hamiltonian.
We start with the basic square root construction.

\begin{lemma}
Suppose that $H$
is a many-body Hamiltonian and $\lambda\in\Real$.
Suppose also that $A\in\PsiSc^{-\infty,0}(X,\calC)$ is self-adjoint,
and for some $c>0$ and $\psi\in\Cinf_c(\Real)$ which
is identically $1$ near $\lambda$,
\begin{equation}
\psi(H)A\psi(H)\geq c\psi(H)^2.
\end{equation}
Then for any $c'\in(0,c)$ and $\phi\in\Cinf_c(\Real)$ such that
\begin{equation}
\supp\phi\cap\supp(1-\psi)=\emptyset,
\end{equation}
there exists $B\in\PsiSc^{-\infty,0}(X,\calC)$ such that
\begin{equation}\label{eq:fc-49}
\phi(H)(A-c')\phi(H)=\phi(H)B^*B\phi(H).
\end{equation}
\end{lemma}

\begin{proof}
Let
\begin{equation}\label{eq:fc-50}
P=\psi(H)A\psi(H)+c(\Id-\psi(H)^2)\in\PsiSc^{0,0}(X,\calC).
\end{equation}
Thus, $P\geq c$, so
$P-c'\geq c-c'>0$. Since the spectrum of $P-c'$ is a subset of
$[c-c',\infty)$ and $c-c'>0$,
we have $(P-c')^{1/2}=f(P-c')$ where $f\in\Cinf_c(\Real)$
and $f(t)=\sqrt{t}$ if $t$ is in the spectrum of $P-c'$.
By Proposition~\ref{prop:fc-3},
\begin{equation}\label{eq:fc-52}
Q=(P-c')^{1/2}=f(P-c')\in\PsiSc^{0,0}(X,\calC).
\end{equation}
Let $\psi_1$ be identically $1$ near $\supp\phi$ and vanish near
$\supp(1-\psi)$. Then
\begin{equation}\label{eq:fc-53}
\psi_1(H)Q^2\psi_1(H)=\psi_1(H)(P-c')\psi_1(H)=\psi_1(H)(A-c')\psi_1(H).
\end{equation}
Now let $\phi\in\Cinf_c(\Real)$ be identically $1$ near $\lambda$ and
vanish near $\supp(1-\psi_1)$. Let
\begin{equation}\label{eq:fc-59}
B=Q\psi_1(H)\in\PsiSc^{-\infty,0}(X,\calC).
\end{equation}
Multiplying \eqref{eq:fc-53} from both sides by $\phi(H)$ then proves
\eqref{eq:fc-49}.
\end{proof}

We now show that under certain additional assumptions,
the positivity of the indicial operators implies positivity of the operator
modulo lower order (hence compact) terms in the calculus. We start by
assuming strict positivity of the indicial operators when localized
in the spectrum of $H$.

\begin{prop}\label{prop:str-pos-ops}
Suppose that $H$
is a many-body Hamiltonian and $\lambda\in\Real$.
Suppose also that $A,C\in\PsiSc^{-\infty,0}(X,\calC)$ are self-adjoint
and $\Ch_{a,0}(\zeta)=c_a(\zeta)\psi_0(\Hh_a(\zeta))^2$ for every
$a$ and $\zeta\in\sct(\Ct_a;X)$ where $c_a(\zeta)$ is a function
with $c_a(\zeta)>0$, $\psi_0\equiv 1$ near $\lambda\in
\Real$, $\psi_0\in\Cinf_c(\Real)$.
Assume in addition that
there exists $\psi\in\Cinf_c(\Real)$ which
is identically $1$ near $\lambda$, $\supp\psi\cap\supp(1-\psi_0)=
\emptyset$, such that
\begin{equation}
\psi(\Hh_a(\zeta))\Ah_a(\zeta)\psi(\Hh_a(\zeta))\geq \psi(\Hh_a(\zeta)
)c_a(\zeta)\psi(\Hh_a(\zeta))
\end{equation}
for every $a$ and $\zeta\in\sct(\Ct_a;X)$.
Then for any $\ep\in(0,1)$ and
$\phi\in\Cinf_c(\Real)$ with
\begin{equation}
\supp\phi\cap\supp(1-\psi)=\emptyset,
\end{equation}
there exists
$R\in\PsiSc^{-\infty,1}(X,\calC)$
such that
\begin{equation}
\phi(H)A\phi(H)\geq (1-\ep) \phi(H)C\phi(H)+R.
\end{equation}
\end{prop}

\begin{proof}
We apply a parameter dependent version of the previous lemma to
the indicial operators to conclude that for each $\zeta$ there
exists $\Bh_a(\zeta)$ with
\begin{equation}\label{eq:fc-63}
\phi(\Hh_a(\zeta))(\Ah_a(\zeta)-(1-\ep)\Ch_a(\zeta))\phi(\Hh_a(\zeta))
=\phi(\Hh_a(\zeta))\Bh_a(\zeta)^*\Bh_a(\zeta)\phi(\Hh_a(\zeta)).
\end{equation}
It follows from the Cauchy integral formula construction of the
square root in the calculus and the explicit formulae
\eqref{eq:fc-50}, \eqref{eq:fc-52} and \eqref{eq:fc-59}
that the indicial operators $\Bh_a(\zeta)$
match up so that there exists $B\in\PsiSc^{-\infty,0}(X,\calC)$ with
indicial operators $\Bh_a(\zeta)$. Here note that the set where
$\psi(\Hh_a(\zeta))$ does not vanish has compact closure, hence $c$ is bounded
below on it by a positive constant. Thus, we can take the same smooth
function $f$ in the expression \eqref{eq:fc-52} for the square root for every
$a$ and $\zeta$.
By \eqref{eq:fc-63},
\begin{equation}
\phi(H)(A-(1-\ep)C)\phi(H)=\phi(H)B^*B\phi(H)+R
\end{equation}
with $R\in\PsiSc^{-\infty,1}(X,\calC)$. Since $\phi(H)B^*B\phi(H)\geq 0$,
rearranging this proves the proposition.
\end{proof}

Similar conclusions hold if we assume a two-sided estimate on the
indicial operators of $A$. In essence, this forces the indicial operators,
hence their square roots, to vanish to infinite order when $c$ vanishes.

\begin{prop}\label{prop:pos-ops}
Suppose that $H$
is a many-body Hamiltonian and $\lambda\in\Real$.
Suppose also that $A,C\in\PsiSc^{-\infty,0}(X,\calC)$ are self-adjoint
and $\Ch_{a,0}(\zeta)=c_a(\zeta)\psi_0(\Hh_a(\zeta))^2$ for every
$a$ and $\zeta\in\sct(\Ct_a;X)$ where $c_a(\zeta)$ is a function
with $c_a(\zeta)\geq 0$ which vanishes with all derivatives at each $\zeta$
with $c_a(\zeta)=0$,
$\psi_0\equiv 1$ near $\lambda\in
\Real$, $\psi_0\in\Cinf_c(\Real)$, $\Ah_a(\zeta)=0$ if $c_a(\zeta)=0$, and
for any differential
operator $Q\in\Diff(\sct(\Ct_a;X))$, all seminorms of
$Q(c_a(\zeta)^{-1}\Ah_a(\zeta))$ in
$\PsiSc^{-\infty,0}(\rho_a^{-1}(p),T_p\calC^a)$,
$\zeta\in\sct_p(\Ct_a;X)$,
are uniformly bounded on the set of $\zeta$'s with $c_a(\zeta)> 0$.
(This is almost, but not quite,
a statement about the seminorms of $c_a(\zeta)^{-1}\Ah_a(\zeta)$
in $\PsiScra^{-\infty,0}(\rho_a^*\sct(\Ct_a;X),\calCt_a)$,
because we restrict our attention
to the region where $c_a(\zeta)> 0$, and do so uniformly.)

Assume in addition that
there exists $\psi\in\Cinf_c(\Real)$ which
is identically $1$ near $\lambda$, $\supp\psi\cap\supp(1-\psi_0)=
\emptyset$, such that
\begin{equation}
\psi(\Hh_a(\zeta))\Ah_a(\zeta)\psi(\Hh_a(\zeta))\geq \psi(\Hh_a(\zeta)
)c_a(\zeta)\psi(\Hh_a(\zeta))
\end{equation}
for every $a$ and $\zeta\in\sct(\Ct_a;X)$.
Then the conclusion of the previous proposition holds, i.e.\ for any
$\ep\in(0,1)$ and
$\phi\in\Cinf_c(\Real)$ with
\begin{equation}
\supp\phi\cap\supp(1-\psi)=\emptyset,
\end{equation}
there exists
$R\in\PsiSc^{-\infty,1}(X,\calC)$, with seminorms bounded by those of $A$ and
$C$ in $\PsiSc^{-\infty,0}(X,\calC)$, and with
$\WFScp(R)\subset\WFScp(A)\cup\WFScp(C)$
such that
\begin{equation}
\phi(H)A\phi(H)\geq (1-\ep) \phi(H)C\phi(H)+R.
\end{equation}
\end{prop}

\begin{proof}
We define $\Bh_a(\zeta)=0$ if $c_a(\zeta)=0$, otherwise we define
$\Bh_a(\zeta)$ as in the previous proposition.
The only additional ingredient is the analysis of $\Bh_a(\zeta)$ near
$\zeta$ with $c_a(\zeta)=0$. To do this analysis, we follow the construction
of $\Bh_a(\zeta)$ in detail. So let
\begin{equation}
\Ph_a(\zeta)=\psi(\Hh_a(\zeta))\Ah_a(\zeta)\psi(\Hh_a(\zeta))+c_a(\zeta)
(\Id-\psi(\Hh_a(\zeta))^2),
\end{equation}
and let
\begin{equation}
c'_a(\zeta)=(1-\ep)c_a(\zeta).
\end{equation}
Thus, $\Ph_a(\zeta)-c'_a(\zeta)\geq\ep c_a(\zeta)$. Let
\begin{equation}
\Qh_a(\zeta)=(\Ph_a(\zeta)-c'_a(\zeta))^{1/2}=c_a(\zeta)^{1/2}(c_a(\zeta)^{-1}
\Ph_a(\zeta)-(1-\ep))^{1/2}.
\end{equation}
By our assumption, there exists $M>0$
such that the norm of $\Ph_a(\zeta)$ in $\bop(L^2,L^2)$ is bounded by
$Mc_a(\zeta)$. Now choose $f\in\Cinf_c(\Real)$ such that $f(t)=\sqrt{t}$ on
$[1-\ep,M]$. Then $M\geq c_a(\zeta)^{-1}\Ph_a(\zeta)-1+\ep\geq\ep$,
so
\begin{equation}
\Qh_a(\zeta)=c_a(\zeta)^{1/2}f(c_a(\zeta)^{-1}\Ph_a(\zeta)-(1-\ep)).
\end{equation}
By our assumptions, the seminorms of $c_a(\zeta)^{-1}\Ph_a(\zeta)$ in
$\PsiSc^{0,0}(
\rho_a^{-1}(p),T_p\calC^a)$, $\zeta\in\sct_p(\Ct_a;X)$, remain uniformly
bounded as $c_a(\zeta)\to 0$, so the Cauchy integral representation of $f$,
via an almost analytic extension, shows that
$f(c_a(\zeta)^{-1}\Ph_a(\zeta)-(1-\ep))$ remains uniformly bounded.
Thus, $\Qh_a(\zeta)$ is continuous
as a function on $\sct(\Ct_a;X)$ with values in $\PsiSc^{0,0}(
\rho_a^{-1}(p),T_p\calC^a)$. A similar argument also holds for the
derivatives of $\Qh_a(\zeta)$.
Let $\psi_1$ be identically $1$ near $\supp\phi$ and vanish near
$\supp(1-\psi)$, and let
\begin{equation}
\Bh_a(\zeta)=\Qh_a(\zeta)\psi_1(H).
\end{equation}
Again, the $\Bh_a(\zeta)$ match up so there exists $B\in\PsiSc^{-\infty,0}(X,
\calC)$ with these indicial operators. We can also make sure that the
lower order terms also vanish where $c$ does, i.e.\ that
$\WFScp(B)\subset\supp c$. Then the indicial operators of
$\phi(H)(A-(1-\ep)C)\phi(H)$ and $\phi(H)B^*B\phi(H)$ are the same, so
\begin{equation}
\phi(H)(A-(1-\ep)C)\phi(H)=\phi(H)B^*B\phi(H)+R
\end{equation}
with $R\in\PsiSc^{-\infty,1}(X,\calC)$, proving the proposition.
\end{proof}

\section{Commutators}\label{sec:commutators}
In this section we discuss the basic technical tool underlying the
propagation estimates of the following sections. Thus, we show how
an estimate of the commutator $[A,H]$ at $\Ct_0$,
which is essentially obtained by
a symbolic calculation in the scattering calculus, can give a
positive commutator
estimate under the additional assumption that $\Hh_{a,0}(\zeta)$
has no $L^2$ eigenfunctions for any $a\neq 0$ and $\zeta\in\sct(\Ct_a;X)$.
In the Euclidean setting this means simply that the subsystems have no
bound states.

To do so,
we extend the notion of a function being $\pi$-invariant to functions
on $\sct X$ in a trivial way: $q\in\Cinf(\sct X)$ is $\pi$-invariant
if $q|_{\sct_{\bX}X}$ is $\pi$-invariant. Since the analysis
of classical dynamics,
i.e.\ of generalized broken bicharacteristics of $\Delta-\lambda$,
broken at $\calC$, is based on the properties of $\pi$-invariant
functions, we will be interested in quantizing $\pi$-invariant
symbols. More specifically, we are essentially
interested in operators of the form $A=Q\psi_0(H)$,
$\psi_0\in\Cinf_c(\Real)$, where $Q$ is obtained by quantizing a
$\pi$-invariant function $q\in\Cinf(\sct X)$. Since such $Q$ would not
be in our calculus, we construct $A$ directly.

All considerations
in what follows will be local, i.e.\ we will assume that the projection
of the support of $q$ to $X$ lies near a fixed $p\in\bX$, so we can always
work in local coordinates and identify $X$ with $\Snp$.
The problem with such $q\in\Cinf(\Snp\times\Rn)$ is that they are rarely
in $\Cinf(\Snp\times\Snp)$, i.e.\ they are not symbols in $\xi$, so
$Q$ will not be in $\Psisc^{0,0}(\Snp)$ or indeed in
$\PsiSc^{0,0}(\Snp,\calC)$. This, however, is not a major difficulty.
Fix $\psi_0\in\Cinf_c(\Real;[0,1])$ which is identically $1$ in a neighborhood
of a fixed $\lambda$.
Thus,
$\psi_0(H)\in\PsiSc^{-\infty,0}(X,\calC)$, so it is smoothing.
At the symbol level,
$\psi_0(H)$ is locally
the right quantization of some
\begin{equation}\label{eq:comm-3}
p\in\Cinf([\Snp;\calC]\times\Snp)
\end{equation}
which vanishes to infinite order at $[\Snp;\calC]\times\partial\Snp$,
which will enable us to write down $A$ directly.

We are thus interested in the following class of symbols $q$. We assume
that $q\in\Cinf(\Rn_w\times\Rn_\xi)$ and that for
every multiindex $\alpha$, $\beta\in\Nat^n$ there exist constants
$C_{\alpha,\beta}$ and $m_{\alpha,\beta}$ such that
\begin{equation}\label{eq:comm-4-a}
|(D_w^\alpha
D_\xi^\beta q)(w,\xi)|\leq C_{\alpha,\beta}\langle w\rangle^{-|\alpha|}
\langle\xi\rangle^{m_{\alpha,\beta}}.
\end{equation}
This implies, in particular, that
\begin{equation}
q\in\bcon^0(\Snp\times\Rn),
\end{equation}
i.e.\ that $q$ is a 0th order symbol in $w$, though it may blow up
polynomially in $\xi$. Indeed, in the compactified notation,
\eqref{eq:comm-4-a} becomes that for every $P\in\Diffb(\Snp)$,
acting in the base ($w$) variables, and for every
$\beta\in\Nat^n$ there exist $C_{P,\beta}$ and $m_{P,\beta}$ such that
\begin{equation}\label{eq:comm-4-b}
|(P D_\xi^\beta) q| \leq C_{P,\beta}
\langle\xi\rangle^{m_{P,\beta}}.
\end{equation}

It is convenient to require that $q$ be polyhomogeneous on $\Snp\times\Rn$:
\begin{equation}\label{eq:comm-5}
q\in\Cinf(\Snp\times\Rn);
\end{equation}
this stronger statement automatically holds for the $\pi$-invariant
symbols we are interested in.

We next introduce the product symbol
\begin{equation}\label{eq:comm-7}
a(w,w',\xi)=q(w,\xi)p(w',\xi),
\end{equation}
where $\psi_0(H)$ is given locally by the right quantization of $p$.
The main point is

\begin{lemma}\label{lemma:comm-1}
The symbol $a$ defined by \eqref{eq:comm-7} is in $\Cinf(\Snp
\times[\Snp;\calC]
\times\Snp)$ and it vanishes with all derivatives at
$[\Snp;\calC]\times\partial\Snp$. Hence, it defines an operator
$A\in\PsiSc^{-\infty,0}(X,\calC)$ by the oscillatory integral
\eqref{eq:Psop-37}.
\end{lemma}

\begin{proof}
First, $a\in\Cinf(\Snp\times[\Snp;\calC]\times\Rn)$ follows from
\eqref{eq:comm-3} and \eqref{eq:comm-5}. Moreover, the infinite
order vanishing of $p$ at $[\Snp;\calC]\times\Sn$
implies that for every $P'\in\Diffb([\Snp;\calC])$, $\beta\in\Nat^n$ and
$N\in\Nat$,
\begin{equation}
|P'D^\beta_{\xi}p|\leq C_{P',\beta,N}\langle\xi\rangle^{-N}.
\end{equation}
Thus, Leibniz' rule shows that for $P\in\Diffb(\Snp)$ acting in $w$,
$P'\in\Diffb([\Snp;\calC])$ acting in $w'$, $\beta\in\Nat^n$ and $N$
\begin{equation}
|PP'D^\beta_{\xi}a|\leq C_{PP',\beta,N}\langle\xi\rangle^{-N}.
\end{equation}
But this means precisely that $a\in\Cinf(\Snp
\times[\Snp;\calC]\times\Snp)$ and it vanishes to infinite order at
the boundary in the last factor.
\end{proof}

The indicial operators of $A$ are just given by the quantization of
the appropriate restriction of $a$ similarly to \eqref{eq:Psop-46}
(except that now $a$ depends on the base variables from
both the left and the right factors of $\Snp$).
This takes a particularly simple form if $q$ is $\pi$-invariant, for
then, in the notation of \eqref{eq:Psop-46}, $q$ is independent of
both $Y$ and $\xi^a$. Thus, we can take $q$ outside the integral
in \eqref{eq:Psop-46}, i.e.\ it simply multiplies the indicial
operator of $\psi_0(H)$ by a constant.

\begin{lemma}\label{lemma:comm-2}
Suppose that $q\in\Cinf(\sct\Snp)$ is $\pi$-invariant and it
satisfies \eqref{eq:comm-4-b}. Let $A\in\PsiSc^{-\infty,0}(X,\calC)$
be as in the previous lemma. If $\zeta\in\sct(\Ct_a;X)$, then
$\Ah_a(\zeta)=q(\zeta)\widehat{\psi_0(H)}_a(\zeta)$.
\end{lemma}

Combining this lemma with Proposition~\ref{prop:WF-3} gives

\begin{cor}
Suppose that $\zeta\in\sct_{C'_a}(C_a;X)$ and $u\in\dist(X)$. If
$A$ is as in Lemma~\ref{lemma:comm-2}, $q(\zeta)\neq 0$, $Au\in\dCinf(X)$ and
$\zeta\nin\WFSc((H-\lambda)u)$ then $\zeta\nin\WFSc(u)$.
\end{cor}

Since the indicial operator of $[A,H]=AH-HA$ in $\PsiSc^{-\infty,0}(X,\calC)$
is just
\begin{equation}
\widehat{[A,H]}_{a,0}(\zeta)=[\Ah_{a,0}(\zeta),\Hh_{a,0}(\zeta)]
=q(\zeta)[\psi_0(\Hh_{a,0}(\zeta)),\Hh_{a,0}(\zeta)]=0
\end{equation}
for every $a$ and $\zeta\in\sct (\Ct_a;X)$, we see that for every
$A$ as in Lemma~\ref{lemma:comm-1}, $[A,H]\in\PsiSc^{-\infty,1}(X,\calC)$.
The additional order of decay corresponds to the one
in the scattering calculus. Moreover, the indicial operator of $[A,H]$
at $\Ct_0$,
as an operator in $\PsiSc^{-\infty,1}(X,\calC)$
(so this indicial operator is just a function on $\sct(\Ct_0;X)$), is given by
the Poisson bracket formula from the scattering calculus. Since
$V$ vanishes at $\Ct_0$, this gives
\begin{equation}\label{eq:comm-23}
\widehat{i[A,H]}_{1,0}=-\scHg (q\psi_0(g))=-\psi_0(g)\scHg q.
\end{equation}
If the indicial operators of $H$ at the other faces have no $L^2$
eigenfunctions, then this estimate combined with a compactness argument
suffice to prove an estimate for $[A,H]$ modulo lower operators
(i.e. modulo $\PsiSc^{-\infty,2}(X,\calC)$). However, to make the
compactness argument work, we need to estimate the indicial operators,
$\widehat{[A,H]}_{a,1}$, for all $a$. This is facilitated by the
following lemma.

\begin{lemma}\label{lemma:comm-4}
Let $q$ and $A$ be as in Lemma~\ref{lemma:comm-2}. For every seminorm
in
\begin{equation*}
\PsiSc^{-\infty,0}(\rho_a^{-1}(p),T_p\calC^a)
\end{equation*}
and for every $l\in\Nat$
there exist
$C>0$ and $m\in\Nat$ such that for every $a$ and every
$\zeta\in\sct_p(\Ct_a;X)$, $p\in\Ct_a$,
the seminorm
of $\widehat{[A,H]}_{a,1}(\zeta)$
in $\PsiSc^{-\infty,0}(\rho_a^{-1}(p),T_p\calC^a)$ is bounded by
\begin{equation}
C(|q(\zeta)|+\sum_{|\alpha|\leq m}\sup_{\xi^a}|\langle\xi^a\rangle^{-l}
(\partial^\alpha_{\xi^a} dq)
(\zeta,\xi^a)|)
\end{equation}
where the differential $dq$ is taken with respect to all variables,
in $\sct \Snp$, i.e.\ it is the differential of $q\in\Cinf(\sct \Snp)$.
\end{lemma}

\begin{rem}
Similar conclusions hold for every seminorm in
$\PsiScra^{-\infty,0}(\rho_a^*\sct(\Ct_a;X),\calCt_a)$,
which can be seen directly
from our calculations in the following proof.
\end{rem}

\begin{proof}
This can be proved directly from the definition of the
indicial operators, i.e.\ by computing $x^{-1}e^{-i\ft}[A,H]e^{i\ft}u'$
where $\ft\in\Cinf(X)$ and $u'\in\Cinf([X;\calC])$, similarly to
\cite[Sections~7,13]{Vasy:Propagation-2}.
Since this is equal to $x^{-1}[e^{-i\ft}Ae^{i\ft},e^{-i\ft}He^{i\ft}]u'$,
and $e^{-i\ft}Ae^{i\ft}\in\PsiSc^{-\infty,0}(X,\calC)$, we can assume
that $\ft=0$, the calculation being very similar in the general case.
To compute the commutator, it suffices to commute both $Av$ and $Hv$ for
every $v\in\Cinf([X;\calC])$ modulo terms that vanish with their first
derivatives in $\betaSc^*C_a$. A straight-forward
calculation can be performed
just as in \eqref{eq:Psop-40}-\eqref{eq:Psop-46},
where only the 0th order terms were kept. That shows with our coordinates
that
\begin{equation}\begin{split}\label{eq:comm-43}
\widehat{[A,H]}_{a,1}(\zeta)=&[\widehat{(\partial_x A)}_{a,0}(\zeta)
,\Hh_{a,0}(\zeta)]\\
&+(-(D_\tau q)(\zeta)([Y,\Hh_{a,0}(\zeta)]\partial_Y
+Y\partial_Y\Hh_{a,0}(\zeta))\\
&+(D_\nu q)(\zeta)(\partial_z\Hh_{a,0})(\zeta)-(\partial_z q)(\zeta)
(D_\nu \Hh_{a,0})(\zeta)\\
&+(\partial_\tau q)(\zeta)(\nu\cdot D_\nu \Hh_{a,0})(\zeta)
-(\nu\cdot D_\nu q)(\zeta)(\partial_\tau\Hh_{a,0}))\psi(\Hh_{a,0}(\zeta)).
\end{split}\end{equation}
Here $\partial_x A$ denotes the operator with kernel given by $\partial_x$
applied to that of $A$. Since in our notation the kernel of $A$ is
\begin{equation}
\int e^{i(w-w')\cdot\xi}q(w,\xi)p(w',\xi)\,d\xi,
\end{equation}
with the integral being convergent, rewriting this with the coordinates
on the compactification $[\Snp;\calC_a]$, \eqref{eq:ff-coords}, so
that $q$ takes the form $q(x,xY,z,\xi)$ proves
that all terms of \eqref{eq:comm-43} satisfy the stated estimate,
completing the proof.

Another approach to compute $a$-indicial operators
is to use that near $C'_a$, $A$ can be regarded as a
(non-classical!) pseudo-differential operator in the free variables
$(w_a,\xi_a)$ with values in bounded operators on
$L^2(X_a)$ (in fact, with values in $\PsiSc^{-\infty,0}(\Xb^a,\calC^a)$).
More precisely, $A\in\Psiscc^{-\infty,0}(\Xb_a;\bop(L^2(X^a),L^2(X^a)))$.
This allows us to use the scattering calculus for the computation
of the commutators to give the stated result.
\end{proof}

As an application of these estimates, we now show how, under the assumption
that the subsystems have no bound states, a positive Poisson bracket with
$g$ can give rise of a positive operator estimate.
We thus assume that
\begin{equation}\label{eq:hypo-1a}
\Hh_{a,0}(\xi)\ \text{has no}\ L^2\ \text{eigenvalues for any}
\ a\neq 0\Mand \zeta\in\sct(\Ct_a;X).
\end{equation}

To simplify the notation in the following proposition, we introduce the
notation $\supp_a e\subset\sct(\Ct_a;X)$
for $\pi$-invariant functions $e\in\Cinf(\sct_{\bX} X)$. This is defined
as the support of the function on $\sct(\Ct_a;X)$ induced by $e$. Indeed,
as $e$ is $\pi$-invariant, its restriction to $\sct_{C_a}X$
can be regarded as a function
on $\sct (C_a;X)$. Then $\supp_a e$ is the support of the pull-back
of this function to $\sct (\Ct_a;X)$.

\begin{prop}\label{prop:comm-1}
Suppose that $H$ is a many-body Hamiltonian satisfying \eqref{eq:hypo-1a},
and $\lambda\in\Real$.
Suppose also that $q,b,e\in\Cinf(\sct X;\Real)$ are $\pi$-invariant,
satisfy
the bounds \eqref{eq:comm-4-b}, $q,b\geq 0$,
and that there
exist $\delta>0$, $C>0$,
$C_\alpha>0$, such that for all $\xi\in\sct_{\bX}X$,
\begin{equation}\label{eq:comm-35}
|g(\xi)-\lambda|<\delta\Rightarrow
\scHg q(\xi)\leq -b(\xi)^2+e(\xi)
\end{equation}
and
\begin{equation}\label{eq:comm-37}
|g(\xi)-\lambda|<\delta\Mand \xi\nin\supp e\Rightarrow
q(\xi)\leq Cb(\xi)^2\Mand|(\partial_\mu^\alpha dq)(\xi)|\leq C_\alpha
b(\xi)^2.
\end{equation}
Let $A\in\PsiSc^{-\infty,0}(X,\calC)$ be as in Lemma~\ref{lemma:comm-1}.
For any $\ep'>0$, $a\in I$ and for any
$K_a\subset\sct (\Ct_a;X)$ compact with $\supp_a e\cap K=\emptyset$
there exists $\delta'>0$
such that
if $\psi\in\Cinf_c(\Real)$ is supported in $(\lambda-\delta',\lambda+\delta')$
and $\zeta\in K_a$ then
\begin{equation}\label{eq:ind-est-1}
i\widehat{(\psi(H)[A^*A,H]\psi(H))}_{a,1}(\zeta)
\geq (2-\ep')b^2 q\psi(\Hh_{a,0}(\zeta))^2.
\end{equation}
\end{prop}

\begin{proof}
Note that the estimate \eqref{eq:ind-est-1} is trivial
if $\tau_a^2+|\nu_a|^2_{z_a}>\lambda+1$ (with $\zeta=(z_a,\tau_a,\nu_a)$,
$\delta'<1$ arbitrary) since
then both sides vanish as
\begin{equation}
\widehat{\psi(H)}_a(\zeta)=\psi(h_a(z_a)+\tau_a^2+|\nu_a|^2_{z_a}),
\end{equation}
$h_a$ denoting the subsystem Hamiltonian as in \eqref{eq:h_a(p)-def},
and $H_a\geq 0$ by the
assumption on the absence of bound states of {\em all} subsystem Hamiltonians
(including $H_c$ with $C_a\subset C_c$).

We prove \eqref{eq:ind-est-1} by induction on $a$.
First, \eqref{eq:ind-est-1} is
certainly satisfied for $a=0$. In fact,
as
$A\in\Psisc^{-\infty,0}(X,\calC)$, we can use the commutator
formula in the scattering calculus, \eqref{eq:comm-23}, to find
$\widehat{[A,H]}_{0,1}$.
Since $V$ vanishes at
the free face, $\betaSc^* C_0$, it does not contribute to
$\widehat{[A,H]}_{0,1}$, so we indeed have, by \eqref{eq:comm-35},
\begin{equation}
i\widehat{\psi(H)[A^*A,H]\psi(H)}_{0,1}
=-2q(\scHg q)\psi(g)^2\geq 2 b^2 q\psi(g)^2=2b^2 q\psi(\Hh_{0,0})^2
\end{equation}
away from $\supp_0 e$ under the assumption that
\begin{equation}
\supp\psi\subset(\lambda-\delta,\lambda+\delta).
\end{equation}

So suppose now that \eqref{eq:ind-est-1} has been proved for all $c$ with
$C_a\subset C_c$, $C_a\neq C_c$. This implies that all indicial operators of
$i\widehat{[\psi(H)A^*A\psi(H),H]}_a(\zeta)$, $\zeta=(z_a,\tau_a,\nu_a)\in K_a$
satisfy an inequality like
\eqref{eq:ind-est-1}. In fact, the indicial operators are of the
form $i\widehat{[\psi(H)A^*A\psi(H),H]}_c(\zetat)$ with $\betaSc(\zetat)
=(0,z_a)\in C_a$, $\pit_{ca}(\zetat)=\zeta$.
Such a $\zetat$ is of the form $\zetat=
(\Yh''_a,z_a,\tau_a,
\mu''_a,\nu_a)$ where $C_c$ is given by $x=0$, $y'=0$, so $(\Yh''_a,z_a)$
give coordinates along $\Ct_c$. Note that as $K_a$
is compact, so is
\begin{equation}
K_c=\{\zetat=(\Yh''_a,z_a,\tau_a,\mu''_a,\nu_a):\ (z_a,\tau_a,\nu_a)
\in K_a,\ \betaSc(\zetat)\in C_a,\ |\mu''_a|\leq\lambda+1\}
\end{equation}
and as $e$ is independent of $\mu''_a$ at
$C_a$, $K_c\cap\supp_c e=\emptyset$, so we can apply the inductive
hypothesis. Taking into account that the estimate \eqref{eq:ind-est-1} is
trivial at $C_c$ for $\zetat$ with $|\mu''_a|>\lambda+1$, we see that
for all $\zetat=(0,z_a,\tau_a,\mu''_a,\nu_a)$ with $(z_a,\tau_a,\nu_a)\in K_a$,
we have
\begin{equation}\label{eq:ind-est-1c}
i\widehat{(\psi(H)[A^*A,H]\psi(H))}_{c,1}(\zetat)
\geq (2-\ep')b^2 q\widehat{\psi(H)^2}_{c,0}(\zetat).
\end{equation}
Since $b^2q$ is $\pi$-invariant
on $\sct X$, it is independent of $\zetat$ for each fixed $\zeta$, and
if it vanishes at $\zeta$, then so does
$\widehat{[\psi(H)A^*A\psi(H),H]}_{a,1}(\zeta)$ by
Lemmas~\ref{lemma:comm-2}-\ref{lemma:comm-4} and \eqref{eq:comm-37}.
Thus, by Proposition~\ref{prop:str-pos-ops},
\begin{equation}
i\widehat{[\psi(H)A^*A\psi(H),H]}_{a,1}(\zeta)\geq (2-\ep')
b^2 q\widehat{\psi(H)}_{a,0}^2(\zeta)
+R(\zeta)
\end{equation}
where the seminorms of
\begin{equation*}
R(\zeta)\in\PsiSc^{-\infty,1}(\rho_a^{-1}(p),T_p\calC^a),\quad
\zeta\in\sct_p(\Ct_a;X),
\end{equation*}
are bounded by those of $\widehat{[\psi(H)A^*A\psi(H),H]}_{a,1}(\zeta)$
and by $b(\zeta)^2 q(\zeta)$. By assumption \eqref{eq:comm-37}
and Lemma~\ref{lemma:comm-4} the former
are bounded by the latter, so $R(\zeta)$ satisfies
the estimate
\begin{equation}\label{eq:p-of-s-4}
\|R(\zeta)\|
_{\bop(L^2_\scl(\rho_a^{-1}(p)),\Hsc^{1,1}(\rho_a^{-1}(p)))}
\leq C''q(\zeta) b(\zeta)^2
\end{equation}
with $C''$ independent of $q$ and $b$.

We now use our hypothesis on the absence of bound states.
So suppose that $\psi_1,\psi_2\in\Cinf_c(\Real)$,
$\psi\equiv 1$ near $\supp\psi_1$, $\psi_1\equiv 1$ near $\supp\psi_2$.
By assumption, $\lambda-\tau_a^2-|\nu_a|^2_{z_a}$
is not an eigenvalue of the subsystem Hamiltonian, $h_a(z)$.
Thus,
\begin{equation}\label{eq:psit(H)-to-0}
\psi_1(\Hh_a(\zeta))=\psi_1(h_a(z)+\tau_a^2+|\nu_a|_{z_a}^2)\to 0
\end{equation}
strongly as $\supp\psi_1\to\{\lambda\}$.
Since $K_a$ is compact,
and the inclusion map
\begin{equation}
T:\Hsc^{1,1}(\rho_a^{-1}(p))\hookrightarrow
L^2_{\scl}(\rho_a^{-1}(p))
\end{equation}
is compact,
for $\psi_1$ with sufficiently small support we have
\begin{equation}
\|\widehat{(\psi_1(H)T)}_a(\zeta)\|_{\bop(\Hsc^{1,1}(\rho_a^{-1}(p)),
L^2_{\scl}(\rho_a^{-1}(p)))}\leq\epsilon' (C'')^{-1}
\end{equation}
for all $\zeta\in K_a$.
Thus,
\begin{equation}\label{eq:p-of-s-6}
\widehat{i(\psi_1(H)[A^*A,H]\psi_1(H))}_{a,1}(\zeta)
\geq (2-\ep')b^2 q\widehat{\psi_1(H)}_{a,0}^2(\zeta)-\epsilon'b^2 q,
\ \zeta\in K_a.
\end{equation}
Multiplying by $\psi_2(H)$ from both left and right we finally conclude that
\begin{equation}\label{eq:p-of-s-7}
\widehat{i(\psi_2(H)[A^*A,H]\psi_2(H))}_{a,1}
\geq(2-2\epsilon')b^2 q\widehat{\psi_2(H)}_{a,0}^2.
\end{equation}
Relabelling $\psi_2$ and $2\ep'$ as $\psi$ and $\ep'$ (thereby putting
stronger restrictions on $\supp\psi$) provides the inductive step
and completes the proof of \eqref{eq:ind-est-1}.
\end{proof}

In the following corollary we add an extra term to the commutator that
will enable us to deal with other terms arising in the propagation
estimates.

\begin{cor}\label{cor:comm-7}
Suppose that the assumptions of Proposition~\ref{prop:comm-1} are
satisfied and let $C$ be as in \eqref{eq:comm-37}.
Suppose in addition that for any differential operator $Q$ on
$\sct(\Ct_a;X)$ and multiindex $\alpha$ there exist constant
$C_Q$ and $C_{\alpha,Q}$ such that
\begin{equation}\begin{split}\label{eq:comm-37-p}
|g(\xi)-\lambda|<\delta,&\ b(\xi)\neq 0 \Mand \xi\nin\supp e\\
&\Rightarrow
|Q(b^{-2} q)(\xi)|\leq C_Q\Mand|Q(b^{-2}(\partial_\mu^\alpha dq))(\xi)|\leq C_{\alpha,Q}.
\end{split}\end{equation}
For any $\ep'>0$, $M>0$, and for any
$K\subset\scdt X$ compact with $\supp e\cap K=\emptyset$
there exists $\delta'>0$,
$B,E\in\PsiSc^{-\infty,0}(X,\calC)$, $F\in\PsiSc^{-\infty,1}(X,\calC)$
with
\begin{equation}\label{eq:comm-51}
\WFScp(E)\cap K=\emptyset,\ \WFScp(F)\subset\supp q,\ \Bh_{a,0}(\zeta)=
b(\zeta)q(\zeta)^{1/2}
\psi(\Hh_{a,0}(\zeta)),\ \zeta\in K,
\end{equation}
such that
if $\psi\in\Cinf_c(\Real)$ is supported in $(\lambda-\delta',\lambda+\delta')$
then
\begin{equation}\label{eq:comm-53}
i\psi(H)x^{-1/2}[A^*A,H]x^{-1/2}\psi(H)-M\psi(H)A^*A\psi(H)
\geq (2-\ep'-MC)B^*B+E+F.
\end{equation}
\end{cor}

\begin{proof}
Let $p\in\Cinf(\sct X)$ be $\pi$-invariant, $p\geq 0$, satisfy estimates
\eqref{eq:comm-4-b}, and such that $\supp p\cap \supp e=\emptyset$
and $\supp(1-p)\cap K=\emptyset$. (Here
$p$ can be regarded as a function on $\scdt X$.) Let
$\psi_0\in\Cinf_c(\Real;[0,1])$ identically $1$ near
$[\lambda-\delta,\lambda+\delta]$, and let
$P\in\PsiSc^{-\infty,0}(X,
\calC)$ be such that $\Ph_{a,0}(\zeta)=p(\zeta)\psi_0(\Hh_a(\zeta))$
and $\WFScp(\psi_0(H)-P)\cap K=\emptyset$.
For example, $P$ can be constructed as in Lemma~\ref{lemma:comm-1}.

The indicial operators of
\begin{equation*}
i\psi(H)P^*x^{-1/2}[A^*A,H]x^{-1/2}P\psi(H)-M\psi(H)A^*A\psi(H)
\end{equation*}
are
\begin{equation}\begin{split}
&i\widehat{\psi(H)P^*x^{-1/2}[A^*A,H]x^{-1/2}P\psi(H)}_{a,0}(\zeta)
-M\widehat{(\psi(H)P^*A^*AP\psi(H))}_{a,0}(\zeta)\\
&\hspace*{2 cm}=ip(\zeta)^2\widehat{\psi(H)[A^*A,H]\psi(H)}_{a,1}(\zeta)
-Mq(\zeta)^2p(\zeta)^2
\psi(\Hh_{a,0}(\zeta))
\end{split}\end{equation}
since $\psi_0\psi=\psi$. Thus, by Proposition~\ref{prop:comm-1} and as
$Mq\leq MCb^2$, we have
\begin{equation}\begin{split}
&i\widehat{\psi(H)P^*x^{-1/2}[A^*A,H]x^{-1/2}P\psi(H)}_{a,0}(\zeta)
-M\widehat{(\psi(H)P^*A^*AP\psi(H))}_{a,0}(\zeta)\\
&\hspace*{5 cm}\geq(2-\ep'-MC)b^2 q\psi(\Hh_{a,0}(\zeta))^2.
\end{split}\end{equation}
Thus, taking into account \eqref{eq:comm-37-p} and the remark following
Lemma~\ref{lemma:comm-4},
Proposition~\ref{prop:pos-ops} gives
\begin{equation}\label{eq:comm-59}
i\psi(H)P^*x^{-1/2}[A^*A,H]x^{-1/2}P\psi(H)-M\psi(H)A^*A\psi(H)
\geq (2-\ep')B^*B+F,
\end{equation}
with $B\in\PsiSc^{-\infty,0}(X,\calC)$, $F\in\PsiSc^{-\infty,1}(X,\calC)$,
\begin{equation}
\Bh_{a,0}(\zeta)=p(\zeta)b(\zeta)q(\zeta)^{1/2},
\end{equation}
so the second statement of \eqref{eq:comm-51} holds. Moreover,
writing $\psi(H)=P\psi(H)+(\psi_0(H)-P)\psi(H)$, and expanding
the left hand side of \eqref{eq:comm-53}, every term but the one given
in \eqref{eq:comm-59} has operator wave front set disjoint from $K$.
Letting $E$ be the sum of these terms proves the corollary.
\end{proof}

\section{Propagation of singularities}
\label{sec:propagation}
In this section we prove that singularities of generalized eigenfunctions
of the many-body operator $H$ propagate along generalized broken
bicharacteristics under the assumption that
that no (proper) subsystems of $H$ have a bound state. That is, due to
our definition in Section~\ref{sec:Hamiltonian}, we assume that
\begin{equation}\label{eq:hypo-1}
\Hh_{b,0}(\xi)\ \text{has no}\ L^2\ \text{eigenvalues for any}
\ b\neq 0\Mand \xi\in\sct(\Ct_b;X).
\end{equation}
The technical
reason for this assumption lies in
the argument of Proposition~\ref{prop:comm-1} in which a symbolic estimate
is used to deduce positivity estimates for the indicial operators. However, it
is clear that the generalized broken bicharacteristics of $\Delta-\lambda$
cannot be expected to describe propagation if the subsystems have bound
states since in this situation even the characteristic set of $H$
(i.e.\ the set where $\Hh_{b,0}(\zeta)$ is not invertible) changes.

Suppose that $p\in C'_a=C'$ (the regular part of
$C$). As in Section~\ref{sec:Hamiltonian}, let $(x,y,z)=(x_a,y_a,z_a)$
be coordinates on $X$ near $p$ with $x$ defining $\bX$ as
usual, $C$ defined by $x=0$, $y=0$, chosen so that
every $C_b$ with $p\in C_b$ is a
product-linear
submanifold of $\bX$ in these local coordinates, i.e.\ it is of the
form $\{(y,z):\ A_b y=0\}$ where $A=A_b$ is a matrix.
In addition, as in Section~\ref{sec:Hamiltonian}, we arrange
that
at $C$, $\partial_{y_j}=\partial_{(y_a)_j}$
is perpendicular to $T C$ for each $j$ (with respect to $h$) and they are
orthonormal with respect to each other at $p$.
Let
$(\tau,\mu,\nu)=(\tau_a,\mu_a,z_a)$
denote the $\scl$-dual variables, so we write elements
of $\sct X$ as
\begin{equation}
\tau\frac{dx}{x^2}+\mu\cdot\frac{dy}{x}+\nu\cdot\frac{dz}{x}.
\end{equation}
Thus, at $p$ (i.e.\ on $\sct_p X$)
the metric function of $h$ is of the form $|\mu|^2+\htil(z,\nu)$
with $|\mu|$ denoting the Euclidean length of $\mu$ and $\htil$ is the
metric function of the restriction of $h$ to $TC$.
When talking about $C_b$, we sometimes write
the corresponding orthogonal splitting of $y$ as $y=(y',y'')$, so $C_b$
is defined by $A_b y=y'=0$ in $\bX$.

Recall that $\pi_{0a}:\sct_C X\to \sct(C;X)$
is the (orthogonal) projection given by the metric at $C$. Thus, in
our local coordinates $(y,z,\tau,\mu,\nu)$ on $\sct_{\bX}X$,
$\pi_{0a}(0,z,\tau,\mu,\nu)=(z,\tau,\nu)$.
We use composition with the projection $\sct_{\bX} X$ to $\sct_C X$
given by our choice of local coordinates,
$(y,z,\tau,\mu,\nu)\mapsto (z,\tau,\mu,\nu)$,
to extend $\pi_{0a}$ to a map,
denoted by $\pi^e_{0a}$, from $\sct_{\bX}X$ to
$\sct(C;X)$.
Thus, $\pi_{0a}^e(y,z,\tau,\mu,\nu)=(z,\tau,\nu)$.

The propagation of singularities estimate
in directions tangential to $C$ proceeds much as in the 3-body case.
In fact, essentially the same operator as there gives a positive commutator,
see Propositions~\ref{prop:tot-geod-tgt-prop}-\ref{prop:tgt-prop};
the functional analysis part of the argument is much as in the normal
case which we proceed to examine.
Recall that the normal part of the characteristic set of $H-\lambda$
over $C'$ is
\begin{equation}
\Sigma_n(\lambda)\cap\sct_{C'}(C;X)=\{(z,\tau,\nu):\ \tau^2
+\htil(z,\nu)<\lambda\}.
\end{equation}
Since the characteristic set $\Sigma_{\Delta-\lambda}$
of $\Delta-\lambda$ is given by $\tau^2+|\nu|^2_z
+|\mu|^2=\lambda$ at $p$, the condition
$\pi(\xit)\in\Sigma_n(\lambda)\cap\sct_p(C;X)$,
$\xit\in\Sigma_{\Delta-\lambda}$
implies that $\mu\neq 0$. Since the rescaled Hamilton vector field
$\scHg$ of $\Delta$ (restricted to $\sct_{\bX}X$) is given by
\begin{equation}
\scHg=2\tau(\mu\cdot\partial_\mu+\nu\cdot\partial_\nu)-2h
\partial_\tau+H_h,
\end{equation}
the $\partial_y$ component of $\scHg$ at $p$ is $2\mu\cdot\partial_y$, meaning
that bicharacteristics of $\Delta$ through $\xit$ are normal to
$\sct_C X$. In addition, with $\eta=y\cdot\mu$,
$\eta$ is $\pi$-invariant and can be used to parameterize bicharacteristic
curves near $\xi=\pi(\xit)$.
In fact, at each $C_b$ with $p\in C_b$, $\eta=\mu\cdot y$
has the property
that if we split $y=(y',y'')$ so that $x=0$, $y'=0$ defines $C_b$ then
$\mu\cdot y=\mu'\cdot y'+\mu''\cdot y''$ is independent of $\mu'$ at $y'=0$,
so $\eta$ is $\pi$-invariant.
Moreover, $\scHg\eta(\xit)=2|\mu|^2>0$,
so $\eta$ can be used to parametrize the generalized broken bicharacteristics
near $\xi$ as claimed. We remark that $\tau$ is another possible variable
to use for the parameterization, as usual.

We now proceed to prove two normal propagation estimates. The first
one will be less precise, but it works under our most general assumptions.
On the other hand, the second estimate requires that all elements
of $\calC$ be totally geodesic, but it locates the incoming
singularities more precisely. Although the consequences are the same,
as far as propagation along generalized broken bicharacteristics
is concerned (due to the geometry of these bicharacteristics),
the finer estimate is worth proving since it
is closer to the tangential estimates in spirit
and it applies in the setting of
most interest, Euclidean many-body scattering.

We only state the following propagation result for propagation in the
forward direction along the generalized broken bicharacteristics.
A similar result holds in the backward direction, i.e.\ if we replace
$\eta(\xi)<0$ by $\eta(\xi)>0$
in \eqref{eq:prop-9a};
the proof in this case only requires changes in some signs
in the argument given below.

\begin{prop}\label{prop:normal-prop}
Suppose that $H$ is a many-body Hamiltonian satisfying \eqref{eq:hypo-1}.
Let $u\in\dist(X)$,
$\lambda>0$.
Let $\xi_0=(z_0,\tau_0,\nu_0)\in\Sigma_n(\lambda)\cap\sct_{C'}(C;X)$
and let $\eta=y\cdot\mu$ be the $\pi$-invariant function defined
in the local coordinates discussed above.
If there exists a neighborhood $U$ of $\xi_0$ in $\dot\Sigma$ such that
\begin{equation}\begin{split}\label{eq:prop-9a}
\xi\in U\Mand \eta(\xi)<0\Rightarrow\xi\nin\WFSc(u)
\end{split}\end{equation}
then $\xi_0\nin\WFSc(u)$.
\end{prop}

\begin{rem}
Note that $\eta(\xi)<0$ implies $y\neq 0$, so $\xi\nin\sct_{C'}(C;X)$.
\end{rem}

\begin{proof}
The main step in the proof is the construction of an operator which has
a microlocally positive commutator with $H$ near $\xi_0$. In fact, we
construct the symbol of this operator. This symbol will not be a
scattering symbol, i.e.\ it will not be in $\Cinf(\Snp\times\Snp)$,
only due to its behavior as $\mu\to\infty$ corresponding to its
$\pi$-invariance. This
will be accommodated by composing its quantization with a cutoff in the
spectrum of $H$, $\phi(H)$, $\phi\in\Cinf_c(\Real)$ supported near $\lambda$,
as discussed in Lemma~\ref{lemma:comm-1}. This approach simply
extends the one taken in
the three-body scattering proof of \cite{Vasy:Propagation-2}, though the
actual construction is different due to the more complicated geometry.

Employing an iterative argument as usual, we may assume
that $\xi_0\nin\WFSc^{*,l}(u)$ and we need
to show that $\xi_0\nin\WFSc^{*,l+1/2}(u)$.

First we define a distance function to $\xi_0$. Thus, we let
\begin{equation}
\omega=|y|^2+|z-z_0|^2+|\tau-\tau_0|^2+|\nu-\nu_0|^2,
\end{equation}
$|.|$ denoting the Euclidean norm.
Then $\omega$ vanishes quadratically at $\xi_0$, so $|d\omega|\leq C'_1
\omega^{1/2}$. In particular,
\begin{equation}\label{eq:omega-est-a}
|\scHg\omega|\leq C_1\omega^{1/2}.
\end{equation}

Next, we use the variable $\eta=y\cdot\mu$ to measure propagation.
Let
\begin{equation}
c_0=\lambda-\tau_0^2-|\nu_0|^2_{z_0}>0.
\end{equation}
Since the $\partial_y$ component of $\scHg$ at $(0,z_0,\tau,\mu,\nu)$ is
$2\mu$, we see that
\begin{equation}
|\scHg \eta-2|\mu|^2|\leq C'_2(|y|+|z-z_0|)\leq C_2\omega^{1/2}.
\end{equation}
In addition,
\begin{equation}\begin{split}\label{eq:n-prop-30aa}
|\lambda&-\tau_0^2-|\nu_0|^2_{z_0}-|\mu|^2|
\leq|\lambda-g|+|g-\tau_0^2-|\nu_0|^2_{z_0}-|\mu|^2|\\
&\leq|\lambda-g|
+C'(|y|+|z-z_0|+|\tau-\tau_0|+|\nu-\nu_0|)
\leq|\lambda-g|+C_3\omega^{1/2}
\end{split}\end{equation}
so we conclude that
\begin{equation}\label{eq:eta-est-a}
|\scHg\eta-2c_0|\leq C_4(|\lambda-g|+\omega^{1/2}).
\end{equation}

For $\beta>0$, $\delta>0$, with other restrictions to be imposed later on,
let
\begin{equation}
\phi=\eta+\frac{\beta}{\delta}\omega,
\end{equation}
so $\phi$ is a $\pi$-invariant function.
Let $\chi_0\in\Cinf(\Real)$ be equal to $0$ on $(-\infty,0]$ and
$\chi_0(t)=\exp(-1/t)$ for $t>0$. Thus, $\chi_0'(t)=t^{-2}\chi_0(t)$.
Let $\chi_1\in\Cinf(\Real)$ be $0$
on $(-\infty,0]$, $1$ on $[1,\infty)$, with $\chi_1'\geq 0$ satisfying
$\chi_1'\in\Cinf_c((0,1))$.
Furthermore, for $A_0>0$ large, to be determined, let
\begin{equation}\label{eq:prop-22}
q=\chi_0(A_0^{-1}(2-\phi/\delta))\chi_1(y\cdot\mu/
\delta+2).
\end{equation}
Thus, on $\supp q$ we have $\phi\leq 2\delta$ and $y\cdot\mu\geq-2\delta$.
Since $\omega\geq 0$, the first of these inequalities implies that
$y\cdot\mu\leq 2\delta$, so on $\supp q$
\begin{equation}
|y\cdot\mu|\leq 2\delta.
\end{equation}
Hence,
\begin{equation}
\omega\leq (\delta/\beta)(2\delta-y\cdot\mu)\leq4\delta^2\beta^{-1}.
\end{equation}
We now proceed to estimate $\scHg\phi$. First, by \eqref{eq:eta-est-a}
and \eqref{eq:omega-est-a},
\begin{equation}
|\scHg\phi-2c_0|<C_4(|\lambda-g|+\omega^{1/2})
+\frac{C_1 \beta}{\delta}\omega^{1/2}.
\end{equation}
So let
\begin{equation}
\beta=\frac{c_0^2}{(8C_1)^2}\Mand \delta_0=\frac{c_0\sqrt\beta}{8 C_4}.
\end{equation}
Under the additional assumptions
\begin{equation}
\delta<\delta_0
\Mand |\lambda-g|<\frac{c_0}{4C_4}
\end{equation}
we have $\omega^{1/2}\leq c_0/(4C_4)$, so
we conclude that
$|\scHg\phi-2c_0|\leq c_0$, hence
\begin{equation}
\scHg \phi\geq c_0>0.
\end{equation}

This at once gives a positivity estimate for $\scHg q$ near $\xi_0$.
Namely,
\begin{equation}\begin{split}\label{eq:scHg q-calc}
\scHg q=-A_0^{-1}&\delta^{-1}\chi'_0(A_0^{-1}(2-\phi/\delta))
\chi_1(y\cdot\mu/
\delta+2)\scHg \phi\\
&+\delta^{-1}
\chi_0(A_0^{-1}(2-\phi/\delta))\chi_1'(y\cdot\mu/
\delta+2)\scHg \eta.
\end{split}\end{equation}
Thus,
\begin{equation}\label{eq:scHg q-pos}
\scHg q=-\bt^2+e
\end{equation}
with
\begin{equation}
\bt^2=A_0^{-1}\delta^{-1}\chi'_0(A_0^{-1}(2-\phi/\delta))
\chi_1(y\cdot\mu/\delta+2)\scHg \phi.
\end{equation}
Hence, with
\begin{equation}
b^2=c_0 A_0^{-1}\delta^{-1}\chi'_0(A_0^{-1}(2-\phi/\delta))
\chi_1(y\cdot\mu/\delta+2),
\end{equation}
we have
\begin{equation}
\scHg q\leq -b^2+e.
\end{equation}
Moreover,
\begin{equation}\label{eq:prop-33}
b^2\geq (c_0 A_0/16) q
\end{equation}
since $\phi\geq y\cdot\mu\geq -2\delta$ on $\supp q$, so
\begin{equation}\begin{split}\label{eq:prop-34}
\chi'_0(A_0^{-1}(2-\phi/\delta))
&=A_0^2(2-\phi/\delta)^{-2}\chi_0(A_0^{-1}(2-\phi/\delta))\\
&\geq (A_0^2/16)
\chi_0(A_0^{-1}(2-\phi/\delta)).
\end{split}\end{equation}
On the other hand, $e$ is supported where
\begin{equation}\label{eq:supp-e-a}
-2\delta\leq y\cdot\mu \leq-\delta,\quad \omega^{1/2}\leq2\beta^{-1/2}
\,\delta,
\end{equation}
so,
for $\delta>0$ sufficiently small,
in the region which we know is disjoint from $\WFSc(u)$.
Moreover, on $\supp q$,
\begin{equation}\label{eq:supp-q-a}
-2\delta\leq y\cdot\mu \leq 2\delta,\ \omega^{1/2}\leq2\beta^{-1/2}
\,\delta,
\end{equation}
so, for $\delta>0$ sufficiently small, we deduce from the
inductive hypothesis that $\supp q$ (hence $\supp b$) is disjoint
from $\WFSc^{*,l+1/2}(u)$.
In addition,
by choosing $\delta>0$ sufficiently small, we can assume that the support
of $q$, $e$ and $b$ are all disjoint from $\WFSc((H-\lambda)u)$.

Moreover, with $\partial$ denoting a partial derivative with respect to
one of $(y,z,\tau,\mu,\nu)$,
\begin{equation}\begin{split}
\partial q=-A_0^{-1}&\delta^{-1}\chi'_0(A_0^{-1}(2-\phi/\delta))
\chi_1(\eta/
\delta+2)\partial \phi\\
&-\delta^{-1}
\chi_0(A_0^{-1}(2-\phi/\delta))\chi_1'(\eta/
\delta+2)\partial \eta.
\end{split}\end{equation}
As $y=0$ is outside the support of the second term, and as $\partial_\mu\phi$
vanishes at $y=0$, we conclude that for any multiindex $\beta$,
\begin{equation}\label{eq:prop-38}
|\partial_\mu^\beta dq|\leq C_\beta b^2\ \text{at}\ y=0.
\end{equation}
More generally, at any $C_b$ with $p\in C_b$, defined by $x=0$, $y'=0$, as
above, $\phi$ is independent of $\mu'$ at $y'=0$ so outside $\supp e$
\begin{equation}
|\partial_{\mu'}^\beta dq|\leq C_\beta b^2\ \text{at}\ y'=0.
\end{equation}
In fact, outside $\supp e$, but in the set where $b$ is positive,
\begin{equation}\label{eq:prop-39}
b^{-2}\partial q=c_0^{-1}\partial\phi,
\end{equation}
so the uniform bounds of \eqref{eq:comm-37-p} also follow.

Let $\psit\in\Cinf_c(\Real)$ be identically $1$ near $0$ and supported
sufficiently close to $0$ so that the product decomposition of $X$ near
$\bX$ is valid on $\supp\psit$. We also define
\begin{equation}
\qt=\psit(x)q.
\end{equation}
Thus, $\qt\in\Cinf(\sct X)$ is a $\pi$-invariant function satisfying
\eqref{eq:comm-4-b}. Let $A$ be the operator given by Lemma~\ref{lemma:comm-1}
with $\qt$ in place of $q$,
so in particular its indicial operators are $q(\zeta)\psi_0(\Hh_{b,0}(\zeta))$.
Note that \eqref{eq:comm-37} holds with $C=16 c_0^{-1}A_0^{-1}$. So suppose
that $M>0$ and $\ep'>0$. Choose $A_0$ so large that $MC<\ep'$.
By Corollary~\ref{cor:comm-7} and the hypothesis \eqref{eq:hypo-1},
we deduce the following statement.
For any
$K'\subset\scdt X$ compact with $\supp e\cap K'=\emptyset$
there exists $\delta'>0$,
$B,E\in\PsiSc^{-\infty,0}(X,\calC)$, $F\in\PsiSc^{-\infty,1}(X,\calC)$
with
\begin{equation}\label{eq:comm-51a}
\WFScp(E)\cap K'=\emptyset,\ \WFScp(F)\subset\supp \qt,\ \Bh_{a,0}(\zeta)=
b(\zeta)q(\zeta)^{1/2}
\psi(\Hh_{a,0}(\zeta)),\ \zeta\in K',
\end{equation}
such that
if $\psi\in\Cinf_c(\Real)$ is supported in $(\lambda-\delta',\lambda+\delta')$
then
\begin{equation}\label{eq:comm-53a}
i\psi(H)x^{-1/2}[A^*A,H]x^{-1/2}\psi(H)-M\psi(H)A^*A\psi(H)
\geq (2-2\ep')B^*B+E+F.
\end{equation}

Let
\begin{equation}
\Lambda_r=x^{-l-1/2}(1+r/x)^{-1},\quad r\in(0,1),
\end{equation}
so $\Lambda_r\in\PsiSc^{0,-l+1/2}(X,\calC)$ for $r\in(0,1)$ and it is
uniformly bounded in $\PsiSc^{0,-l-1/2}(X,\calC)$. The last statement
follows from $(1+r/x)^{-1}$ being uniformly bounded as a 0th order
symbol, i.e.\ from $(x\partial_x)^k(1+r/x)^{-1}\leq C_k$ uniformly
($C_k$ independent of $r$).
We also define
\begin{equation}\label{eq:prop-57}
A_r=A\Lambda_r x^{-1/2}\psi(H),\ B_r=B\Lambda_r,\ E_r=\Lambda_r E\Lambda_r.
\end{equation}
Then, with $\psi_0\in\Cinf_c(\Real;[0,1])$ identically $1$ near $\supp\psi$,
\begin{equation}\begin{split}
ix^{l+1/2}&[A^*_rA_r,H]x^{l+1/2}\\
&=i(1+r/x)^{-1}\psi(H)x^{-1/2}[A^*A,H]x^{-1/2}\psi(H)(1+r/x)^{-1}\\
&\qquad+i\psi(H)A^*x^{l+1/2}[\Lambda_r x^{-1/2},H](1+r/x)^{-1}x^{-1/2}\psi_0(H)
A\psi(H)\\
&\qquad
+i\psi(H)A^*\psi_0(H)
x^{-1/2}(1+r/x)^{-1}[\Lambda_r x^{-1/2},H]x^{l+1/2}A\psi(H)
+H_r,
\end{split}\end{equation}
where $H_r$ is uniformly bounded in $\PsiSc^{-\infty,1}(X,\calC)$.
Note that $H_r$ arises by commuting $A$, powers of $x$ and $\Lambda_r$ through
other operators, but as the indicial operators of $A$ and $x$
are a multiple of the identity, $A$, $x$ and $\Lambda_r$ commute with these
operators to top order, and in case of $\Lambda_r$, the commutator
is uniformly bounded as an operator of one lower order.
Then, multiplying \eqref{eq:comm-53a} by $(1+r/x)^{-1}$ from the left
and right and rearranging the terms we obtain the following estimate of
bounded self-adjoint operators on $L^2_\scl(X)$:
\begin{equation}\label{eq:prop-61}\begin{split}
ix^{l+1/2}[A^*_rA_r&,H]x^{l+1/2}
-\psi(H)A^*(G_r^*+G_r)A\psi(H)-M\psi(H)A^*A\psi(H)\\
&\geq x^{l+1/2}((2-\epsilon')B_r^*B_r+E_r+F_r)x^{l+1/2}
\end{split}\end{equation}
where
\begin{equation}
G_r=i\psi_0(H)x^{-1/2}(1+r/x)^{-1}[\Lambda_r x^{-1/2},H]x^{l+1/2},
\end{equation}
and $F_r\in\PsiSc^{-\infty,-2l+1}(X,\calC)$ is
uniformly bounded in $\PsiSc^{-\infty,-2l}(X,\calC)$ as $r\to 0$.
Now, $G_r\in\PsiSc^{-\infty,2}(X,\calC)$ is uniformly bounded
in $\PsiSc^{0,0}(X,\calC)$, hence as a bounded operator on $L^2_\scl(X)$.
Thus, if $M>0$ is chosen sufficiently large, then $G_r+G_r^*\geq- M$ for
all $r\in(0,1)$, so
\begin{equation}
\psi(H)A^*(G_r+G_r^*+M)A\psi(H)\geq 0.
\end{equation}
Adding this to \eqref{eq:prop-61} shows that
\begin{equation}\label{eq:prop-61a}
ix^{l+1/2}[A^*_rA_r,H]x^{l+1/2}
\geq x^{l+1/2}((2-\epsilon')B_r^*B_r+E_r+F_r)x^{l+1/2}.
\end{equation}
The point of the commutator calculation is that in
$L^2_\scl(X)$
\begin{equation}\begin{split}
\langle u,&[A_r^*A_r,H]u\rangle\\
&=\langle u,A_r^*A_r(H-\lambda)u\rangle
-\langle u,(H-\lambda)A_r^*A_r u\rangle\\
&=2i\im\langle u,A_r^*A_r(H-\lambda)u\rangle;
\end{split}\end{equation}
the pairing makes sense for $r>0$ since $A_r\in\PsiSc^{-\infty,-l}(X,\calC)$.
Now apply \eqref{eq:prop-61} to $x^{-l-1/2}u$ and pair it with
$x^{-l-1/2}u$ in $L^2_{\scl}(X)$. Then for $r>0$
\begin{equation}\label{eq:pos-comm-est}
\|B_r u\|^2\leq|\langle u,E_r u\rangle|+|\langle u,F_ru\rangle|
+2|\langle u,A^*_rA_r(H-\lambda)u\rangle|.
\end{equation}

Letting $r\to 0$ now keeps the right hand side of \eqref{eq:pos-comm-est}
bounded.
In fact,
$A_r(H-\lambda)u\in\dCinf(X)$ remains bounded in $\dCinf(X)$ as $r\to 0$.
Similarly, by \eqref{eq:comm-51a},
$E_r u$ remains bounded in $\dCinf(X)$
as $r\to 0$ if we chose $K'$ so large that $\WFSc(u)\subset K'$.
Also, $F_r$ is bounded in $\bop(\Hsc^{m,l}(X),\Hsc^{-m,-l}(X))$, so
$\langle u,F_r u\rangle$ stays bounded by \eqref{eq:comm-51a} as well.
These estimates
show that $B_r u$ is uniformly bounded in $L^2_\scl(X)$.
Since $(1+r/x)^{-1}\to \Id$ strongly
on $\bop(\Hsc^{m',l'}(X),\Hsc^{m',l'}(X))$, we conclude that
$Bx^{-l-1/2}u\in L^2_\scl(X)$. By \eqref{eq:comm-51a}
and Proposition~\ref{prop:WF-3} this implies that for every $m$,
\begin{equation}
\xi_0\nin\WFSc^{m,l+1/2}(u).
\end{equation}
This is exactly the iterative step we wanted to prove. In the next step
we decrease $\delta>0$ slightly to ensure that $\WFScp(F)\subset \supp\qt$
is disjoint from $\WFSc^{m,l+1/2}(u)$.
\end{proof}

To state and prove the finer estimate under the assumption that
all elements of $\calC$ are totally geodesic,
first note that in geodesic normal
coordinates around $p\in C'$, $h-(|\mu|^2+|\nu|^2)$ vanishes together
with its first derivatives at $p=(0,0)$. Thus, $\scHg$ agrees with
$W^\flat$ on $\sct_p X$ where
\begin{equation}
W^\flat=2\tau(\mu\cdot\partial_\mu+\nu\cdot\partial_\nu)-2(|\mu|^2+|\nu|^2)
\partial_\tau+2\nu\cdot\partial_z+2\mu\cdot\partial_y.
\end{equation}
We will use $W^\flat$
to model the bicharacteristic flow of $H$. Note that $W^\flat$ is
the (rescaled) Hamilton vector field of the metric function
$\tau^2+|\nu|^2+|\mu|^2$, i.e.\ where we replace the actual metric
$h$ by a flat one.

We remark that it is the $\partial_\mu$ and
$\partial_\nu$ components
of $\scHg$ that differ from $W^\flat$ on $\sct_p X$
if we do not assume that the elements
of $\calC$ are totally geodesic. The former is inconsequential since
we only consider $\pi$-invariant functions (in particular, the only
$\mu$-dependence is via $\eta=y\cdot\mu$), but the latter rules out
the more precise location of the singularities given in the following
proposition.

\begin{prop}\label{prop:normal-prop-fine}
Suppose that $H$ is a many-body Hamiltonian satisfying \eqref{eq:hypo-1}
and that every element of $\calC$ is totally geodesic with respect to $h$.
Let $u\in\dist(X)$,
$\lambda>0$. Given $K\subset\Sigma_n(\lambda)\cap\sct_{C'}(C;X)$ compact
with $K\cap\WFSc((H-\lambda)u)=\emptyset$
there exist constants $C_0>0$,
$\delta_0>0$ such that the following holds. If
$\xi_0=(0,\tau_0,\nu_0)\in K$ and
for some $0<\delta<\delta_0$, $C_0\delta^{1/2}\leq\epsilon<1$ and
for all $\alpha=(y,z,\tau,\mu,\nu)\in\sct_{\partial X} X\cap
\Sigma_{\Delta-\lambda}$
\begin{equation}\begin{split}\label{eq:prop-9}
\alpha\in\sct_{C'_b}X &\Mand
|\pi_{0a}^e(\exp(\delta W^\flat)(\alpha))-\xi_0|\leq\epsilon\delta\Mand
|y(\exp(\delta W^\flat)(\alpha))|\leq\epsilon\delta\\
&\Rightarrow
\pi_{0b}(\alpha)\nin\WFSc(u)
\end{split}\end{equation}
then $\xi_0\nin\WFSc(u)$.
\end{prop}

\begin{proof}
We again employ an iterative argument, so we assume
that $\xi_0\nin\WFSc^{*,l}(u)$ and we need
to show that $\xi_0\nin\WFSc^{*,l+1/2}(u)$.

We first construct a $\Cinf$ function $\omega$ of $z$, $\tau$, $\nu$,
$\eta=\mu\cdot y$ and
$s=|y|^2$ which measures the distance of bicharacteristics
of $\Delta$ in $\Sigma_{\Delta-\lambda}$ from $\pi_{a0}^{-1}(\xi_0)\cap\Sigma
_{\Delta-\lambda}$.
Thus, $\tau^2+|\nu|^2+|\mu|^2-\lambda$ will be
small along these bicharacteristics.
We will take $\omega$ of the form
\begin{equation}\label{eq:n-pr-8}
\omega=\omega_0^2+(|y|^2-\frac{(y\cdot\mu)^2}
{\lambda-\tau_0^2-|\nu_0|^2})^2
\end{equation}
where $\omega_0$ only depends on $z$, $\tau$, $\nu$ and $\eta=y\cdot\mu$.
Note that
\begin{equation}
|y|^2-\frac{(y\cdot\mu)^2}{|\mu|^2}=|y-\frac{y\cdot\mu}{|\mu|^2}\mu|^2
\end{equation}
is the squared distance of the integral curves of $H_{|\mu|^2}$, which
are just straight lines, from $y=0$, so near $\Sigma_{\Delta-\lambda}$
the second term in $\omega$ gives the fourth power of this distance.

Pushing forward $W^\flat$ by the map $F:(y,z,\tau,\mu,\nu)\mapsto (z,\tau,\nu,
\mu\cdot y)$ at some point $\alpha=(y,z,\tau,\mu,\nu)$,
we obtain the vector
\begin{equation}
F_*|_\alpha W^\flat=2(\tau\eta+|\mu|^2)\partial_\eta
+2\tau\nu\cdot\partial_\nu-2(|\mu|^2+|\nu|^2)
\partial_\tau+2\nu\cdot\partial_z.
\end{equation}
Since we are interested in what happens near $\Sigma_{\Delta-\lambda}
\cap\sct_p X$,
where $\lambda=\tau^2+|\nu|^2+|\mu|^2$,
we are led to consider the constant coefficient vector field
\begin{equation}\label{eq:def-W_0}
W_0=2(\lambda-\tau_0^2-|\nu_0|^2)
\partial_\eta
+2\tau_0\nu_0\cdot\partial_\nu-2(\lambda-\tau_0^2)
\partial_\tau+\nu_0\cdot\partial_z
\end{equation}
in the variables $(z,\tau,\nu,\eta)$,
so
\begin{equation}\label{eq:n-pr-17}
F_*|_\alpha W^\flat=W_0+2(\lambda-\tau_0^2-|\mu_0|^2-|\nu_0|^2))
(-\partial_\eta+\partial_\tau).
\end{equation}

Note that the $\partial_\eta$ component of $W_0$ is nonzero.
Let
\begin{equation}
z_0(t)=\frac{W_0 z}{W_0 \eta}t,
\ \tau_0(t)=\tau_0+\frac{W_0 \tau}{W_0 \eta}t,
\ \nu_0(t)=\nu_0+\frac{W_0 \nu}{W_0 \eta}t,
\end{equation}
so
\begin{equation}
\gamma:t\mapsto(z_0((W_0\eta) t),\tau_0((W_0\eta) t),
\nu_0((W_0\eta) t),(W_0 \eta) t)
\end{equation}
gives a curve through $(\xi_0,0)$ with tangent vector $W_0$.
Now we define $\omega_0$ by
\begin{equation}
\omega_0=(z-z_0(\eta))^2+(\tau-\tau_0(\eta))^2+
(\nu-\nu_0(\eta))^2
\end{equation}
so $\omega_0$ vanishes exactly quadratically along $\gamma$ and
is positive elsewhere, and
\begin{equation}
W_0\omega_0=0.
\end{equation}
Note that by the triangle inequality
\begin{equation}
|z|+|\tau-\tau_0|+|\nu-\nu_0|+|\eta|\leq C(\omega_0^{1/2}+|\eta|)
\end{equation}
for sufficiently large $C$.

Since for $\alpha\in\pih^{-1}(\xi_0)$ we have $F_*|_\alpha\scHg=
F_*|_\alpha W^\flat=W_0$,
we see that
\begin{equation}
\scHg(z-z_0(\eta))=0\ \text{at}\ \pih^{-1}(\xi_0),
\end{equation}
i.e.\ when $y=0$, $z=0$, $\tau=\tau_0$, $\nu=\nu_0$, $g=\lambda$, so
\begin{equation}
|\scHg(z-z_0(\eta))|\leq C(|y|+\omega_0^{1/2}+|\lambda-g|).
\end{equation}
Hence,
\begin{equation}
|\scHg(z-z_0(\eta))^2|\leq 2C\omega_0^{1/2}(|y|+\omega_0^{1/2}+|\lambda-g|).
\end{equation}
Similar conclusions hold for $\tau-\tau_0(\eta)$ and $\nu-\nu_0(\eta)$, so
\begin{equation}
|\scHg \omega_0|\leq C_1(|y|+\omega_0^{1/2}+|\lambda-g|)\omega_0^{1/2}.
\end{equation}

Next, we calculate $\scHg(|y|^2-(y\cdot\mu)^2
/(\lambda-\tau_0^2-|\nu_0|^2))$.
Since the function we are differentiating
vanishes quadratically at $y=0$, the same follows for its derivatives
with respect to any vector field tangent to $y=0$. Since the $\partial_y$
component of $\scHg$ (and of $W^\flat$) is of the form $2\mu\cdot\partial_y
+\sum \beta_j\partial_{y_j}$ with $\beta_j$ vanishing at $y=0$, $z=0$
(i.e. at $p$), we
conclude that
\begin{equation}
|(\scHg-2\mu\cdot \partial_y)
(|y|^2-(y\cdot\mu)^2/(\lambda-\tau_0^2-|\nu_0|^2))|
\leq C_2|y|(|y|+|z|).
\end{equation}
On the other hand,
\begin{equation}
(2\mu\cdot \partial_y)(|y|^2-(y\cdot\mu)^2/(\lambda-\tau_0^2-|\nu_0|^2))
=4(\mu\cdot y)\frac{\lambda-\tau_0^2-|\nu_0|^2-|\mu|^2}
{\lambda-\tau_0^2-|\nu_0|^2}.
\end{equation}
But, as in \eqref{eq:n-prop-30aa},
\begin{equation}\begin{split}\label{eq:n-prop-30}
|\lambda&-\tau_0^2-|\nu_0|^2-|\mu|^2|
\leq|\lambda-g|
+C'(|y|+|z|+|\tau-\tau_0|+|\nu-\nu_0|)\\
&\leq C_3(|\lambda-g|+|y|+\omega_0^{1/2}).
\end{split}\end{equation}
Thus,
\begin{equation}
|\scHg(|y|^2-(y\cdot\mu)^2/(\lambda-\tau_0^2-|\nu_0|^2))|
\leq C_4|y|(|\lambda-g|+|y|+\omega_0^{1/2}).
\end{equation}

Our results thus far imply that
\begin{equation}
|\scHg\omega|\leq C_{5}\omega^{1/2}(|y|+|\lambda-g|+\omega_0^{1/2})^2.
\end{equation}

Now let $1>\epsilon>0$, $\delta>0$, with other restrictions to be imposed on
these later, and let
\begin{equation}\label{eq:phi-def-normal}
\phi=\tau_0-\tau+\frac{1}{\epsilon^4\delta^3}\omega.
\end{equation}
We use $\tau_0-\tau$ to measure propagation along the bicharacteristics;
$\eta=y\cdot\mu$ would also work. We again
let $\chi_0\in\Cinf(\Real)$ be equal to $0$ on $(-\infty,0]$ and
$\chi_0(t)=\exp(-1/t)$ for $t>0$ and we
let $\chi_1\in\Cinf(\Real)$ be $0$
on $(-\infty,0]$, $1$ on $[1,\infty)$, with $\chi_1'\geq 0$ satisfying
$\chi_1'\in\Cinf_c((0,1))$.
Furthermore, for $A_0>0$ large, to be determined, $t\in(0,1)$, let
\begin{equation}\label{eq:prop-22-fine}
q_t=q=\chi_0(A_0^{-1}(1+t-\phi/\delta))\chi_1((\tau_0-\tau+\delta)/
(\epsilon\delta)+t).
\end{equation}
We usually simply write $q$ in place of $q_t$. We only use $t$ to slightly
shrink the support of $q$ in our inductive proof (i.e.\ as $l$ is increasing),
instead of adjusting $\delta$ as in the proof of
Proposition~\ref{prop:normal-prop}.
Thus, on $\supp q$ we have $\phi\leq 2\delta$ and $\tau_0-\tau\geq-2\delta$.
Since $\omega\geq 0$, the first of these inequalities implies that
$\tau_0-\tau\leq 2\delta$, so
\begin{equation}
|\tau-\tau_0|\leq 2\delta \Mand
\omega\leq \epsilon^4\delta^3(2\delta+\tau-\tau_0)\leq 4\epsilon^4\delta^4.
\end{equation}
Hence, $\omega_0\leq 2\epsilon^2\delta^2$, which together with $|\tau-\tau_0|
\leq2\delta$ gives $|\eta|=|\mu\cdot y|\leq C_{6}\delta$ since the
$\partial_\tau$ component of $W$ in non-zero. Since we also have
\begin{equation}
||y|^2-(y\cdot\mu)^2/(\lambda-\tau_0^2-|\nu_0|^2)|
\leq 2\epsilon^2\delta^2,
\end{equation}
we conclude that $|y|\leq C_{7}\delta$. Thus, under the additional
assumption
\begin{equation}\label{eq:prop-27}
|\lambda-g|<\delta
\end{equation}
we deduce that
$|\scHg\omega|\leq C_{8}\epsilon^2\delta^4$, so
\begin{equation}
|\scHg\phi-2h|\leq C_{8}\delta/\epsilon^2.
\end{equation}
Hence, for $c_0>0$, $C_0>0$ appropriately chosen and for $\epsilon\in(0,1)$,
$\delta>0$ satisfying $\delta/\epsilon^2<C_0$, we have
\begin{equation}
\scHg\phi>c_0>0.
\end{equation}

Again, this directly gives a positivity estimate for $\scHg q$ near $\xi_0$.
Now
\begin{equation}\begin{split}\label{eq:scHg q-calc-fine}
\scHg q=-A_0^{-1}&\delta^{-1}\chi'_0(A_0^{-1}(1+t-\phi/\delta))
\chi_1((\tau_0-\tau+\delta)/
(\epsilon\delta)+t)\scHg \phi\\
&-(\epsilon\delta)^{-1}
\chi_0(A_0^{-1}(1+t-\phi/\delta))\chi_1'((\tau_0-\tau+\delta)/
(\epsilon\delta)+t)\scHg \tau.
\end{split}\end{equation}
Hence, with
\begin{equation}\begin{split}\label{eq:scHg q-pos-fine}
b^2&=c_0 A_0^{-1}\delta^{-1}\chi'_0(A_0^{-1}(1+t-\phi/\delta))
\chi_1((\tau_0-\tau+\delta)/(\epsilon\delta)+t),\\
e&=-(\epsilon\delta)^{-1}
\chi_0(A_0^{-1}(1+t-\phi/\delta))\chi_1'((\tau_0-\tau+\delta)/
(\epsilon\delta)+t)\scHg \tau
\end{split}\end{equation}
we have
\begin{equation}
\scHg q\leq -b^2+e.
\end{equation}

In addition, similarly to \eqref{eq:prop-33}-\eqref{eq:prop-34}, we see that
\begin{equation}\label{eq:prop-33-fine}
b^2\geq (c_0 A_0/16) q.
\end{equation}
Moreover, with $\partial$ denoting a partial derivative with respect to
one of $(y,z,\tau,\mu,\nu)$,
\begin{equation}\begin{split}
\partial q=-A_0^{-1}&\delta^{-1}\chi'_0(A_0^{-1}(1+t-\phi/\delta))
\chi_1((\tau_0-\tau+\delta)/
(\epsilon\delta)+t)\partial \phi\\
&-(\epsilon\delta)^{-1}
\chi_0(A_0^{-1}(1+t-\phi/\delta))\chi_1'((\tau_0-\tau+\delta)/
(\epsilon\delta)+t)\partial \tau.
\end{split}\end{equation}
Thus, \eqref{eq:prop-38}-\eqref{eq:prop-39} hold, and hence
the uniform bounds of \eqref{eq:comm-37-p} also follow.
Now $e$ is supported where
\begin{equation}\label{eq:supp-e-fine}
-\delta-t\epsilon\delta\leq\tau_0-\tau\leq-\delta+(1-t)\ep\delta,
\ \omega^{1/4}\leq\sqrt 2\,\ep\delta,
\end{equation}
so near the backward direction along bicharacteristics through $\xi_0$,
in the region which we know is disjoint from $\WFSc(u)$. In addition,
by choosing $\delta>0$ sufficiently small, we can assume that the support
of $q$, $e$ and $b$ are all disjoint from $\WFSc((H-\lambda)u)$.

From this point we can simply follow the proof of
Proposition~\ref{prop:normal-prop}.
Thus, we conclude that for every $m$,
\begin{equation}
\xi_0\nin\WFSc^{m,l+1/2}(u).
\end{equation}
This is exactly the iterative step we wanted to prove. In the next step
we decrease $t$ slightly to ensure that $\supp\qt_t$
is disjoint from $\WFSc^{m,l+1/2}(u)$.
\end{proof}

Before proving the general tangential propagation estimate, we first do
it in the totally geodesic case ($\calC$ totally geodesic).
Proposition~\ref{prop:tgt-bichar-tot-geod}
shows that for sufficiently short times there is a unique generalized
broken bicharacteristic through
any point in $\Sigma_t(\lambda)$, namely the integral curve of $\scHg$.
The simplicity of this description may already give a hint that it is
particularly easy to prove the corresponding propagation estimate for
singularities.
Indeed, in the proof of the aforementioned proposition,
we have essentially already
constructed the pseudo-differential operator $A$ to commute through $H$
by defining the $\pi$-invariant function $\phi$ (which will play an
analogous role to that of $\phi$ in the proof of normal propagation).
The following argument may also clarify the close relationship between 
proving results about the geometry
of the generalized broken bicharacteristics and the positive commutator
proof of propagation estimates. Again, we only state it for forward
propagation.

\begin{prop}\label{prop:tot-geod-tgt-prop}
Suppose that $H$ is a many-body Hamiltonian satisfying \eqref{eq:hypo-1}.
Suppose also that every element of $\calC$ is totally geodesic with
respect to $h$.
Let $u\in\dist(X)$,
$\lambda>0$. Let $\xi_0\in\Sigma_t(\lambda)\cap\sct_{C'}(C;X)$, $C=C_a$,
satisfy $\xi_0\nin\WFSc((H-\lambda)u)$.
Then there exists $\ep'>0$ such that if in addition for some $s\in(-\ep',0)$
we have
\begin{equation*}
\pi_{0a}(\exp(s\scHg)(\hat\pi^{-1}(\xi_0)))\nin\WFSc(u)
\end{equation*}
then
$\xi_0\nin\WFSc(u)$.
\end{prop}

\begin{proof}
First note that there is nothing to prove if $\xi_0\in R_+(\lambda)
\cup R_-(\lambda)$, so from now on we assume that $\xi_0\nin R_+(\lambda)
\cup R_-(\lambda)$.
The proof is very similar to the previous one and the positive commutator
construction is exactly the same as in three-body scattering
\cite[Proposition~15.4]{Vasy:Propagation-2}, based on the $\pi$-invariant
function $\phi$ used here in the proof of
Proposition~\ref{prop:tgt-bichar-tot-geod}.
Thus, we take local coordinates centered at $C$ as above, i.e.\ of the form
$(y,z)$, and let $\phi=\phi^{(\ep)}$ be defined by \eqref{eq:def-phi-tot-geod},
so in particular $\phi$ is $\pi$-invariant.
In the proof of Proposition~\ref{prop:tgt-bichar-tot-geod} we showed that
there exists $\delta_0\in (0,1)$ such that for any $\delta\in(0,\delta_0)$
and any $\ep\in(0,1)$
\begin{equation}
\phi(\xit)\leq 2\delta,\ \tau(\xit)-\tau_0\leq 2\delta\Mand
|\tau^2(\xit)+h(\xit)-\lambda|<\epsilon\delta
\end{equation}
imply that $\scHg\phi$ satisfies \eqref{eq:scHg-phi-est-tot-geod}, so
\begin{equation}
\scHg\phi(\xit)\geq c_0>0.
\end{equation}
We define $q$ as in \eqref{eq:prop-22-fine}. Then \eqref{eq:scHg q-calc-fine},
hence
\eqref{eq:scHg q-pos-fine}-\eqref{eq:supp-e-fine} also hold.
Since $\ep>0$ can be taken
arbitrarily small, we can choose it and $\delta\in(0,\delta_0)$ so that
$\supp e$ is a small neighborhood of $\exp(s\scHg)(\hat\pi^{-1}(\xi_0))$;
in particular, $\pi_{0b}(\supp e)$ is disjoint from $\WFSc(u)$ for each $b$.
We can then apply the compactness argument of
Proposition~\ref{prop:normal-prop} to prove \eqref{eq:prop-61} for the
operators $A$, $B$, etc., defined in that proof, and conclude that
$\xi_0\nin\WFSc(u)$.
\end{proof}

We now return to the general setting of not necessarily totally geodesic
$\calC$.

\begin{prop}\label{prop:tgt-prop}
Suppose that $H$ is a many-body Hamiltonian satisfying \eqref{eq:hypo-1}.
Let $u\in\dist(X)$,
$\lambda>0$. Given
\begin{equation}
K\subset(\Sigma_t(\lambda)\cap\sct_{C'}(C;X))\setminus(R_+(\lambda)\cup
R_-(\lambda)\cup\WFSc((H-\lambda)u)
\end{equation}
compact there exist constants $C_0>0$,
$\delta_0>0$ such that the following holds. If
$\xi_0=(z_0,\tau_0,\nu_0)\in K$ and
for some $0<\delta<\delta_0$, $C_0\delta\leq\epsilon<1$ and
for all $\alpha=(y,z,\tau,\mu,\nu)\in\sct_{\partial X} X\cap
\Sigma_{\Delta-\lambda}$
\begin{equation}\begin{split}
\alpha\in\sct_{C'_b}X &\Mand
|\pi_{0a}^e(\exp(\delta W^\flat)(\alpha))-\xi_0|\leq\epsilon\delta\Mand
|y(\exp(\delta W^\flat)(\alpha))|\leq\epsilon\delta\\
&\Rightarrow
\pi_{0b}(\alpha)\nin\WFSc(u)
\end{split}\end{equation}
then $\xi_0\nin\WFSc(u)$.
\end{prop}

\begin{proof}
The proof is very similar to the previous ones and now
the positive commutator
construction follows that of \cite[Proposition~15.2]{Vasy:Propagation-2}
in three-body scattering.
Thus, we take local coordinates as above, i.e.\ of the form
$(y,z)$ with $C_b$ defined by linear equations in $y$.
Then we
construct $\omega_0\in\Cinf(\sct_{C'}(C;X))$ (defined near $\xi_0$) to
measure
the squared distance from integral curves of
\begin{equation}
W^\sharp=2\tau\nu\cdot\partial_\nu-2\htil\partial_\tau+H_{\htil};
\end{equation}
this is achieved by solving a Cauchy problem as in
\cite{Vasy:Propagation-2} and in \eqref{eq:tot-geod-Cauchy} here. (Indeed,
an approximate construction, like that
of $\omega_0$ in the normal case discussed above, would also work). Then we
extend $\omega_0$ to a function on $\sct_{\bX} X$ (using the coordinates
$(y,z,\tau,\mu,\nu)$ near $\bX$),
let
\begin{equation}
\omega=\omega_0+|y|^2,\ \phi=\tau_0-\tau+\frac{1}{\ep^2\delta}\omega,
\end{equation}
and define $q$ as in \eqref{eq:prop-22}.
The difference in the powers of $\ep$ and $\delta$ in this definition
of $\phi$ in the (general) tangential setting and 
that in the normal case (given in \eqref{eq:phi-def-normal})
arises since in the normal
setting $\omega$ approximates the fourth power of the distance from
the generalized bicharacteristics while here it approximates the squared
distance.
The estimates on $\scHg\phi$
are just as in \cite[Proposition~15.2]{Vasy:Propagation-2},
see also the proof
of Proposition~\ref{prop:tgt-bichar-tot-geod} here in the similar
totally geodesic setting (the estimates are simply better but no different
in nature under the totally geodesic assumption since now we do not have
\eqref{eq:tot-geod-metric-2b}), giving a slightly
better result than in the totally geodesic normal case: it is $\delta/\ep$, not
$\delta/\ep^2$, that has to be
bounded below by an appropriate positive constant. The difference arises
as the model integral curves in the tangential setting are closer to the
actual ones than in the normal setting. Thus, one obtains
\eqref{eq:prop-33} here as well.
The functional analysis part, under the assumption
that there are no bound states, is exactly as in the normal case.
\end{proof}

An argument of Melrose-Sj\"ostrand \cite{Melrose-Sjostrand:I,
Melrose-Sjostrand:II}, see also \cite[Chapter~XXIV]{Hor} and
\cite{Lebeau:Propagation} allows us to conclude our main result
concerning the singularities of generalized eigenfunctions of $H$.
Here we concentrate on totally geodesic $\calC$ (since that is the
case in Euclidean scattering), in which case the more delicate
tangential propagation argument of Melrose-Sj\"ostrand is not
necessary. The proof presented below essentially follows Lebeau's
paper \cite[Proposition~VII.1]{Lebeau:Propagation}.
We thus have the following theorem.

\begin{thm}\label{thm:prop-sing}
Let $(X,\calC)$ be a locally locally linearizable many-body space, and
ssuppose that $H$ is a many-body Hamiltonian satisfying \eqref{eq:hypo-1}.
Let $u\in\dist(X)$,
$\lambda>0$. Then $\WFSc(u)\setminus\WFSc((H-\lambda)u)$
is a union of maximally extended
generalized broken bicharacteristics of $\Delta-\lambda$.
\end{thm}

\begin{proof}
We only need to prove that for every $a$, if
$\xi_0\in\WFSc(u)\setminus\WFSc((H-\lambda)u)$ and $\xi_0\in\sct_{C'_a}
(C_a;X)$
then there exists a generalized broken bicharacteristic
$\gamma:[-\ep_0,\ep_0]\to\dot\Sigma$, $\ep_0>0$,
with $\gamma(0)=\xi_0$ and such that
$\gamma(t)\in\WFSc(u)\setminus\WFSc((H-\lambda)u)$ for $t\in[-\ep_0,\ep_0]$.
In fact, if this statement holds
for all $a$ with $C_c\subset C_a$, let
\begin{equation}\begin{split}
\calR=&\{\text{generalized broken bicharacteristics}\\
&\ \gamma:(\alpha,\beta)
\to (\WFSc(u)\setminus\WFSc((H-\lambda)u))\cap\bigcup
\{\sct_{C'_a}(C_a;X):\ C_c\subset C_a\},\\
&\ \alpha<0<\beta,\ \gamma(0)=\xi_0\},
\end{split}\end{equation}
and put the natural partial order on $\calR$, so $\gamma\leq\gamma'$
if the domains satisfy $(\alpha,\beta)\subset(\alpha',\beta')$ and
$\gamma=\gamma'|_{(\alpha,\beta)}$. Then $\calR$ is not empty
and every non-empty
totally ordered subset
of $\calR$ has an upper bound,
so an application of Zorn's lemma gives
a maximal generalized broken bicharacteristic of $\Delta-\lambda$
in the intersection of
$\WFSc(u)\setminus\WFSc((H-\lambda)u)$ with $\cup_{C_c\subset C_a}
\sct_{C'_a}(C_a;X)$ which passes through $\xi_0$. A similar maximal
statement holds if we replace $C_c\subset C_a$ by $C_c\subsetneq C_a$.

Indeed, it suffices to show that
for any $a$, if
\begin{equation}\label{eq:prop-103}
\xi_0\in\WFSc(u)\setminus\WFSc((H-\lambda)u)\Mand\xi_0\in\sct_{C'_a}
(C_a;X)
\end{equation}
then
\begin{equation}\label{eq:prop-104}\begin{split}
&\text{there exists a generalized broken bicharacteristic}
\ \gamma:[-\ep_0,0]\to\dot\Sigma,\ \ep_0>0,\\
&\qquad\qquad \gamma(0)=\xi_0,
\ \gamma(t)\in\WFSc(u)\setminus\WFSc((H-\lambda)u),\ t\in[-\ep_0,0],
\end{split}\end{equation}
for the existence of a generalized broken bicharacteristic
on $[0,\ep_0]$ can be demonstrated similarly
by replacing the forward propagation
estimates by backward ones, and, directly from
Definition~\ref{Def:gen-br-bichar}, piecing together the two
generalized broken bicharacteristics gives one defined on $[-\ep_0,\ep_0]$.

Note that if every element of $\calC$ is totally geodesic and $\xi_0
\in\Sigma_t(\lambda)$ then,
due to
Proposition~\ref{prop:tot-geod-tgt-prop},
\eqref{eq:prop-103}$\Rightarrow$\eqref{eq:prop-104}.
Indeed, to prove this statement
for $\xi_0\in \sct_{C'_a}(C_a;X)$, we only need that $C_b$ be totally
geodesic for $C_b$ with $C_a\subset C_b$ (since the result is local);
in particular, it always holds at $C'_0$.

We proceed to prove that \eqref{eq:prop-103} implies
\eqref{eq:prop-104} by induction on $a$. As remarked above, this is
certainly true for $a=0$. We only prove the implication here under
the assumption that all elements of $\calC$ are totally geodesic,
the general case repeats the argument of Melrose-Sj\"ostrand.
Thus, we already know the implication if $\xi_0\in\Sigma_t(\lambda)
\cap\sct_{C'_a}(C_a;X)$.

So suppose that \eqref{eq:prop-103}$\Rightarrow$\eqref{eq:prop-104}
has been proved for all $b$ with
$C_a\subsetneq C_b$ and that $\xi_0\in\Sigma_n(\lambda)\cap\sct_{C'_a}
(C_a;X)$ satisfies \eqref{eq:prop-103}.
As noted in the first paragraph, we thus know
that
the intersection of
$\WFSc(u)\setminus\WFSc((H-\lambda)u)$ with $\cup_{C_a\subsetneq C_b}
\sct_{C'_b}(C_b;X)$ is a union of maximally extended
generalized broken bicharacteristics of $\Delta-\lambda$.
We use the notation of the proof of
Proposition~\ref{prop:normal-prop} below.
Let $U$ be a neighborhood of $\xi_0=(0,z_0,\tau_0,\nu_0)$ in
$\dot\Sigma$ which is given
by equations of the form $|y|<\delta'$, $|z-z_0|<\delta'$,
$|\tau-\tau_0|<\delta'$, $|\nu-\nu_0|<\delta'$, $\delta'>0$,
such that $\scHg\eta>0$
on $\pih^{-1}(U)$ and $U\cap\WFSc((H-\lambda)u)=\emptyset$.
Such a neighborhood exists since $\xi_0\nin\WFSc((H-\lambda)u)$
and $\scHg\eta(\xit_0)
=\lambda-\tau_0^2-\htil(z_0,\nu_0)>0$ for every $\xit_0\in\pih^{-1}(\xi_0)$.
Also let $U'$ be a subset of $U$ defined by replacing $\delta'$ by a
smaller $\delta''>0$.
By Proposition~\ref{prop:normal-prop}, there is a sequence of points
$\xi_n\in\dot\Sigma$ such that $\xi_n\in\WFSc(u)$,
$\xi_n\to\xi_0$ as $n\to\infty$, and
$\eta(\xi_n)<0$ for all $n$, so we may assume that
$\xi_n\in U'$ for all $n$.
By the inductive hypothesis, there exist generalized broken bicharcteristics
$\gamma_n:[-\ep_0,0]\to\dot\Sigma$ with $\gamma_n(0)=\xi_n$ and such
that if for all $t\in[-\ep',0]$, $0<\ep'\leq\ep_0$, we have
$\gamma_n(t)\nin\sct_{C'_a}(C_a;X)$, then
$\gamma_n(t)\in\WFSc(u)$ for $t\in[-\ep',0]$. But $\eta$ is increasing on
generalized broken bicharacteristics in $U$
since $\scHg\eta>0$ there, so we conclude that
$y(\gamma_n(t))\cdot\mu(\gamma_n(t))=\eta(\gamma_n(t))\leq\eta(\gamma_n(0))<0$
for $t\in[-\ep_0,0]$, hence $y(\gamma_n(t))\neq 0$, so $\gamma_n(t)\in\WFSc(u)$
for all $t\in[-\ep_0,0]$. By Proposition~\ref{prop:Lebeau-compactness},
applied with $K=\WFSc(u)$, there
is a subsequence of $\gamma_n$ converging uniformly to a generalized broken
bicharacteristic $\gamma:[-\ep_0,0]\to\WFSc(u)$. In particular,
$\gamma(0)=\xi_0$ and $\gamma(t)\in\WFSc(u)$ for all $t\in[-\ep_0,0]$,
providing the inductive step.

Note that this argument for normal points
$\xi_0$ does not use that the elements of $\calC$ are totally geodesic; it
works equally well in the general case. Thus, for non-totally geodesic
$\calC$ now we only need to consider $\xi_0\in\Sigma_t(\lambda)$, and,
as mentioned above, this can be done by the Melrose-Sj\"ostrand
argument.
\end{proof}

We remark that the result is optimal as can be seen by considering the
Euclidean setting, taking potentials
singular at a specified $C_a$, thereby placing ourselves into the three-body
framework. As \cite{Vasy:Structure} shows, singularities do reflect in all
permissible directions in general, the reflection being governed to top
order by the (two-body) S-matrix of the subsystem.

\section{The resolvent}\label{sec:resolvent}
Before we can turn Theorem~\ref{thm:prop-sing} into a result on the
wave front relation of the S-matrix, we need to analyze the resolvent.
More precisely, we need to understand the boundary values
\begin{equation}
R(\lambda\pm i0)=(H-(\lambda\pm i0))^{-1}
\end{equation}
of the resolvent at the real axis in a microlocal sense.
To do so, we also need estimates at the radial sets $R_\pm(\lambda)$.
Since the Hamilton vector field of the metric $g$ vanishes at
$R_+(\lambda)\cup R_-(\lambda)$, the estimates must utilize
the weights $x^{-l-1}$ themselves. In this sense they are delicate,
but on the other hand they only involve $x$ and its sc-microlocal
dual variable $\tau$, so they do not need to reflect the
geometry of $\calC$. The best known positive commutator estimate is
the Mourre estimate, originally proved by Perry, Sigal and Simon
in Euclidean many-body scattering \cite{Perry-Sigal-Simon:Spectral},
in which one takes $q=x^{-1}\tau$ with the notation
of Section~\ref{sec:commutators}. Since it is easy to analyze the commutator
of powers of $x$ with $H$ (in particular, they commute with $V$),
the functional calculus allows one to obtain microlocal estimates
from these, as was done by G\'erard, Isozaki and Skibsted
\cite{GerComm, GIS:N-body}.
Thus, nearly all the technical results in this section can be proved, for
example, by using the Mourre estimate and Theorem~\ref{thm:prop-sing}.
In particular, apart from the propagation statements, they are well-known
in Euclidean many-body scattering. The generalization of these Euclidean
results to our
geometric setting is straightforward;
the arguments essentially follow those in three-body scattering
that were used in \cite{Vasy:Propagation-2}.

We first state the weak form of the limiting absorption principle, namely
that for $f\in\dCinf(X)$,
$R(\lambda\pm it)f$, $t>0$, has a limit in $\Hsc^{m,l}(X)$,
$m$ arbitrary, $l<-1/2$, as $t\to 0$. To simplify
the asymptotic expansions of $R(\lambda\pm i0)f$ which we also
describe, for $\lambda>0$ we introduce
the functions
\begin{equation}
\alpha_\pm=\alpha_{\pm,\lambda}=\pm\frac{V}{2\sqrt{\lambda}\,x}
\in\Cinf(X\setminus C_{0,\sing}),
\end{equation}
and the set of polyhomogeneous functions
$\bconc^{\indfam}(X\setminus C_{0,\sing})$ on $X\setminus C_{0,\sing}$ with
index set
\begin{equation}
\indfam=\{(m,p):\ m,p\in\Nat,\ p\leq 2m\}.
\end{equation}
Recall from \cite{RBMCalcCon}
that $v\in\bconc^{\indfam}(X\setminus C_{0,\sing})$
means that $v$ is $\Cinf$ in the interior of $X$ and it has a full
asymptotic expansion at $C'_0$ which in local coordinates $(x,y)$
take the form
\begin{equation}
v(x,y)\sim\sum_{j=0}^\infty\sum_{r\leq 2j}
x^j(\log x)^r a_{j,r}(y),\quad a_{j,r}\in\Cinf(C'_0).
\end{equation}
Thus,
$v\in \calC^0(X\setminus C_{0,\sing})$ and $|v(x,y)-a_{0,0}(y)|\leq
Cx|\log x|^2$.

\begin{thm}\label{thm:lim-absorb}
Suppose that $H$ is a many-body Hamiltonian satisfying \eqref{eq:hypo-1},
$\lambda>0$.
Let $f\in\dCinf(X)$, $u_t^\pm=R(\lambda\pm it)f$, $t>0$. Then
$u_t^\pm$ has a limit $u_\pm=R(\lambda\pm i0)f$ in $\Hsc^{m,l}(X)$, $l<-1/2$,
as $t\to 0$. In addition,
\begin{equation}\label{eq:lim-absorb}
\WFSc(u_\pm)\subset R_\mp(\lambda).
\end{equation}
If $V$ is short-range, i.e.\ $V\in x^2\Cinf(X\setminus C_{0,\sing})$, then
\begin{equation}
u_\pm=e^{\pm i\sqrt\lambda/x}x^{(n-1)/2}v_\pm,\quad v_\pm\in\Cinf(X\setminus
C_{0,\sing}),
\end{equation}
while if $V$ is long-range, i.e.\ $V$ merely satisfies \eqref{eq:mb-pot},
then
\begin{equation}
u_\pm=e^{\pm i\sqrt\lambda/x}x^{(n-1)/2+i\alpha_\pm}v_\pm,\quad
v_\pm\in\bconc^{\indfam}(X\setminus
C_{0,\sing}).
\end{equation}
\end{thm}

\begin{rem}
The first statement in the theorem also holds if we merely assume
$f\in\Hsc^{m,l'}(X)$ with $l'>1/2$, but then $\WFSc(u_\pm)$ has to be replaced
by the filtered wave front set $\WFSc^{m,l'-1}(u_\pm)$. Moreover,
$R(\lambda\pm i0)$ give continuous operators from $\Hsc^{m,l'}(X)$ to
$\Hsc^{m+2,l}(X)$.
\end{rem}

\begin{proof}
This result is a weak form of the limiting absorption principle and
can be proved by a Mourre-type estimate. In the Euclidean setting,
it is a combination of the Mourre estimate, proved by
Perry, Sigal and Simon \cite{Perry-Sigal-Simon:Spectral}, and its
microlocalized version obtained by G\'erard, Isozaki and Skibsted
\cite{GerComm}. In the geometric setting, the Mourre estimate
describes the commutator of $H$ with a self-adjoint first order differential
operator $A\in x^{-1}\Diffsc^1(X)$
such that $A-xD_x\in\Diffsc^1(X)$ (this is of course a restriction only at
$\bX$). Namely, it says that for $\phi\in\Cinf_c(\Real;[0,1])$ supported
sufficiently close to $\lambda$, we have
\begin{equation}\label{eq:Mourre-1}
i\phi(H)[A,H]\phi(H)\geq 2(\lambda-\ep)\phi(H)^2+R,\quad\ep>0,
\end{equation}
where $R\in\PsiSc^{-\infty,1}(X,\calC)$, hence compact on $L^2_\scl(X)$.
It was proved in the geometric three-body setting (with an appropriate
adjustment to allow bound states of subsystems) in
\cite{Vasy:Propagation-2}, following the Euclidean argument of
Froese and Herbst \cite{FroMourre}. The proof 
given there goes through essentially unchanged
for more than three bodies.
Under our assumption \eqref{eq:hypo-1}, the symbolic commutator calculation
in the scattering calculus, $\scHg(x^{-1}\tau)+2g\in x\Cinf(\sct X)$,
and a slight modification of
Corollary~\ref{cor:comm-7}, prove the
Mourre estimate. The argument of \cite{Perry-Sigal-Simon:Spectral} then
proves the existence of the limits $u_\pm$ in $\Hsc^{0,l}(X)$, $l<-1/2$,
and $(H-\lambda)u_\pm=f\in\dCinf(X)$ shows that the same holds
in $\Hsc^{m,l}(X)$ for every $m$ and for every $l<-1/2$.

To show the flavor of the arguments, we prove here a version
of the estimate of G\'erard, Isozaki and Skibsted
\cite{GerComm}. Such arguments as this can be combined to prove
the limiting absorption principle without a direct use of the Mourre
estimate as was done in the geometric two-body type setting by Melrose
\cite{RBMSpec} and in the geometric three-body setting in
\cite{Vasy:Propagation-2}. Here, however, we concentrate on proving
the wave front set result. The major difference between the propagation
estimates of the previous section and the ones near $R_{\pm}(\lambda)$
is that $\scHg$ is radial at $R_+(\lambda)\cup R_-(\lambda)$: it has
the form $2\tau x \partial_x$. Thus, we need to use a weight $x^{-l-1}$
to obtain a positive symbol estimate. So for $l>-1$, let
\begin{equation}
q=x^{-l-1}\chi(\tau)\psit(x)\geq 0
\end{equation}
where $\psit\in\Cinf_c(\Real)$ is identically $1$ near $0$ and
is supported in a bigger neighborhood of $0$ (it is simply a cutoff
near $\bX$), $\chi\in\Cinf_c(\Real;[0,1])$ vanishes on
$(-\infty,\sqrt{\lambda}-2\ep)$, identically $1$ on $(\sqrt{\lambda}-\ep,
\infty)$, $\ep>0$, $\chi'\geq 0$,
and $\chi$ vanishes with all derivatives at every $t$
with $\chi(t)=0$. Then for sufficiently small $\delta>0$,
$|g-\lambda|=|\tau^2+h-\lambda|<\delta$ implies
\begin{equation}\begin{split}
&\scHg q=-2((l+1)\tau\chi(\tau)+h\chi'(\tau))x^{-l-1}
\leq -b^2x^{-l-1},\\
&\quad b=(2(l+1)\tau\chi(\tau)+(\lambda-\tau^2)\chi'(\tau)/2)^{1/2}.
\end{split}\end{equation}
Thus, both $x^{l+1}q$ and $x^{l+1}b$ are $\pi$-invariant. Let $A\in
\PsiSc^{-\infty,-l-1}(X,\calC)$ be
a quantization of $q$ as in Lemma~\ref{lemma:comm-1}, except that now $q$ is
not supported in a single coordinate chart, so we need to define $A$ as
the sum of localized operators (of course, this is not necessary in the
actual Euclidean setting). Thus, roughly speaking,
$A$ is the product of a quantization of $q$ and $\psi_0(H)$, $\psi_0\in
\Cinf_c(\Real)$. The fact that $q\in x^{-l-1}\Cinf(\sct X)$ does not cause any
trouble, and the argument of Corollary~\ref{cor:comm-7} shows that
for $\psi\in\Cinf_c(\Real;[0,1])$ supported sufficiently close to
$\lambda$ we have
\begin{equation}\label{eq:GIS-7}
ix^{l+1/2}\psi(H)[A^*A,H]\psi(H)x^{l+1/2}\geq x^{l+1/2}((2-\ep')B^*B+F)
x^{l+1/2},\quad\ep'>0,
\end{equation}
where
\begin{equation}
F\in\PsiSc^{-\infty,-2l}(X,\calC),\ \WFScp(F)\subset\supp (x^{l+1}q),
\end{equation}
\begin{equation}\label{eq:GIS-9}
B\in\PsiSc^{-\infty,-l-1/2}
(X,\calC),\ \Bh_{a,-l-1/2}(\zeta)=b(\zeta)q(\zeta)^{1/2}\psi(\Hh_{a,0}
(\zeta)).
\end{equation}
Let
\begin{equation}
A_0=A\psi(H)\in\PsiSc^{-\infty,-l-1}(X,\calC).
\end{equation}
Since
\begin{equation}
\langle u^+_t,i[A_0^*A_0,H]u^+_t\rangle
=-2\im\langle u^+_t,A^*_0A_0(H-(\lambda+it))u^+_t
\rangle-2t\|A_0 u^+_t\|^2,
\end{equation}
we conclude that
\begin{equation}\label{eq:GIS-17}
\|B u^+_t\|^2+2t\|A_0 u^+_t\|^2\leq
|\langle u^+_t,F u^+_t\rangle|
+2|\langle u^+_t,A_0^*A_0(H-(\lambda+it))u^+_t\rangle|.
\end{equation}
Since $t>0$, the second term on the left hand side can be dropped.
Since $u^+_t\to u_+$ in $\Hsc^{0,l'}(X)$ for $l'<-1/2$, we conclude
that for $l\in(-1,-1/2)$ the right hand side stays bounded as $t\to 0$.
Thus, $B u^+_t$ is uniformly bounded in $L^2_\scl(X)$, and as $u^+_t\to
u_+$ in $\Hsc^{0,l'}(X)$, we conclude that $Bu_+\in L^2_\scl(X)$.
But then \eqref{eq:GIS-9} shows that for any $\zeta$ with $q(\zeta)\neq 0$,
we have $\zeta\nin\WFSc^{m,l+1/2}(u_+)$ for every $m$.
This proves that the incoming radial set,
$R_+(\lambda)$, is disjoint from $\WFSc^{m,l+1/2}(u_+)$, $l+1/2\in(-1/2,0)$.
Iterating the argument, as in the proof of Proposition~\ref{prop:normal-prop},
gives that $\WFSc(u_+)\cap R_+(\lambda)=\emptyset$. Since $\WFSc(u_+)$
is closed, the same conclusion holds for a neighborhood of $R_+(\lambda)$.
Finally, as all generalized broken bicharacteristics of $\Delta-\lambda$
tend to $R_+(\lambda)$ as $t\to-\infty$ and
$(H-\lambda)u_+=f\in\dCinf(X)$, the propagation of singularities
theorem, Theorem~\ref{thm:prop-sing}, implies that $\WFSc(u_+)\subset
R_-(\lambda)$. The existence of the asymptotic expansions is a local
question, so at $C'_0$ we can work in the scattering calculus
to prove it, see \cite{Vasy:Asymptotic} for details of the proof.
\end{proof}

A pairing argument immediately shows $R(\lambda\pm i0)v$ also exists
for distributions $v\in\dist(X)$ with wave front set disjoint from the
incoming and outgoing radial set respectively. Combining it with the
propagation theorem, Theorem~\ref{thm:prop-sing}, we
can deduce the following result; as usual, we assume that $(X,\calC)$ is
locally linearizable.

\begin{thm}\label{thm:res-WF}
Suppose that $H$ is a many-body Hamiltonian satisfying \eqref{eq:hypo-1},
$\lambda>0$.
Suppose also that $v\in\dist(X)$ and $\WFSc(v)\cap R_+(\lambda)=\emptyset$.
Let $u_t^+=R(\lambda+ it) v$, $t>0$. Then $u^+_t$
has a limit
$u_+=R(\lambda+ i0)v$ in $\dist(X)$ as $t\to 0$ and
$\WFSc(u_+)\cap R_+(\lambda)=\emptyset$. Moreover, if
$\xi\in\dot\Sigma\setminus R_-(\lambda)$ and
every maximally backward extended generalized broken
bicharacteristic, $\gamma:(-\infty,t_0]\to\dot\Sigma$,
with $\gamma(t_0)=\xi$ is disjoint from $\WFSc(v)$, then $\xi\nin\WFSc(u_+)$.
The result also holds with $R_+(\lambda)$ and $R_-(\lambda)$ interchanged,
$R(\lambda+it)$ replaced by $R(\lambda-it)$, $(-\infty,t_0]$
by $[t_0,\infty)$ and correspondingly `backward extended' by
`forward extended'.
\end{thm}

\begin{proof}
As mentioned above, the first part follows from the self-adjointness of $H$,
so that for $t>0$, $v\in\dist(X)$, $f\in\dCinf(X)$, we have
$v(R(\lambda+it)f)=R(\lambda+it)v(f)$; recall that the distributional
pairing is the real pairing, not the complex (i.e.\ $L^2$) one. The
wave front statement of Theorem~\ref{thm:lim-absorb} and the assumption on
$v$ show the existence of the limit $u_+=R(\lambda+i0)v$ in $\dist(X)$
and that in addition $\WFSc^{m,l}(u_+)\cap R_+(\lambda)=\emptyset$ for
every
$l<-1/2$. The positive commutator argument of Theorem~\ref{thm:lim-absorb}
then applies and shows that $\WFSc(u_+)\cap R_+(\lambda)=\emptyset$.
In the Euclidean setting these results follow from a microlocalized
version of the Mourre estimate due to G\'erard, Isozaki and Skibsted
\cite{GIS:N-body}; see \cite{Hassell-Vasy:Symbolic} for a detailed
argument.

Finally,
since $\WFSc(u_+)$ is closed, a neighborhood of $R_+(\lambda)$
in $\dot\Sigma$ is disjoint from $\WFSc(u_+)$. Since all generalized
broken bicharacteristics approach $R_+(\lambda)$ as $t\to -\infty$ by
Proposition~\ref{prop:param-bichar},
the last part follows from $(H-\lambda)u_+=v$ and
Theorem~\ref{thm:prop-sing}. It can be also proved by modifying the argument
of Propositions~\ref{prop:normal-prop}-\ref{prop:tgt-prop} along the
lines of our proof of Theorem~\ref{thm:lim-absorb}. Namely,
we consider the family $u^+_t\in\dist(X)$, $t>0$, and note that
for $t>0$, $R(\lambda+it)\in\PsiSc^{-2,0}(X,\calC)$, so $\WFSc(u^+_t)
\subset\WFSc(v)$. Let $A_0$, etc., be defined as $A_r$ with $r=0$ where
$A_r$ is given by \eqref{eq:prop-57} (i.e.\ we do not need to use
the approximating factor $(1+r/x)^{-1}$). Then
\begin{equation}
\langle u^+_t,i[A^*_0A_0,H]u^+_t\rangle
=-2\im\langle u^+_t,A_0^*A_0(H-(\lambda+it))u^+_t
\rangle-2t\|A_0u^+_t\|^2.
\end{equation}
Note that the pairings make sense since now $\WFScp(A_0)$ is disjoint from
$\WFSc(u^+_t)$, $t>0$.
Thus,
\begin{equation}\begin{split}\label{eq:res-17}
\|B_0 u^+_t\|^2+2t\|A_0 u^+_t\|^2\leq
|\langle u^+_t,E_0 u^+_t\rangle|&+|\langle u^+_t,F_0u^+_t\rangle|\\
&+2|\langle u^+_t,A^*_0A_0(H-(\lambda+it))u^+_t\rangle|.
\end{split}\end{equation}
Since $t>0$, the second term can be dropped from the left hand side.
Thus, knowing that $u^+_t\to u_+$ in $\dist(X)$ as $t\to 0$, and
assuming that $\xi_0\nin\WFSc^{m,l}(u_+)$ has already been proved
and \eqref{eq:prop-9} is satisfied by $u_+$, we conclude that
$\xi_0\nin\WFSc^{m,l+1/2}(u_+)$. The
iteration of this argument of Proposition~\ref{prop:normal-prop} and the
similar arguments for tangential propagation allow us to conclude
the forward propagation estimates which can then be turned into maximal
statements as we did in Theorem~\ref{thm:prop-sing}. This argument
also shows the influence of the sign of $t$: if $t<0$, the inequality
\eqref{eq:res-17} cannot be used to derive results on $u_+$. Instead,
the signs are then correct in the backward estimate, just as expected.
\end{proof}

We conclude this section with the following uniqueness theorem on
solutions of $(H-\lambda)u=0$. It is essentially a geometric version
of Isozaki's uniqueness theorem
\cite[Theorem~1.3]{IsoUniq}, though we allow arbitrary growth of $u$ away
from one of the radial sets, say $R_+(\lambda)$.

\begin{thm}\label{thm:uniqueness}
Suppose that $H$ is a many-body Hamiltonian satisfying \eqref{eq:hypo-1},
$\lambda>0$.
Suppose also that $u\in\dist(X)$, $(H-\lambda)u=0$ and $\WFSc^{m,l}(u)
\cap R_+(\lambda)=\emptyset$ for some $m$ and some $l>-1/2$. Then
$u=0$. The same conclusion holds if we replace $R_+(\lambda)$ by
$R_-(\lambda)$.
\end{thm}

\begin{proof}
Just as in the proof of Theorem~\ref{thm:res-WF},
the positive commutator estimate of
Theorem~\ref{thm:lim-absorb} (but now applied with a regularizing factor
in $x$)
shows that
$\WFSc(u)\cap R_+(\lambda)=\emptyset$, and then Theorem~\ref{thm:prop-sing}
shows that
\begin{equation}
\WFSc(u)\subset R_-(\lambda).
\end{equation}
We remark that although we need a regularizing factor here which
requires some changes in the proof, e.g.\ see the argument of the
paragraph below, the regularizing factor (whether
$(1+r/x)^{-1}$ or another one) commutes with $V$, so the additional
arguments for dealing with it are essentially the same as
the two-body ones. Thus,
the regularization part of the
proof of $\WFSc(u)\cap R_+(\lambda)=\emptyset$ essentially follows
\cite[Proposition~10]{RBMSpec}.

We proceed to show that
\begin{equation}\label{eq:unique-5}
m\in\Real,\ l<-1/2\Rightarrow\WFSc^{m,l}(u)\cap R_-(\lambda)=\emptyset.
\end{equation}
We give the details below since regularity arguments for distributions
which are large at infinity seem to appear less often in the literature than
the `finer ones';
in particular, \cite[Theorem~1.3]{IsoUniq} assumes
$u\in\Hsc^{m,l}(X)$ with $l>-1$.
We essentially
follow the proof of \cite[Proposition~9]{RBMSpec} below.

So suppose that \eqref{eq:unique-5} has been shown for some $l<-1$; we
now show it with $l$ replaced by $l+1/2$.
This time we consider
\begin{equation}\label{eq:unique-7}
q=x^{-l-1}\chi(\tau)\psit(x),\ l<-1,
\end{equation}
where $\psit\in\Cinf_c(\Real)$ is identically $1$ near $0$ and
is supported in a bigger neighborhood of $0$ (it is simply a cutoff
near $\bX$), $\chi\in\Cinf_c(\Real;[0,1])$ identically $1$ on
$(-\infty,-\sqrt{\lambda}+\ep)$, vanishes on $(-\sqrt{\lambda}+2\ep,
\infty)$, $\ep>0$, and $\chi$ vanishes with all derivatives at every $t$
with $\chi(t)=0$.
Then
\begin{equation}
\scHg q=-2((l+1)\tau\chi(\tau)+h\chi'(\tau))x^{-l-1}=(-b^2+e)x^{-l-1},
\end{equation}
\begin{equation}
b^2=2(l+1)\tau\chi(\tau).
\end{equation}
The first
key point now is that on $\WFSc(u)$ we have $\tau=-\sqrt{\lambda}$,
so $\WFSc(u)\cap\pi(\supp e)=\emptyset$.
Let $A\in\PsiSc^{-\infty,-l-1}(X,\calC)$ as in Lemma~\ref{lemma:comm-1}.
Corollary~\ref{cor:comm-7} again shows that
for $\psi\in\Cinf_c(\Real;[0,1])$ supported sufficiently close to
$\lambda$ we have
\begin{equation}\label{eq:unique-21}
ix^{l+1/2}\psi(H)[A^*A,H]\psi(H)x^{l+1/2}\geq x^{l+1/2}((2-\ep')B^*B+E+F)
x^{l+1/2},\quad\ep'>0,
\end{equation}
where
\begin{equation}\begin{split}
&B\in\PsiSc^{-\infty,-l-1/2}
(X,\calC),\ \Bh_{a,-l-1/2}(\zeta)=b(\zeta)q(\zeta)^{1/2}\psi(\Hh_{a,0}
(\zeta)),\\
&E\in\PsiSc^{-\infty,-2l-1}(X,\calC),\ \WFScp(E)\cap\WFSc(u)=\emptyset,\\
&F\in\PsiSc^{-\infty,-2l}(X,\calC),\ \WFScp(F)\subset\supp (x^{l+1}q).
\end{split}\end{equation}

Let
\begin{equation}
A_r=A(1+r/x)^{-1}\psi(H),\ B_r=B(1+r/x)^{-1},\ E_r=(1+r/x)^{-1}E(1+r/x)^{-1},
\end{equation}
so
\begin{equation}
A_r\in\PsiSc^{-\infty,-l}(X,\calC)\ \text{for}\ r>0,\ A_r
\ \text{is uniformly bounded in}\ \PsiSc^{-\infty,-l-1}(X,\calC);
\end{equation}
analogous statements also hold for $B_r$ and $E_r$. Thus,
\begin{equation}\begin{split}\label{eq:unique-25}
ix^{l+1/2}&[A_r^*A_r,H]x^{l+1/2}\\
&=i(1+r/x)^{-1}x^{l+1/2}\psi(H)[A^*A,H]\psi(H)x^{l+1/2}(1+r/x)^{-1}\\
&\quad+i\psi(H)A^*x^{l+1}(G_r+G_r^*)x^{l+1}A\psi(H)+H_r
\end{split}\end{equation}
where $H_r$ is uniformly bounded in $\PsiSc^{-\infty,1}(X,\calC)$ and
\begin{equation}
G_r=i\psi_0(H)^2 x^{-1}
(1+r/x)^{-1}[(1+r/x)^{-1},H],
\end{equation}
$\psi_0\in\Cinf_c(\Real;[0,1])$, $\psi_0\equiv 1$ on $\supp\psi$,
so $G_r$ is uniformly bounded in $\PsiSc^{-\infty,1}(X,\calC)$.
Thus, we need to estimate the commutator $[(1+r/x)^{-1},H]$, and now we
do not have a large $M$ as in the proof of
Proposition~\ref{prop:normal-prop} to help us deal with it.

The other key point is thus that we have
$i[(1+r/x)^{-1},H]=i[(1+r/x)^{-1},\Delta]$ and
\begin{equation}
\scHg(1+r/x)^{-1}=2\tau\frac{r}{x+r}
=-c_r^2+d_r,\ c_r=\chi_1(\tau)
(-\tau)^{1/2}\left(\frac{r}{x+r}\right)^{1/2},
\end{equation}
$\chi_1\in\Cinf_c(\Real;[0,1])$ identically $1$ on
$(-\infty,-\sqrt{\lambda}+3\ep)$, vanishes on $(-\sqrt{\lambda}+4\ep,
\infty)$, $\ep>0$. Let $C_r$ be the quantization of $c_r$ multiplied
by $\psi_0(H)$ as in Lemma~\ref{lemma:comm-1}, and define $D_r$
similarly but with $\psi_0(H)$ replaced by $\psi_0(H)^2$. Thus,
as $(1+r/x)^{-1}$ is uniformly bounded in the symbol class $S^0(X)$,
\begin{equation}\label{eq:unique-31}
i\psi_0(H)x^{-1/2}[(1+r/x)^{-1},H]x^{-1/2}\psi_0(H)=C_r^*C_r+D_r+H'_r
\end{equation}
with $C_r$ and $D_r$ uniformly bounded in $\PsiSc^{-\infty,0}(X,\calC)$,
$C_r \in \PsiSc^{-\infty,1/2}(X,\calC)$ for $r>0$,
$D_r \in \PsiSc^{-\infty,1}(X,\calC)$ for $r>0$,
and $H'_r$ uniformly
bounded in $\Psisc^{-\infty,1}(X,\calC)$. Moreover,
$D_r A\in\PsiSc^{-\infty,\infty}(X,\calC)$ uniformly due to the disjoint
operator wave front sets. Thus,
\begin{equation}
G_r+G_r^*= 2(1+r/x)^{-1/2}(C_r^*C_r+D_r)(1+r/x)^{-1/2}
+H''_r
\end{equation}
with $H''_r$ uniformly bounded in $\PsiSc^{-\infty,1}(X,\calC)$, so
\begin{equation}\begin{split}\label{eq:unique-35}
\psi(H)&A^*x^{l+1}(G_r+G_r^*)x^{l+1}A\psi(H)\\
&= 2\psi(H)
A^*x^{l+1}(1+r/x)^{-1/2}(C_r^*C_r+D_r)(1+r/x)^{-1/2}x^{l+1}A\psi(H)
+H^\flat_r\\
&\geq H^\sharp_r,
\end{split}\end{equation}
$H^\flat_r$, $H^\sharp_r$ uniformly bounded in $\PsiSc^{-\infty,1}(X,\calC)$.
Combining \eqref{eq:unique-21}, \eqref{eq:unique-25} and \eqref{eq:unique-35},
we see that
for $\ep'>0$ we have
\begin{equation}
ix^{l+1/2}[A^*_rA_r,H]x^{l+1/2}
\geq x^{l+1/2}((2-\ep')B_r^*B_r+E_r+F_r)
x^{l+1/2}.
\end{equation}
We deduce as at the end of the proof
of Proposition~\ref{prop:normal-prop} that
$\WFSc^{m,l+1/2}(u)\cap R_-(\lambda)=\emptyset$ for every
$m$ and for every $l+1/2<-1/2$, so \eqref{eq:unique-5} holds. In
particular, $u\in\Hsc^{m,l}(X)$ for every $m$ and for every $l<-1/2$.

In the Euclidean setting we can now simply refer to Isozaki's
uniqueness theorem \cite[Theorem~1.3]{IsoUniq} to conclude that $u=0$.
Here we give some details to indicate how this conclusion can be reached
in general.
The crucial step is improving the estimate past the critical
regularity $\Hsc^{*,-1/2}(X)$. In the Euclidean setting this was done
by Isozaki \cite[Lemma~4.5]{IsoRad} and his argument was adapted to the
geometric setting in \cite[Proposition~17.8]{Vasy:Propagation-2}. We
thus
conclude that $\WFSc^{m,l}(u)\cap R_-(\lambda)=\emptyset$ for
$l\in(0,-1/2)$. This is the point where $(H-\lambda)u=0$, and not just
$(H-\lambda)u\in\dCinf(X)$ is used.
Now we can apply a commutator estimate like that
of Theorem~\ref{thm:lim-absorb} but near $R_-(\lambda)$.
Thus, we conclude that $\WFSc(u)\cap
R_-(\lambda)=\emptyset$, so $u\in\dCinf(X)$. The theorem
of Froese and Herbst \cite{FroExp} on the absence of bound states
with positive energy adapted to the geometric setting, as discussed in
\cite[Appendix~B]{Vasy:Propagation-2}, concludes that
$u=0$.
\end{proof}

\section{The Poisson operator and the scattering matrix}\label{sec:S-matrix}
Just as in \cite{Vasy:Propagation, Vasy:Propagation-2} where three-body
scattering was analyzed, the propagation
of singularities for generalized eigenfunctions of $H$ implies the
corresponding result for the (free-to-free)
scattering matrix, $S(\lambda)$, of $H$. Note that this is the only S-matrix
under our assumption of the absence of bound states of the subsystems.
We start by discussing the Poisson operator, then we use it to relate
the propagation of singularities for generalized eigenfunctions to
the wave front relation of the S-matrix.

The result that allows us to define the Poisson operator is that
if $V$ is short-range, i.e.\ $V\in x^2\Cinf(X\setminus C_{0,\sing})$, then
for $\lambda\in(0,\infty)$
and $g\in\Cinf_c(C_0')$,
there is a unique $u\in\dist(X)$
such that $(H-\lambda)u=0$, and $u$ has the form
\begin{equation}\label{eq:intro-21-P}
u=e^{-i\sqrt{\lambda} x}x^{(n-1)/2}v_-
+R(\lambda+i0)f,
\end{equation}
where $v_-\in\Cinf(X)$, $v_-|_{\bX}=g$, and
$f\in\dCinf(X)$. For long-range $V$ the same statement is valid with
the asymptotic expansion replaced by one similar to that of
Theorem~\ref{thm:lim-absorb}:
\begin{equation}\label{eq:intro-31-P}
u=e^{-i\sqrt{\lambda} x}x^{(n-1)/2+i\alpha_-}v_-+R(\lambda+i0)f,\quad
v_-\in\bconc^{\indfam}(X).
\end{equation}
The Poisson operator with initial state in the free-cluster
is then the map
\begin{equation}
P_{+}(\lambda):\Cinf_c(C'_0)\to\dist(\Snp),\quad P_{+}(\lambda)g=u.
\end{equation}
(Note that the subscript $0$ for the free cluster has been
dropped here in contrast to the introduction and \cite{Vasy:Scattering}.)
To see that such a $u$ is unique, note that the difference $v=u-u'$ of two
distributions $u$ and $u'$ with the above properties satisfies
$(H-\lambda)v=0$ and $\WFSc^{0,0}(v)\cap R_+(\lambda)=\emptyset$ by
Theorem~\ref{thm:lim-absorb}, so $v=0$ due to Theorem~\ref{thm:uniqueness}.
To see the existence of such $u$,
note that as $\supp g\subset C'_0$, we can construct
\begin{equation}
u_-=e^{-i\sqrt{\lambda} x}x^{(n-1)/2}v_-,\ v_-\in\Cinf(X),\ v_-|_{\bX}
=g,\ -f=(H-\lambda)u_-\in\dCinf(X),
\end{equation}
by a local calculation as in \cite{RBMSpec}, i.e.\ essentially
in a two-body type setting. (We need to make slight changes in the
asymptotic expansion for long-range $V$ as described above.) Thus, we
construct the Taylor series of $v_-$
at $\bX$ explicitly, so we can even arrange that $\supp v_-\cap C_{0,\sing}
=\emptyset$. Then $u=u_-+R(\lambda+i0)f$ is of the form \eqref{eq:intro-21-P}
and satisfies $(H-\lambda)u=0$ indeed.

We need to understand the Poisson operator better before
we can extend it to distributions.
So first recall from \cite[Section~19]{Vasy:Propagation-2}
that the Melrose-Zworski \cite{RBMZw} construction of a parametrix for
the Poisson operator in the two-body type setting ($\calC$ is empty) gives
`the initial part' of a parametrix $\Pt_+(\lambda)$ for the Poisson operator
with free initial state
in the many-body setting (for three bodies in that paper, but this makes no
difference).
Although the construction is performed there for short range potentials,
it can be easily adjusted to long range potentials decaying like $x$, see
\cite[Appendix~A]{Vasy:Propagation-2} and \cite[Section~3]{Vasy:Geometric}.
In particular, the kernel of
$\Pt_\pm(\lambda)$ is of the form
\begin{equation}
\Kf_\pm(x,y,y')=e^{\mp i\sqrt{\lambda}\cos\distance(y,y')/x}
x^{i\alpha_\mp(y')}a_\pm(x,y,y')|dh|,
\end{equation}
where $\distance$ is the distance function of the boundary metric $h$,
$|dh|$ is the Riemannian density associated with it, $\alpha_\pm$ are
given by \eqref{eq:intro-31}, and
$a_\pm\in\Cinf(X\times C'_0)$ are cut off to be
supported near $y=y'$. Here $y'$ is the `initial point' of the
plane waves, so $y'\in C'_0$ corresponds to considering free
incoming particles. In
Euclidean scattering $\Kf_\pm$
takes the form $e^{\mp i w\cdot y'}a_\pm(w,y')|dh|$,
$w=y/x$ is the Euclidean variable and $|dh|$ the standard measure on the
sphere; and e.g.\ if the potentials $V_b$
are Schwartz
then $a_\pm$ are just cutoff functions supported near $y=y'$ which are
constant in a smaller neighborhood of $y=y'$. In general, $a_\pm(0,y,y)$
is determined by the condition that
\begin{equation}\label{eq:Pt-asymp}
\Pt_\pm(\lambda)g=e^{\mp i\sqrt{\lambda}/x}x^{i
\alpha_\mp+(n-1)/2}v_\pm,
\end{equation}
$v_\pm\in\bconc^{\indfam}(X)$,
$v_\pm|_{\bX}=g$, and then $a_\pm(0,y,y')$, as well as the other terms
of the Taylor series of $a_\pm$ at $x=0$ can be calculated from
transport equations near $y=y'$. Finally, we cut off the solutions
to the transport equations close to $y=y'$ before reaching $C_{0,\sing}$.
Thus, for $y'$ in a fixed compact subset in $C'_0$, $K(x,y,y')$ is
supported away from $C_{0,\sing}$, so for $g\in\dist_c(C'_0)$,
$\supp(\Pt_\pm(\lambda)g)$ is disjoint from $C_{0,\sing}$.

The most important properties of $\Pt_\pm(\lambda)$ are summarized in the
following proposition. Although we state them for $\Pt_+(\lambda)$ only,
they also hold for $\Pt_-(\lambda)$ with the appropriate sign changes.
Here we use $\sim'_\pm$ as the relation on $S^*\bX\times\Sigma
_{\Delta-\lambda}$ defined analogously to
$\sim_\pm$ (see Definition~\ref{Def:bichar-rel}), but with `generalized
broken bicharacteristics' replaced by `bicharacteristics of $\Delta-\lambda$'.
Note that generalized broken bicharacteristics are simply bicharacteristics
in $\sct_{C'_0}X$ which is where we will apply to following result.

\begin{prop}\label{prop:Pt}(\cite[Proposition~A.1]{Vasy:Propagation-2})
$\Kf_+\in\dist(X\times C'_0;\Omega_R)$, constructed above, is the kernel of an
operator
$\Pt_+(\lambda):\Cinf_c(C'_0)\to\dist(X)$, which extends to an operator
$\Pt_+(\lambda):\dist_c(C'_0)\to\dist(X)$, and for $g\in\dist_c(C'_0)$
\begin{equation}
\supp(\Pt_+(\lambda)g)\cap C_{0,\sing}=\emptyset,
\end{equation}
\begin{equation}\begin{split}\label{eq:Pt-WF}
\WFsc(\Pt_+(\lambda)g)\subset&\{(y,\sqrt{\lambda},0):\ y\in\supp g\}\\
&\quad\cup\{\xi\in\dot\Sigma\setminus (R_+(\lambda)\cup R_-(\lambda))
:\ \exists\zeta\in\WF(g),\ \xi\sim'_-\zeta\},
\end{split}\end{equation}
\begin{equation}\begin{split}\label{eq:Pt-error}
\WFsc((H&-\lambda)\Pt_+(\lambda)g)\\
&\quad\subset\{\xi\in\dot\Sigma\setminus (R_+(\lambda)\cup R_-(\lambda)):
\ \exists\zeta\in\WF(g),\ \xi\sim'_-\zeta\}.
\end{split}\end{equation}
\end{prop}

The actual Poisson operator is then given by
\begin{equation}\label{eq:P_+-def}
P_+(\lambda)=\Pt_+(\lambda)-R(\lambda+i0)(H-\lambda)\Pt_+(\lambda),
\end{equation}
with a similar definition of $P_-(\lambda)$:
\begin{equation}\label{eq:P_--def}
P_-(\lambda)=\Pt_-(\lambda)-R(\lambda-i0)(H-\lambda)\Pt_-(\lambda),
\end{equation}
Indeed, if $g\in\Cinf_c(C'_0)$ then
$(H-\lambda)\Pt_+(\lambda)g
\in\dCinf(X)$ and $\Pt_+(\lambda)g$ has an asymptotic expansion
as in \eqref{eq:Pt-asymp}, so by Theorem~\ref{thm:lim-absorb},
$(H-\lambda)P_+(\lambda)g=0$ and
$P(\lambda)g$ has the form \eqref{eq:intro-21-P}
(with changes as indicated in \eqref{eq:intro-31-P} if
$V$ is long-range).
In addition, for $g\in\dist_c(C'_0)$,
$\WFSc((H-\lambda)\Pt_\pm(\lambda)g)$ is disjoint from
$R_\pm(\lambda)$ by Proposition~\ref{prop:Pt}.
Hence, by Theorem~\ref{thm:res-WF}, \eqref{eq:P_+-def}-\eqref{eq:P_--def}
indeed make sense.
We also immediately deduce from Theorem~\ref{thm:res-WF}

\begin{prop}\label{prop:P-WF}
Suppose that $H$ is a many-body Hamiltonian satisfying \eqref{eq:hypo-1}.
Then the Poisson operator $P_+(\lambda):\Cinf_c(C'_0)\to\dist(X)$
extends by continuity to an operator
$\Pt_+(\lambda):\dist_c(C'_0)\to\dist(X)$.
Moreover, for $g\in\dist_c(C'_0)$ we have
\begin{equation}\begin{split}
\WFSc(P_+(\lambda)g)\subset&\{(y,\sqrt{\lambda},0):\ y\in\supp g\}
\cup R_-(\lambda)\\
&\quad\cup\{\xi\in\dot\Sigma(\lambda)\setminus R_+(\lambda)
:\ \exists\zeta\in\WF(g),\ \xi\sim_-\zeta\}.
\end{split}\end{equation}
\end{prop}

Our definition of the free-to-free
S-matrix is based on asymptotic expansions of
generalized eigenfunctions. So let $g\in\Cinf_c(C'_0)$ and let
$u=P_+(\lambda)g$. By \eqref{eq:intro-21-P} (modified as in
\eqref{eq:intro-31-P} for long-range $V$) and Theorem~\ref{thm:lim-absorb},
$u$ has the form
\begin{equation}
u=e^{-i\sqrt{\lambda} x}x^{(n-1)/2}v_-+e^{i\sqrt{\lambda} x}x^{(n-1)/2}v_+
\end{equation}
with $v_-\in\Cinf(X)$, $v_+\in\Cinf(X\setminus C_{0,\sing})$, $v_-|_{\bX}=g$.
We then define the free-to-free S-matrix by
\begin{equation}
S(\lambda):\Cinf_c(C'_0)\to\Cinf(C'_0),\ S(\lambda)g=v_+|_{C'_0}.
\end{equation}

We need a better description of the S-matrix to describe its structure.
This can be done via a boundary pairing formula analogous to
\cite[Proposition~13]{RBMSpec}.
It gives the following
alternative description of the S-matrix, see
\cite[Equation~(5.7)]{Vasy:Scattering} (or its analogue from
\cite{Vasy:Propagation} in the non-Euclidean setting):

\begin{prop}\label{prop:S-matrix}
For $\lambda>0$ the
scattering matrix is given by
\begin{equation}\label{eq:S-m-2}
S(\lambda)=\frac{1}{2i\sqrt\lambda}((H-\lambda)\Pt_-(\lambda))^*P_+(\lambda).
\end{equation}
\end{prop}

\begin{proof}
The following pairing formula was proved by Melrose
\cite[Proposition~13]{RBMSpec} for short-range $V$, but the same
proof also applies when $V$ is long-range. Also, the proof can be easily
localized, see \cite[Proposition~3.3]{Vasy:Scattering}.
Suppose that for $j=1,2$, $u_j\in\dist(X)$,
\begin{equation}\begin{split}
&u_j=e^{i\sqrt\lambda/x}x^{(n-1)/2+i\alpha_+}v_{j,+}
+e^{-i\sqrt\lambda/x}x^{(n-1)/2+i\alpha_-}v_{j,-},\\
&\quad v_{j,\pm}\in\bconc^{\indfam}(X\setminus C_{0,\sing}),
\ \supp(v_{2,\pm})\Subset X\setminus C_{0,\sing},
\end{split}\end{equation}
and $f_j=(H-\lambda)u_j\in\dCinf(X)$. Let $a_{j,\pm}=v_{j,\pm}|_{\bX}$.
Then
\begin{equation}\label{eq:bdary-pair-1}
2i\sqrt\lambda\int_{\bX}(a_{1,+}\,\overline{a_{2,+}}-a_{1,-}\,\overline{a_{2,-}})
\,dh=\int_X(u_1\,\overline{f_2}-f_1\,\overline{u_2})\,dg.
\end{equation}
We apply this result with $u_1=P(\lambda)g$, $u_2=\Pt(-\lambda)f$.
By the construction of $\Pt(-\lambda)$ we conclude that $a_{2,+}=f$,
$a_{2,-}=0$, while for $u_1$ we see directly from the definition of
$S(\lambda)$ and $P(\lambda)$ that $a_{1,-}=g$, $a_{1,+}=S(\lambda)
g$. Substitution into \eqref{eq:bdary-pair-1} proves the proposition.
\end{proof}

Propositions~\ref{prop:Pt} and \ref{prop:P-WF}, when combined
with \eqref{eq:S-m-2},
allow us to deduce the structure of the S-matrix.

\begin{thm}
Let $(X,\calC)$ be a locally linearizable many-body space.
Suppose that $H$ is a many-body Hamiltonian satisfying \eqref{eq:hypo-1}.
Then the scattering matrix, $S(\lambda)$, extends to a continuous linear map
$\dist_c(C'_0)\to\dist(C'_0)$. The wave front relation of $S(\lambda)$ is
given by the generalized broken geodesic flow at time $\pi$.
\end{thm}

\begin{proof}
Let $f,g\in\dist_c(C'_0)$.
Suppose also that there is no generalized broken geodesic of length $\pi$
starting at some $\zeta\in\WF(g)$ and ending at $\zeta'\in\WF(f)$. That
means that for any $\xi\in\dot\Sigma\setminus(R_+(\lambda)
\cup R_-(\lambda))$ we cannot have $\xi\sim_-\zeta$, $\zeta\in\WF(g)$,
and $\xi\sim_+\zeta'$, $\zeta'\in\WF(f)$, at the same time.
Proposition~\ref{prop:Pt} (with $-$ signs instead of $+$)
implies that
\begin{equation}
\WFSc((H-\lambda)\Pt_-(\lambda)f)\subset
\dot\Sigma\setminus(R_+(\lambda)
\cup R_-(\lambda));
\end{equation}
indeed, we also have $\WFSc((H-\lambda)\Pt_-(\lambda)f)\subset\sct_{C'_0}X$,
so we can even replace $\WFSc$ by $\WFsc$.
Thus, by our assumption on $\WF(f)$ and $\WF(g)$, and by
Propositions~\ref{prop:Pt}-\ref{prop:P-WF}, we have
\begin{equation}
\WFsc((H-\lambda)\Pt_-(\lambda)f)\cap\WFsc(P_+(\lambda)g)=\emptyset.
\end{equation}
But the complex pairing
\begin{equation}
\langle u,u'\rangle_X=\int u\, \overline{u'}\,dg
\end{equation}
extends by continuity from $u,u'\in\dCinf(X)$ to $u,u'\in\dist(X)$
satisfying $\WFsc(u)\cap\WFsc(u')=\emptyset$.
To see this just let $A\in\Psisc^{0,0}(X)$ with
$\WFscp(A)\cap\WFsc(u)=\emptyset$, $\WFscp(\Id-A^*)\cap\WFsc(u')=\emptyset$,
and note that
\begin{equation}
\langle u,u'\rangle_X=\langle Au,u'\rangle_X+\langle u,(\Id-A^*)u'\rangle_X
\end{equation}
extends as claimed.
Hence, the pairing
\begin{equation}
\langle P_+(\lambda)g,(H-\lambda)\Pt_-(\lambda)f\rangle_X=
\langle((H-\lambda)\Pt_-(\lambda))^*P_+(\lambda)g,f\rangle_X
\end{equation}
defined first for $f,g\in\Cinf_c(C'_0)$ extends by continuity to
$f,g\in\dist_c(C'_0)$ satisfying our wave front condition. In other
words, $g$ can be paired with every distribution whose wave front set
has no elements related to $\WF(g)$ by the generalized broken geodesic flow
at time $\pi$. Thus, for any $A\in\Psop_c^{0}(C'_0)$ with $\WFp(A)$ disjoint
from the image of $\WF(g)$ under the generalized broken geodesic flow
at time $\pi$, and for any $f\in\dist_c(C'_0)$, $\langle AS(\lambda)g,f
\rangle_{\bX}=\langle S(\lambda)g,A^*f\rangle_{\bX}$ is defined by continuity
from $f\in\Cinf_c(C'_0)$, so $AS(\lambda)g\in\Cinf(C'_0)$. But this
states exactly that $\WF(S(\lambda)g)$ is contained in the image
of $\WF(g)$ under the generalized broken geodesic flow at time $\pi$.
\end{proof}

\appendix

\section{The proof of Proposition~\ref{prop:Lebeau-2-sided}}

In this appendix we prove Proposition~\ref{prop:Lebeau-2-sided}
under the assumption that $\calC$ is totally geodesic,
roughly following
Lebeau's original proof in \cite{Lebeau:Propagation}. As noted
after the statement of the proposition we can proceed inductively,
using the order on $\calC$. So assume that $\gamma(t_0)=\xi_0\in
\Sigma_n(\lambda)\cap\sct_{C'_a}(C_a;X)$. The inductive hypothesis is that
we have already proved the proposition for $b$ with $C_a\subset C_b$.
Thus,
by Definition~\ref{Def:gen-br-bichar}, part (ii), there exists $\delta'>0$
such that the conclusion of the proposition
holds if we replace $t_0$
replaced by $t\neq t_0$, assuming $|t-t_0|<\delta'$. Let $\xit_\pm(t)\in
\Sigma_{\Delta-\lambda}$,
$t\neq t_0$, be the points given by the inductive hypothesis.
We often write
\begin{equation}
\xit_\pm(t)=(y(t),z(t),\tau(t),\mu_\pm(y),\nu(t))
\end{equation}
in local coordinates, so e.g. $\tau(\xit_\pm(t))=\tau(t)$.
Note that $\pi(\xit_\pm(t))=\gamma(t)$, hence the
independence of the $\pi$-invariant coordinates,
$y$, $z$, $\tau$ and $\nu$, of the $\pm$ signs.

Notice first that $\tau$ is $\pi$-invariant, so for $t\neq t_0$ we
have
\begin{equation}
d(\tau\circ\gamma)/dt|_{t\pm}=\scHg\tau(\xit_\pm(t))=-2h(\xit_\pm(t))
=2(\tau(\xit_\pm(t))^2-\lambda)=2(\tau(\gamma(t))^2-\lambda)
\end{equation}
where we used that $\tau^2+h=\lambda$ in $\Sigma_{\Delta-\lambda}$.
Thus, $\tau(t)=\tau(\gamma(t))$ is differentiable on
$(t_0-\delta', t_0+\delta')$
except possibly at $t_0$, it is continuous at $t_0$, and its derivative
$\tau'(t)$ extends to a continuous function on $(t_0-\delta', t_0+\delta')$.
Hence $\tau(t)$ is differentiable at $t_0$ and
$\tau'(t_0)=2(\tau(t_0)^2-\lambda)=\scHg\tau(\xit_0)$ for {\em any}
$\xit_0\in\Sigma_{\Delta-\lambda}$. Notice also that, with the
notation of \eqref{eq:def-W_0} in the proof of
Proposition~\ref{prop:normal-prop}, $\tau'(t_0)=W_0\tau=(\scHg\tau)(\xit_0)$.
In particular,
\begin{equation}\label{eq:app-7}
|\tau(t)-\tau_0|\leq C_1|t-t_0|\Mif|t-t_0|<\delta'.
\end{equation}
In fact, the ODE $\tau'(t)=2(\tau(t)^2-\lambda)$, satisfied for $|t-t_0|<
\delta'$, has a unique $\Cinf$ solution, so on $(t_0-\delta',t_0+\delta')$,
$\tau(t)$ is $\Cinf$ and
\begin{equation}
|\tau(t)-(\tau_0+(W_0\tau)(t-t_0))|\leq C|t-t_0|^2.
\end{equation}

From now on we only consider differentiability issues from the left at $t_0$;
of course, the situation on the right is similar.
We define the $\pi$-invariant functions
$\eta=y\cdot\mu$, $\omega_0$, $\omega$ and $\phi=\phi^{(\ep,\delta)}$ as in
the proof of Proposition~\ref{prop:normal-prop}. It is shown there
that there exist $C_0>0$ and $\delta_0>0$ such that if $\ep\in(0,1)$,
$\delta\in(0,\delta_0)$, $\delta<C_0\ep^2$
and $\xit=(y,z,\tau,\mu,\nu)\in\Sigma_{\Delta-\lambda}$
satisfies $\tau_0-\tau\geq-2\delta$ and $\phi(\xit)\leq 2\delta$ then
$\scHg\phi\geq c_0>0$.
So suppose that we fixed some
\begin{equation}
0<T<\min(\delta',C_1\delta_0)
\end{equation}
and let
\begin{equation}
\delta=C_1 T,\ \ep=2(\delta/C_0)^{1/2}.
\end{equation}
Thus, for $t\in[t_0-T,t_0)$, $|\tau(t)-\tau_0|<2\delta$.
As $\phi$ is a $\pi$-invariant function which vanishes at $\xi_0$,
we see that $F=\phi_\pi\circ\gamma$
satisfies $F(t)<0$ and
$dF/dt|_{t\pm}=\scHg\phi(\xit_\pm(t))\geq c_0>0$ for $t\in[t_0-T,t_0)$
(cf.\ the proof of
Proposition~\ref{prop:tgt-bichar-tot-geod} after \eqref{eq:tot-geod-36}).
Taking into
account the form of $\phi$ and \eqref{eq:app-7}, we deduce that
for $t\in[t_0-T,t_0)$, $\omega(t)=\omega(\gamma(t))$ satisfies
\begin{equation}
\omega(t)\leq C_1\ep^4\delta^3|t-t_0|.
\end{equation}
Applying this with $t=t_0-T$ we see that
\begin{equation}
\omega(t_0-T)\leq C_2 T^6.
\end{equation}
Since $\omega$ is independent of $\ep$ and $\delta$, we have deduced that
there exists $\delta_1>0$ such that
\begin{equation}\label{eq:app-20}
t_0-\delta_1<t<t_0\Rightarrow \omega(t)\leq C|t-t_0|^6.
\end{equation}
In particular, under the same assumption,
\begin{equation}\label{eq:omega_0-17}
\omega_0(t)\leq C'|t-t_0|^3,
\end{equation}
so
\begin{equation}
|\tau(t)-(\tau_0+(\frac{W_0 \tau}{W_0\eta})\eta(t))|\leq C''|t-t_0|^{3/2}.
\end{equation}
Since $W_0 \tau\neq 0$ and $\tau(t)$ is $\Cinf$, this shows that
$\eta(t)$ is differentiable from the left at $t_0$ and
\begin{equation}\label{eq:eta(t)-diff}
|\eta(t)-(W_0\eta)(t-t_0)|\leq C|t-t_0|^{3/2},\quad W_0\eta=\scHg\eta(\xit),
\ \xit\in\pih^{-1}(\xi_0)\ \text{arbitrary}.
\end{equation}
Using this and the definition of $\omega_0$ we also conclude that
\begin{equation}
|z_j(t)-(W_0 z_j)(t-t_0)|\leq C|t-t_0|^{3/2},
\end{equation}
\begin{equation}
|\nu_j(t)-(W_0\nu_j)(t-t_0)|\leq C|t-t_0|^{3/2}.
\end{equation}
This proves the proposition for the $\pi$-invariant functions $\tau$, $z_j$,
$\nu_j$ and $\eta$, and indeed it provides a better error estimate.
However, we still need to estimate $y_j$.

To do so, we consider the second term in $\omega$, see \eqref{eq:n-pr-8}.
Thus, from \eqref{eq:app-20},
\begin{equation}\label{eq:app-29}
||y(t)|^2-\mu_0^{-2}\eta(t)^2|\leq C|t-t_0|^3,\quad \mu_0=(\lambda
-\tau_0^2-\htil(z_0,\nu_0))^{1/2}.
\end{equation}
Taking into account \eqref{eq:eta(t)-diff}, we deduce that
\begin{equation}
r(t)=|y(t)|
\end{equation}
satisfies
\begin{equation}
|r(t)^2-4\mu_0^2(t-t_0)^2|\leq C|t-t_0|^{5/2}.
\end{equation}
Thus,
\begin{equation}\label{eq:r(t)-diff}
|r(t)+2\mu_0(t-t_0)|\leq C|t-t_0|^{3/2}.
\end{equation}
Hence, $r(t)$ is also differentiable from the left at $t_0$, and in
particular
\begin{equation}\label{eq:y-23}
|y(t)|=r(t)\leq C|t-t_0|.
\end{equation}

Now,
\begin{equation}\begin{split}
|y(t)-\frac{\eta(t)}{\mu_0^2}\mu_\pm(t)|^2&=|y(t)|^2-\frac{\eta(t)^2}{\mu_0^2}
-\eta(t)^2\,\frac{\mu_0^2-|\mu_\pm(t)|^2}{\mu_0^2}.
\end{split}\end{equation}
By \eqref{eq:n-prop-30}, \eqref{eq:omega_0-17} and \eqref{eq:y-23},
\begin{equation}
||\mu_\pm(t)|^2-\mu_0^2|\leq C(|y(t)|+\omega_0(t)^{1/2})\leq C'|t-t_0|.
\end{equation}
Thus, by \eqref{eq:app-29},
\begin{equation}
|y(t)-\frac{\eta(t)}{\mu_0^2}\mu_\pm(t)|^2\leq C|t-t_0|^3
\end{equation}
In particular, for each $j$ we have
\begin{equation}\label{eq:app-24}
|y_j(t)-\frac{\eta(t)}{\mu_0^2}\mu_{j,\pm}(t)|^2\leq C|t-t_0|^3
\end{equation}

Let
\begin{equation}
\theta_j=y_j/r,
\end{equation}
so $\theta_j$ is a $\pi$-invariant function
away from $C_a$, and we have $|\theta_j|\leq 1$. Also let
\begin{equation}
\theta_j(t)=\frac{y_j(t)}{r(t)},\ t_0-\delta_1<t<t_0.
\end{equation}
By the inductive hypothesis, $\theta_j(t)$ is differentiable for
$t\in(t_0-\delta_1,t_0)$ from both the left and the right and
\begin{equation}
\frac{d\theta_j}{dt}|_{t\pm}
=r(t)^{-1}\frac{dy_j}{dt}-y_j(t)r(t)^{-2}\frac{dr}{dt}
\end{equation}
with
\begin{equation}
dy_j/dt|_{t\pm}=2\mu_{j,\pm}(t)
\end{equation}
and
\begin{equation}
dr/dt|_{t\pm}=
\half |y(t)|^{-1}(d|y|^2/dt|_{t_\pm})=2\frac{\eta(t)}{r(t)}.
\end{equation}
Thus,
\begin{equation}
\frac{d\theta_j}{dt}|_{t\pm}=2r(t)^{-1}(\mu_{j,\pm}(t)-\frac{y_j(t)
\eta(t)}{r(t)^{-2}}),
\end{equation}
so by \eqref{eq:app-24} and \eqref{eq:eta(t)-diff},
\begin{equation}
|\frac{d\theta_j}{dt}|_{t\pm}-2r(t)^{-1}y_j(t)(\mu_0^2\eta(t)^{-1}-\eta(t)
r(t)^{-2}|\leq C|t-t_0|^{-1/2}.
\end{equation}
But, by \eqref{eq:r(t)-diff} and \eqref{eq:eta(t)-diff}, this gives
\begin{equation}
|\frac{d\theta_j}{dt}|_{t\pm}|\leq C|t-t_0|^{-1/2}.
\end{equation}
Integrating from $t_0-\delta_1$ to $t_0$ gives that
$\theta_{j,-}(t_0)=\lim_{t\to t_0-}\theta_j(t)$ exists and
\begin{equation}
|\theta_{j,-}(t_0)-\theta_j(t)|\leq C'|t-t_0|^{1/2}.
\end{equation}
Returning to the original notation, $\theta_j=y_j/r$, we see that
\begin{equation}\label{eq:y(t)-diff}
|y_j(t)+2\mu_0\theta_{j,-}(t_0)(t-t_0)|\leq C'|t-t_0|^{3/2},
\end{equation}
so $y_j(t)$ is differentiable at $t_0$ from the left. We then let
\begin{equation}
\xit_{-}(t_0)=(0,z(t_0),\tau(t_0),\nu(t_0),-\mu_0\theta_{j,-}(t_0)).
\end{equation}
Then the compositions of the $\pi$-invariant coordinate functions
$y_j$, $z_j$, $\tau$ and $\nu_j$ with $\gamma$ are all differentiable from
the left at $t_0$ and the derivative is given by $\scHg$ applied to the
appropriate coordinate function, evaluated at $\xit_-(0)$. Note also
that from \eqref{eq:app-24} and \eqref{eq:y(t)-diff} we have
\begin{equation}
|\mu_\pm(t)-\mu_-(t_0)|\leq C|t-t_0|^{1/2},\quad t\in(t_0-\delta_1,t_0).
\end{equation}
Since
a general smooth $\pi$-invariant function $f$ has the form
\begin{equation}
f(y,z,\tau,\mu,\nu)=f_0(z,\tau,\nu)+\sum y_j f_j(z,\tau,\mu,\nu)+\sum
y_j y_k f_{jk}(y,z,\tau,\mu,\nu),
\end{equation}
$f_0$, $f_j$, $f_{jk}$ all $\Cinf$,
this finishes the proof of the proposition.

\bibliographystyle{plain}
\bibliography{sm}

\end{document}